\documentclass[phd,tocprelim]{cornell}

\usepackage{graphicx}
\usepackage{amsmath}
\usepackage{amsthm}
\usepackage{amscd}
\usepackage{epstopdf}
\usepackage{amssymb}
\theoremstyle{plain}

\newtheorem{thm}{Theorem}
\newtheorem{prop}[thm]{Proposition}
\newtheorem{cor}[thm]{Corollary}
\newtheorem{lem}[thm]{Lemma}
\newtheorem{conv}[thm]{Convention}

\theoremstyle{definition}
\newtheorem{df}[thm]{Definition}
\newtheorem{rem}[thm]{Remark}
\newtheorem{eg}[thm]{Example}
\newtheorem{calc}[thm]{Calculation}
\newtheorem{conj}[thm]{Conjecture}

\usepackage{hangcaption}
\renewcommand{\caption}[1]{\singlespacing\hangcaption{#1}\normalspacing}

\title {Eternal solutions and heteroclinic orbits of a semilinear
parabolic equation}
\author {Michael Robinson} 
\conferraldate{May}{2008} 
\degreefield {Ph.D.}  
\copyrightholder{Michael Robinson}
\copyrightyear{2008}

\begin{document}

\maketitle
\makecopyright

\begin{abstract}
This dissertation describes the space of heteroclinic orbits for a
class of semilinear parabolic equations, focusing primarily on the
case where the nonlinearity is a second degree polynomial with
variable coefficients.  Along the way, a new and elementary proof of
existence and uniqueness of solutions is given.  Heteroclinic orbits
are shown to be characterized by a particular functional being finite.
A novel asymptotic-numeric matching scheme is used to uncover delicate
bifurcation behavior in the equilibria.  The exact nature of this
bifurcation behavior leads to a demonstration that the equilibria are
degenerate critical points in the sense of Morse.  Finally, the space
of heteroclinic orbits is shown to have a cell complex structure,
which is finite dimensional when the number of equilibria is finite.
\end{abstract}

\begin{biosketch}
Michael Robinson entered the field of mathematics by a rather
circuitous route.  In high school, he learned computer programming as
a hobby.  It was there that his first exposure to partial differential
equations occurred, when he wrote a solver for the invicid
Navier-Stokes equations in the plane under the direction of Robert
Ryder (who was then with Pratt \& Whitney).  Feeling that he ought to
understand computer hardware at a deeper level, he enrolled at
Rensselaer Polytechnic Institute and completed a Bachelor of Science
in Electrical Engineering.  His Master of Science degree in
Mathematics at Rensselaer Polytechnic Institute was completed under
the direction of Dr. Ashwani Kapila.  His Master's thesis combined the
Navier-Stokes equations and Maxwell's equations to examine the
scattering of plasma waves.  After completing his Master's degree,
Robinson went to work for Syracuse Research Corporation and enrolled
at Cornell University one year later.
\end{biosketch}

\begin{dedication}
This project is dedicated to my wife, Donna, who urged me to do the
obvious thing and continue my graduate studies.
\end{dedication}

\begin{acknowledgements}
I would like to thank all of the members of my thesis committee for
their helpful and insightful discussions concerning this project.  I
would especially like to thank my advisor, Dr. John Hubbard for
suggesting that I study $\frac{\partial u}{\partial t} = \Delta u -
u^2 + \phi$ as a ``summer project.''  This problem has led into a
surprising variety of interesting mathematics.
\end{acknowledgements}

\contentspage
\figurelistpage

\normalspacing \setcounter{page}{1} \pagenumbering{arabic}
\pagestyle{cornell} \addtolength{\parskip}{0.5\baselineskip}

\chapter{Introduction}
\label{intro_ch}
\section{The Grand Plan}

This dissertation presents some recent progress towards a rather lofty
(and very difficult) goal.  Specifically, there is great interest in
understanding the topology of solution spaces of systems of
semilinear parabolic equations,
\begin{equation}
\label{far_too_general_eqn}
\begin{cases}
\frac{\partial u_1}{\partial t} = D_1(u_1, u_2, ... u_p) + f_1
(t,x,u_1, u_2, ... u_p) \\
\frac{\partial u_2}{\partial t} = D_2(u_1, u_2, ... u_p) + f_2
(t,x,u_1, u_2, ... u_p) \\
...
\frac{\partial u_p}{\partial t} = D_p(u_1, u_2, ... u_p) + f_p
(t,x,u_1, u_2, ... u_p) \\
\end{cases}
\end{equation}
where $u_i \in C^1(\mathbb{R},C^{0,\alpha}(\mathbb{R}^n))$, $D_i$ are
densely-defined linear elliptic (diffusion) operators, and $f_i$
satisfy reasonable smoothness conditions.  Of course, a major problem
for anyone interested in \eqref{far_too_general_eqn} is existence and
uniqueness of short-time solutions!  Although the existence and
uniqueness in general is daunting, many interesting and important
problems have the form \eqref{far_too_general_eqn}.  (Fortunately, for
some special cases, existence and uniqueness can be proven, as is done
in Chapter \ref{conv_ch}.)  The applications of
\eqref{far_too_general_eqn} are numerous; for instance:
\begin{itemize}
\item The Navier-Stokes equations, which describe fluid flow, can be
  put in the form of \eqref{far_too_general_eqn} \cite{Henry}.
  Understanding the topology of the solution space of the
  Navier-Stokes equations gives insight into the onset of turbulence.
  This has applications to fluid amplifiers (viscous fluid logic
  gates, with no moving parts) \cite{Hobbs_1963}, in which a
  particular geometry and set of boundary conditions allows several
  semistable equilibria.  The orbits which connect these equilibria
  (the {\it heteroclinic orbits} or {\it heteroclines}) can involve
  turbulent flows.  As a result, understanding the topology of the
  space of heteroclines for such a situation might provide insight
  into turbulence phenomena.
\item Many chemical reactions are of this form, in particular those
  describing combustion.  The precise nature of the ignition of a
  flame is encoded in the topology of the solution space.
\item The combination of the Navier-Stokes equations with combustion
  equations can model turbulent combustion phenomena.  Such turbulent,
  reacting flows are important in modeling the inside of internal
  combustion engines.  The ignition of a turbulent combustion depends
  very delicately upon the exact nature of the flow, and a topological
  description for such events is lacking.
\item Related to the Navier-Stokes and chemical reaction equations are
  nonlinear wave equations.  Many nonlinear wave equations are of the
  form \eqref{far_too_general_eqn}, and traveling waves appear as
  heteroclines connecting two equilibrium states.  Traveling waves are
  often stable, in that perturbations of them tend to ``wash out''
  over time.  However, there are interesting situations where
  traveling waves are suddenly suppressed as a parameter is changed
  slightly.  This is a consequence of an abrupt change in the topology
  of the solution space.
\item In population biology, \eqref{far_too_general_eqn} describes a
  number of competing or cooperating species.  One can ask about the
  kinds of bifurcations in stable populations when a new species is
  introduced, or when harvesting patterns are changed.
\item One of the most pressing issues in population biology is that of
  non-native invasive species, which disrupt the ecology of many parts
  of our planet.  In agriculture, they cause significant crop losses
  and threaten our protected areas.  One of the important problems
  concerning invasive species is how to displace them with minimal
  ecological impact.  Understanding the topology of the space of
  solutions for \eqref{far_too_general_eqn} would help find optimal
  control algorithms for eliminating (or limiting) the spread of
  invasive species.  There is a vast literature on this subject, going
  back to Fisher \cite{Fisher_1937}, Kolmogorov, Petrovski, and
  Piskunov \cite{KPP_1937}.
\end{itemize}

\section{Specialization to a scalar gradient equation}
Since the general setting of \eqref{far_too_general_eqn} is much
too difficult to allow any kind of progress at this time, we must
instead consider more specialized situations.  To this end, we restrict
attention to the case of 
\begin{itemize}
\item a scalar equation (p=1 in \eqref{far_too_general_eqn}),
\item in one spatial dimension, where
\item the diffusion operator is the Laplacian, $D_1=\frac{\partial^2}{\partial x^2}$,
  and
\item the reaction term is a polynomial.
\end{itemize}

In other words, consider
\begin{equation}
\label{pde}
\frac{\partial u(t,x)}{\partial t}=\Delta u(t,x) +  \sum_{i=0}^N
a_i(x)u^i(t,x)=\Delta u + P(u),
\end{equation}
where $t \in \mathbb{R}$ and $x \in \mathbb{R}^n$, and the $a_i$ are
bounded and smooth.  In this dissertation, we consider {\it eternal}
solutions, those $u$ satisfying \eqref{pde} that lie in
$C^1(\mathbb{R},C^{0,\alpha}(\mathbb{R}))$ for $0<\alpha<1$.  In
particular, observe that $\Delta : C^{0,\alpha}(\mathbb{R}) \to
C^{0,\alpha}(\mathbb{R})$ is densely defined when $\alpha$ is not an
integer.

This kind of equation provides a simple model for a number of physical
phenomena.  First, choosing the right side to be $\Delta u - u^2 + a_1
u$ results in an equation which can represent a model of the
population of a single species with diffusion and a spatially-varying
carrying capacity, $a_1(x)$.  As a second application, this equation
is a very simple model of combustion.  If $a_1$ is a positive
constant, then the equation supports traveling waves.  Such traveling
waves can model the propagation of a flame through a fuel source.

\section{Discussion of the literature}

Equations of the form \eqref{pde} have been of interest to researchers
for quite some time.  Existence and uniqueness of solutions on short
time intervals (on strips $(-t_0,t_0)\times\mathbb{R}^n$) can been
shown using semigroup methods and are entirely standard
\cite{ZeidlerIIA}. However, there are obstructions to the existence of
eternal solutions.  Aside from the typical loss of regularity due to
solving the backwards heat equation, there is also a blow-up
phenomenon which can spoil existence in the forward-time solution to
\eqref{pde}.  Blow-up phenonmena in the forward time Cauchy problem
(where one does not consider $t<0$) have been studied by a number of
authors \cite{Fujita} \cite{Kobayashi_1977} \cite{Weissler_1981}
\cite{Klainerman_1982} \cite{Brezis_1984} \cite{Zheng_1986}
\cite{Zheng_1995}.  More recently, Zhang {\it et al.}
(\cite{Zhang_2000} \cite{Souplet_2002} \cite{Wrkich_2007}) studied
global existence for the forward Cauchy problem for
\begin{equation*}
\frac{\partial u}{\partial t} = \Delta u + u^p - V(x) u
\end{equation*}
for positive $u,V$.  Du and Ma studied a related problem in
\cite{DuMa2001} under more restricted conditions on the coefficients
but they obtained stronger existence results.  In fact, they found
that all of the solutions which were defined for all $t>0$ tended to
equilibrium solutions.

Eternal solutions to \eqref{pde} are rather rare.  Most works which
describe blow-up make the assumption that the solution is positive.
Unfortunately, blow-up is much more difficult to characterize in the
general situation, and understanding exactly what kind of initial
conditions are responsible for blow-up in the Cauchy problem for
\eqref{pde} is an important part of the question.

The boundary value problem that results from taking $x
\in \Omega \subset \mathbb{R}^n$ for some bounded $\Omega$ (instead of
$x \in \mathbb{R}^n$) has also been discussed extensively in the literature
\cite{Henry} \cite{Jost_2007} \cite{BrunovskyFiedler1989}.  For the
boundary value problem, all bounded forward Cauchy problem solutions
tend to limits as $|t|\to\infty$, and these limits are equilibrium
solutions.

Almost all of the literature (including this dissertation) describing
eternal solutions to \eqref{pde} is restricted to discussing
heteroclines.  For unbounded domains and certain symmetries among the
coefficients $a_i$, one can find traveling waves.  Since the
propagation of waves in nonlinear models is of great interest in
applications, there is much written on the subject.  The general idea
is that one makes a change of variables $(t,x) \mapsto \xi = x-ct$
which reduces \eqref{pde} to an ordinary differential equation.
This ordinary differential equation describes the profile of a
traveling wave.  Powerful topologically-motivated techniques, such as
the Leray-Schauder degree, can be used to prove existence of wave
solutions to \eqref{pde}.  Asymptotic methods can be used to determine
the wave speed $c$, which is often of interest in applications.  See
\cite{Volpert_1994} for a very thorough introduction to the subject of
traveling waves in \eqref{pde}.

\section{A Morse-theoretic approach}
A somewhat less traditional approach to studying \eqref{pde}
exists.  This method attempts to directly compute topological
invariants for the space of heteroclinic orbits $\mathcal{H}$ of \eqref{pde}.  It
makes use of the fact that equation \eqref{pde} defines the flow of the $L^2$ gradient of a
certain {\it action} functional,
\begin{equation}
\label{action_eqn}
A(f)=\int_{\mathbb{R}^n} \frac{1}{2} \|\nabla f\|^2 + \sum_{i=0}^{N}
\frac{a_i(x)}{i+1} f^{i+1}(x) dx.
\end{equation}
It is then evident that along a solution $u(t)$ to \eqref{pde},
$A(u(t))$ is a monotonic function in $t$.  As an immediate consequence,
nonconstant $t$-periodic solutions to \eqref{pde} do not exist.  
This kind of behavior suggests that a Morse-theoretic framework might
be helpful.

Morse theory is concerned with the computation of homotopy or homology
groups of a Riemannian or Hilbert manifold $M$ by ``exploring'' it
with a suitable scalar function $f:M \to \mathbb{R}$.  The function
$f$ is selected to satisfy the Morse-Smale (-Floer) conditions, namely
\begin{itemize}
\item a nondegeneracy condition: if $x\in M$ is a critical point ($df
  (x) = 0$), then the Hessian at $x$ is nonsingular,
\item the {\it Morse index}, which is the number of negative
  eigenvalues of the Hessian is finite for each critical point,
\item stable and unstable manifolds for the gradient flow of $f$ are
  transverse (the Smale condition), and
\item if $M$ is noncompact, there is a compact isolating neighborhood
  for each pair of critical points under the gradient flow \cite{Floer_morse}.
\end{itemize}
The function $f$ can be thought of as a special kind of ``height
function'' on $M$.  One then examines the topology of sets $M^t=\{x\in
M | f(x) \le t\}$, which form a cover for $M$.  It is
straightforward to show that the homotopy class of $M^t$ remains
constant on $t \in (t_0,t_1)$ when there are no critical points in
$f^{-1}((t_0,t_1))$.  The homotopy class of $M^t$ changes abruptly,
however, when $f^{-1}(t)$ contains a critical point.  Morse theory
describes how this homotopy class changes by the attachment of handles
to $M^t$.  A very readable introduction to Morse theory is
\cite{MilnorMorse}.

There is a dual formulation of Morse theory, which uses Witten's
complex to compute homology instead of homotopy.  This approach is
better suited to understanding differential equations, as it focuses
not on level sets, but rather on the flow of 
\begin{equation}
\label{gradient_morse_eqn}
\frac{dx}{dt} = \nabla f(x(t)).
\end{equation}
Using this flow, one constructs a chain complex $(C_*,\partial_*)$ in which the
$C_k$ are free modules generated by the critical points of Morse index
$k$.  The boundary maps $\partial_k$ are then constructed by the
formula
\begin{equation*}
\partial_k (q) = \sum_{p \in C_{k-1}} n(q,p) p,
\end{equation*}
where $n(q,p)$ is the number of heteroclinic orbits of
\eqref{gradient_morse_eqn} which connect $q$ to $p$, counted with
sign.  The surprising thing is that this chain complex computes the
homology of $M$!  A thorough, modern treatment of Morse theory can be
found in \cite{Banyaga}.

\section{Floer homology}

A similar theory can be made to work even if the Morse index of all
critical points is infinite.  Instead of relying on the critical
points to supply an index directly, one constructs a ``relative''
index based on the structure of connecting manifolds.  This theory was
first assembled by Floer for the purpose of understanding the homology
of the space of orbits for an exact symplectomorphism
\cite{Floer_gradient}.  More recently, Ghrist {\it et al.}
\cite{Ghrist_2003} have done work on a similar theory for a certain
evolution of braids.  What is crucial to Floer theories is that the
manifold of heteroclinic orbits which connect a given pair of
equilibria (a connecting manifold) is finite dimensional, and compact modulo time
translation.  Suppose that the connecting manifold for a given pair of
equilibria $x,y$ has dimension $\mu(x,y)$.  One shows that the
following two relations hold for this dimension:
\begin{equation}
\label{additivity_for_dims}
\mu(x,x)=0,\; \mu(x,z)=\mu(x,y)+\mu(y,z),
\end{equation}
where $x,y,z$ are distinct critical points.  This relation allows one to assign
indices $I$ to the critical points such that $\mu(x,y)=I(x)-I(y)$.
Evidently, $I$ is only defined uniquely up to an additive constant.
However, one can use the index $I$ in place of the Morse index to
construct a Witten complex in the usual way.  However, the ambiguity
in the definition of $I$ means that the degrees of the resulting
complex are only defined up to this additive constant.  In this way,
one obtains a kind of ambiguous-degree homology theory, which is
called ``Floer homology.''  Alternatively, one may take the dual
approach, and use the finite dimensional connecting manifolds to
assemble a cell complex structure for the space of heteroclines.  The
attaching maps of this cell complex are evidently related to the
boundary maps of the Witten complex.

\section{How to construct a Floer theory for parabolic equations}
For the case of semilinear parabolic equations, the following must be
established in order to construct a Floer theory:
\begin{enumerate}
\item One must compute the Morse index of each critical point, or more
  properly, show that the Morse index is not well-defined due to
  degeneracy.
\item One must show that connecting manifolds are finite dimensional,
  and that they form a cell-complex structure for the space of
  heteroclinic orbits $\mathcal{H}$.
\item The connecting manifolds must obey the additivity relation
  \eqref{additivity_for_dims}.
\item The space of heteroclinic orbits $\mathcal{H}$ must be compact
  moduluo time translation.
\item One must construct the boundary maps in the Witten complex, and
  verify that they actually form a chain complex that computes the
  homology of $\mathcal{H}$.
\end{enumerate}
This dissertation contains proofs of items 1, 2, and most of 3 for a
special case of \eqref{pde}.  (In the case of a bounded spatial
domain, all but item 5 are standard \cite{JRobinson_2001},
\cite{Jost_2007}.) Rather than working with equation \eqref{pde} in
full generality, the later chapters use the following special case:
\begin{equation}
\label{limited_pde}
\frac{\partial u(t,x)}{\partial t}=\frac{\partial^2 u(t,x)}{\partial x^2} - u^2(t,x)
+ \phi(x),
\end{equation}
where $\phi$ tends to zero as $|x|\to\infty$.  The resulting questions
and techniques have obvious generalizations to \eqref{pde}, though
there are many technical obstacles.  Higher degree polynomial
nonlinearities may of course have more than two roots, which creates
the possibility for more complicated equilibrium structure than we
analyze here.

This simpler model still provides insight into applications, as it is
still a model of the population of a single species, with a
spatially-varying carrying capacity, $\phi$.  Indeed, one easily finds
that under certain conditions the behavior of solutions to
\eqref{limited_pde} is reminiscent of the growth and (admittedly tenuous) control
of invasive species \cite{Blaustein_2001}.  It is the control of
invasive species that is of most interest, and it is also what the
structure of the boundary maps reveals.  In the example given in Chapter
\ref{example_ch}, there is one more stable equilibrium, and several
other less stable ones.  The more stable equilibrium can be thought of
as the situation where an invasive species dominates.  The task, then,
is to try to perturb the system so that it no longer is attracted to
that equilibrium.  An optimal control approach is to perturb the
system so that it barely crosses the boundary of the stable manifold
of the the undesired equilibrium, and thereby the invasive species is
eventually brought under control with minimal disturbance to the rest
of the environment. 

\section{Outline and prerequisite results}

While Chapters \ref{instability_ch} and \ref{unstable_ch} contain the
results of most interest to constructing a Floer theory for
\eqref{limited_pde}, the other chapters provide a number of results
that are prerequisites.  

\subsection{Prerequisites that concern the higher degree case}

Chapters \ref{conv_ch} and \ref{classify_ch} are devoted to the
general equation \eqref{pde}.  The first of these provides a new proof
of short-time existence and uniqueness for \eqref{pde}, and as a
side-effect provides a numerical method for approximating its
solutions.  While existence and uniqueness for \eqref{pde} is standard
\cite{ZeidlerIIA}, the usual proofs are not suited to computation.

The first novel result for \eqref{pde} is obtained in Chapter
\ref{classify_ch}, where a decay condition on the $a_i$ allows us to
classify heteroclinic orbits of \eqref{pde}.  In particular, if an
eternal solution $u$ exists and converges uniformly to equilibria as
$t\to\pm\infty$ if and only if $u$ has finite energy (supremum of the
difference in the action functional \eqref{action_eqn} over all time).
Without the decay condition on the $a_i$, finite energy classifies
those eternal solutions that connect finite action equilibria.

The result in Chapter \ref{classify_ch} is actually quite important
for connecting the Morse-theoretic results with the analysis.  Morse
theory, and in particular Witten's complex, requires the flow to have
a gradient structure.  As a result, the space on which the flow acts
must have an inner product structure, so a natural solution space
would be $C^1(\mathbb{R},L^2(\mathbb{R}))$.  However, the proofs of
the cell-complex structure (Chapter \ref{unstable_ch}) do not work
with timeslices in $L^2(\mathbb{R})$.  In particular, one needs this
space to have a Banach algebra structure (in which the Laplacian is
densely defined), so the H\"older space $C^{0,\alpha}(\mathbb{R})$ for
$0< \alpha<1$ is more natural.  What follows from the results of
Chapter \ref{classify_ch} is that the space of heteroclines
$\mathcal{H}$ lies in the intersection
$C^1(\mathbb{R},C^{0,\alpha}(\mathbb{R})) \cap
C^1(\mathbb{R},L^2(\mathbb{R}))$, so in fact there is no difficulty
(Corollary \ref{heteroclines_in_L1}).

\subsection{Prerequisites that concern the quadratic case}

In the remaining chapters (Chapters \ref{nonauto_ch} through
\ref{example_ch}), only the special case \eqref{limited_pde} is
considered.  This is quite sufficient to obtain interesting results
about the structure of $\mathcal{H}$. 

It is important to understand the collection of equilibrium solutions
for \eqref{limited_pde}, which are global solutions to
\begin{equation}
\label{eq_pde}
0=\Delta u - u^2 + \phi
\end{equation}
on all of $\mathbb{R}$.  Like \eqref{pde}, there are obstructions to
global existence in \eqref{eq_pde} \cite{Veron_1996}. Indeed, there
are fairly few global solutions to \eqref{eq_pde}.  We examine
solutions to this problem under asymptotic decay conditions for $\phi$
in Chapter \ref{nonauto_ch}.  The solution reveals delicate
bifurcation behavior in the number of equilbria as $\phi$ is varied.
Further, the asymptotic behavior is such that all global solutions to
\eqref{eq_pde} have finite action (see \eqref{action_eqn}).

Since they are rare, it is reassuring to construct an example of
heteroclinic orbits, which is done in Chapter \ref{global_ch}.  This
example makes specific use of the structure of the equilibrium
solutions, in particular, their asymptotic decay is crucial.

In order to construct a Morse theory for \eqref{pde}, understanding
the dimension of the stable, center, and unstable manifolds of
equilibria is important.  In Chapter \ref{instability_ch} it is shown
that the the center/stable manifold's dimension is typically infinite,
and later in Chapter \ref{unstable_ch} it is shown that the unstable
manifold has finite dimension.  Each equilibrium solution is in fact
unstable, even if its linearization is stable.  This implies that each
equilibrium is a degenerate critical point.  This neatly derails any
hope of using a standard Morse theory, or even using any of its
extensions to infinite dimensional dynamical systems
\cite{Palais_1963}.  (In Chapter \ref{conj_ch}, Conjecture
\ref{degeneracy_is_ok} suggests that restriction of the flow to
$\mathcal{H}$ may correct the degeneracy.)

The most important result of this work is obtained in Chapter
\ref{unstable_ch}, where the space of heteroclinic orbits is shown to
have a cell-complex structure (with finite dimensional cells).  The
dimension of each cell is determined, under a standing assumption of
transversality (Conjecture \ref{baire_conj}).  From the formula for
the dimension of the cells, it is clear that an additivity rule like
\eqref{additivity_for_dims} will hold.  This result is further
explained by an example in Chapter \ref{example_ch}.  Finally, in
Chapter \ref{conj_ch}, several important future directions are
outlined.

\chapter{Short-time existence and uniqueness}
\label{conv_ch}
\section{Introduction}

(This chapter has already been published as \cite{RobinsonIMEX}.)

Existence and uniqueness of solutions for \eqref{pde2} under
reasonable initial conditions have been known for some time.  For
instance, \cite{Henry} and \cite{ZeidlerIIA} contain straightforward
proofs using semigroup methods.  The purpose of this chapter is to
show how a {\it more elementary} proof can be obtained from a sequence
of explicitly computed discrete-time approximations.

Due to their theoretical and computational stability, implicit
iteration schemes are often prefered over their easier-to-implement
explicit analogues.  However, in the case of semilinear equations, one
can form a hybrid implicit-explicit (IMEX) method which offers
computational and theoretical benefits.  The use of IMEX methods for
approximating semilinear parabolic equations is well-established
\cite{AscherRuuthWetton}.  Many of the recent works on these methods
employ discretizations in both space and time.  These fully discrete
approximations can be computed directly by a computer.  However, one
can obtain a stronger condition for convergence of the approximation
if only the time dimension is discretized \cite{Crouzeix}.  We show
how an even stronger condition for convergence is met by the Cauchy
problem for
\begin{equation}
\label{pde2}
\frac{\partial u(x,t)}{\partial t}=\Delta u(x,t) + \sum_{i=0}^\infty
a_i(x) u^i(x,t),
\end{equation}
where $a_i \in L^1(\mathbb{R}^n) \cap L^\infty(\mathbb{R}^n)$, and
how convergence of this method provides an elementary proof of
existence and uniqueness of solutions.  

The Cauchy problem for \eqref{pde2} arises in a variety of settings.
Notably, some reaction-diffusion equations are of this
form \cite{FiedlerScheel}.  Another application is the special case
\begin{equation*}
\frac{\partial u(x,t)}{\partial t}=\Delta u(x,t) - u^2(x,t) + a_0(x),
\end{equation*}
where $a_0$ is a nonzero function of $x$.  This situation corresponds
to a spatially-dependent logistic equation with a diffusion term,
which can be thought of as a toy model of population growth with
migration.

Following \cite{Crouzeix}, the approximation to be used is
\begin{equation}
\label{e_i_deriv}
u_{n+1}=(I-h \Delta)^{-1}(u_n + h \sum_{i=0}^\infty a_i u_n^i),
\end{equation}
which is obtained by inverting the linear portion of a discrete
version of \eqref{pde2}.  For brevity, we shall call \eqref{e_i_deriv}
{\it the} implicit-explicit method.  (In the summary paper
\cite{AscherRuuthWetton}, this is called an SBDF method, to
distinguish it from other implicit-explicit methods.)  One can compute
the operator $(I-h\Delta)^{-1}$ explicitly using Fourier transform
methods, and obtain a proof of the numerical stability of the
iteration as a whole.

\section{A version of the fundamental inequality}

In order to simplify the algebraic expressions, we make the following definitions.

\begin{df}
Let 
\begin{equation}
\label{def_F}
F(u(x,t))=\Delta u(x,t) + \sum_{i=0}^\infty a_i(x) u^i(x,t),
\end{equation}
and
\begin{equation}
\label{def_G}
G(u(x,t))=\sum_{i=0}^\infty a_i(x) u^i(x,t).
\end{equation}
\end{df}
\begin{df}
\label{def_g}
Define the analytic functions
\begin{equation}
\label{def_g_1}
g_1(z)=\sum_{i=0}^\infty \|a_i\|_1 z^i,
\end{equation}
and
\begin{equation}
\label{def_g_infty}
g_\infty(z)=\sum_{i=0}^\infty \|a_i\|_\infty z^i.
\end{equation}
\end{df}

Since we do not discretize the spatial dimension, we can employ some
of the theory of ordinary differential equations.  We therefore first
prove a variant of the fundamental (Gronwall) inequality for
\eqref{pde2} as is done in \cite{HubbardWest}.  The fundamental
inequality gives a sufficient condition for approximate solutions to
converge.  A slightly weaker version of Lemma \ref{excond_lem} was
obtained in Theorem 3.1 of \cite{Crouzeix}, where the existence of
solutions was required.

\begin{lem}
\label{excond_lem}
Suppose $\{u_i\}_{i=1}^\infty$ is a sequence of piecewise $C^1$
functions $u_i:[0,T] \rightarrow C^2(\mathbb{R}^n)\cap
L^1(\mathbb{R}^n) \cap L^\infty(\mathbb{R}^n)$, such that 
\begin{enumerate}
\item there exist $A,B>0$ so that for each $i$ and $t \in [0,T]$,
  $\|u_i(t)\|_1 \le A$ and $\|u_i(t)\|_\infty \le B$,
\item for each $i$ and $t\in [0,T]$, the series
  $g_1(\|u_i(t)\|_1)$ and $g_\infty(\|u_i(t)\|_\infty)$ converge,
\item for each $t \in [0,T]$,
  $\|\frac{d}{dt} u_i(t) - F(u_i(t)) \|_\infty < \epsilon_i$ and
  $\lim_{i\rightarrow\infty}\epsilon_i = 0$, and
\item $u_1(0) = u_i(0)$ for all $i \ge 0$
\end{enumerate}
Then for each $t \in [0,T]$, $\{u_i(t)\}_{i=1}^\infty$ is a Cauchy
sequence in $L^2(\mathbb{R}^n)$.
\begin{proof}
Let $i,j>0$ be given.  Let $\eta(t)=\|u_i(t) -
u_j(t)\|_2^2=\int{(u_i(t)-u_j(t))^2 dx}$.  Notice that the fourth
condition in the hypothesis gives $\eta(0)=0$.  
\begin{equation*}
\eta'(t)=2\int{\left(u_i'(t)-u_j'(t)\right)\left(u_i(t)-u_j(t)\right) dx}.
\end{equation*}
But, $\|\frac{d}{dt} u_i(t) - F(u_i(t)) \|_\infty < \epsilon_i$ is
equivalent to the statement that for each $t \in [0,T]$ and $x \in
\mathbb{R}^n$,  
\begin{equation*}
F(u_i(x,t)) - \epsilon_i < u_i'(x,t) < F(u_i(x,t)) + \epsilon_i,
\end{equation*}
giving
\begin{eqnarray*}
\eta'(t) &\le& 2 \int{\left(F(u_i(t))-F(u_j(t)) \right)(u_i(t)-u_j(t))
  dx} \\&&+ 2(\epsilon_i + \epsilon_j)\int{|u_i(t)-u_j(t)| dx}\\
&\le&2 \int{\left(\Delta u_i(t) + G(u_i(t)) - \Delta u_j(t) -
  G(u_j(t))\right)(u_i(t) - u_j(t)) dx} \\
 &&+ 2(\epsilon_i + \epsilon_j)\|u_i(t)-u_j(t)\|_1 \\
&\le&2 \int{\left(\Delta(u_i(t)-u_j(t))\right)(u_i(t)-u_j(t)) dx} \\
 &&+ 2 \int{\left(G(u_i(t)) - G(u_j(t))\right)(u_i(t)-u_j(t)) dx} 
 + 2(\epsilon_i + \epsilon_j)\|u_i(t)-u_j(t)\|_1 \\
&\le& -2 \int{\|\nabla(u_i(t)-u_j(t))\|^2 dx} + 2
  \|G(u_i(t)) - G(u_j(t))\|_2 \|u_i(t)-u_j(t)\|_2 \\
 &&+ 2(\epsilon_i + \epsilon_j)\|u_i(t)-u_j(t)\|_1 \\
&\le& 2 \|G(u_i(t)) - G(u_j(t))\|_2 \|u_i(t)-u_j(t)\|_2 +
 2(\epsilon_i + \epsilon_j)\|u_i(t)-u_j(t)\|_1. \\
\end{eqnarray*}
Now also 
\begin{eqnarray*}
\|G(u_i(t)) &-& G(u_j(t))\|_2 = \left \|\sum_{k=0}^\infty a_k
(u_i^k(t)-u_j^k(t)) \right \|_2\\
&\le&\sum_{k=0}^\infty \|a_k\|_\infty \left \|u_i^k(t)-u_j^k(t) \right
\|_2 \\
&\le&\sum_{k=0}^\infty \|a_k\|_\infty
\sqrt{\int{\left(u_i^k(x,t)-u_j^k(x,t)\right)^2 dx}} \\
&\le&\sum_{k=0}^\infty \|a_k\|_\infty
\sqrt{\int{\left(u_i(x,t)-u_j(x,t)\right)^2 \left(\sum_{m=0}^{k-1} u_i^m(x,t)
    u_j^{k-m-1}(x,t)\right)^2 dx}} \\
&\le&\sum_{k=0}^\infty \|a_k\|_\infty \left \| \sum_{m=0}^{k-1} u_i^m(t)
    u_j^{k-m-1}(t) \right \|_\infty \|u_i(t)-u_j(t)\|_2 \\
&\le&\left( \sum_{k=0}^\infty \|a_k\|_\infty k B^{k-1}
    \right)\|u_i(t)-u_j(t)\|_2 \\
&\le&g_\infty'(B) \|u_i(t)-u_j(t)\|_2, \\
\end{eqnarray*}
which allows
\begin{eqnarray*}
\eta'(t) &\le& 2 g_\infty'(B) \|u_i(t)-u_j(t)\|_2^2 +  2(\epsilon_i +
\epsilon_j)\|u_i(t)-u_j(t)\|_1. \\
&\le&2 g_\infty'(B) \eta(t) +  2(\epsilon_i +
\epsilon_j)\|u_i(t)-u_j(t)\|_1. \\
\end{eqnarray*}
\begin{eqnarray*}
\eta'(t)-2 g_\infty'(B) \eta(t) &\le& 2(\epsilon_i +
\epsilon_j)\|u_i(t)-u_j(t)\|_1 \\
\frac{d}{dt}\left(\eta(t) e^{-2 g_\infty'(B) t } \right) 
&\le& 2(\epsilon_i +\epsilon_j) e^{-2 g_\infty'(B) t}\|u_i(t)-u_j(t)\|_1, \\
\end{eqnarray*}
so (recall $\eta(0)=0$)
\begin{eqnarray*}
\eta(t) &\le& \left[ 2(\epsilon_i +\epsilon_j) \int_0^t e^{-2
    g_\infty'(B) s } \|u_i(s)-u_j(s)\|_1 ds
    \right]e^{2 g_\infty'(B) t} \\
&\le& \left[ 2(\epsilon_i +\epsilon_j) \int_0^t{\|u_i(s)-u_j(s)\|_1
    ds} \right] e^{2 g_\infty'(B) t}.\\ 
&\le& 4(\epsilon_i +\epsilon_j) A t e^{2 g_\infty'(B) t}.\\
\end{eqnarray*}
Hence as $i,j \rightarrow \infty$, $\eta(t) \rightarrow 0$ for each
$t$.  Thus for each $t$, $\{u_i(t)\}_{i=1}^\infty$ is a Cauchy
sequence in $L^2(\mathbb{R}^n)$.
\end{proof}
\end{lem}

\begin{rem}
Since $C^2(\mathbb{R}^n) \cap L^1(\mathbb{R}^n) \cap L^\infty(\mathbb{R}^n)
\subseteq L^2(\mathbb{R}^n) $ and $L^2$ is complete, Lemma \ref{excond_lem}
gives conditions for existence and uniqueness of a short-time solution
to \eqref{pde2}.
\end{rem}

\begin{lem}
\label{derivative_lem}
Suppose $\{u_i(t)\}_{i=1}^\infty$ is the sequence of functions defined
in Lemma \ref{excond_lem}, and that $u(t)=\lim_{i\rightarrow\infty}
u_i(t)$ in $L^2(\mathbb{R}^n)$.  Then
\begin{equation}
\label{derivativelimit_eqn}
u'(t,x)=\lim_{i\rightarrow\infty}u_i'(t,x) \text{ for almost every $x$},
\end{equation}
wherever the limit exists.
\begin{proof}
Notice that since each $u_i(t) \in L^\infty(\mathbb{R}^n)$ and
$\|u_i(t)\|_\infty \le B$, the dominated convergence theorem allows for
each $x \in \mathbb{R}^n$
\begin{eqnarray*}
\int_0^t{\lim_{i\rightarrow\infty} u_i'(\tau,x) d\tau} &=&
  \lim_{i\rightarrow\infty} \int_0^t{u_i'(\tau,x) d\tau}\\
&=& \lim_{i\rightarrow\infty}(u_i(t,x) - u_i(0,x))\\
&=& u(t,x)-u(0,x) \text{ for almost every $x$}. 
\end{eqnarray*} 
Hence, by differentiating in $t$, 
\begin{equation*}
u'(t,x)=\lim_{i\rightarrow\infty}u_i'(t,x) \text{ for almost every $x$}.
\end{equation*}
\end{proof}
\end{lem}

\section{The implicit-explicit approximation}

In this section, we consider the case of a 1-dimensional
spatial domain, that is, $x \in \mathbb{R}$.  There is no obstruction
to extending any of these results to higher dimensions, though it
complicates the exposition unnecessarily.

As is usual, the first task is to define the function spaces to be
used.  Initial conditions will be drawn from a subspace of
$L^1(\mathbb{R})\cap L^\infty(\mathbb{R})$, as suggested by Lemma
\ref{excond_lem}, and the first four spatial derivatives will be
prescribed, for use in Lemma \ref{hubbardbnd_lem}.

\begin{df}
\label{spaces_df}
Let
\begin{equation*}
W=L^1 \cap C^4(\mathbb{R}),
\end{equation*}
where we interpret $C^4(\mathbb{R})$ as being the space of bounded functions
with four continuous bounded derivatives.  For the remainder of this
chapter, we consider the case where each of the coefficients $a_i \in
W$.  Then let $X=\{f\in W | g_1(\|f\|_1) < \infty\text{ and }
g_\infty(\|f\|_\infty) < \infty\}$.  We consider the case where the
initial condition is drawn from $X$.
\end{df}

An approximate solution given by the implicit-explicit iteration will
be the piecewise linear interpolation through the iterates computed by
\eqref{e_i_deriv}.  A smoother approximation will prove to be
unnecessary, as will be shown in Lemma \ref{existunique_h_thm}.

\begin{df}
\label{hubbard_df}
Suppose $f_0$ and $h > 0$ are given.  Put 
\begin{equation}
\label{hubbard_df_pts}
f_{n+1}=(I-h \Delta)^{-1} ( f_n + h G(f_n) ).
\end{equation}
The function 
\begin{equation}
\label{hubbard_df_fcn}
u(t)=\left(1-\left(\frac{t}{h}-n(t)\right)\right) f_{n(t)} +
\left(\frac{t}{h} - n(t)\right ) f_{n(t)+1},
\end{equation}
where $n(t)=\lfloor \frac{t}{h} \rfloor$, is called the {\bf
  implicit-explicit iteration of size $h$ beginning at $f_0$}.
\end{df}

\begin{calc}
\label{invert_calc}
We explicitly compute the operator $(I-h\Delta)^{-1}$ using Fourier
transforms.  Suppose
\begin{equation*}
(I-h\Delta)u(x)=u(x)-h\Delta u(x)=f(x).
\end{equation*}
Taking the Fourier transform (with transformed variable $\omega$)
gives
\begin{equation*}
\hat{u}(\omega)+h \omega^2 \hat{u}(\omega)=\hat{f}(\omega),
\end{equation*}
\begin{equation*}
\hat{u}(\omega)=\frac{\hat{f}(\omega)}{1+h\omega^2}.
\end{equation*}
The Fourier inversion theorem yields
\begin{eqnarray*}
u(x)&=&\frac{1}{2 \pi} \int{\frac{e^{i\omega x}}{1+h\omega^2} \int{f(y)
    e^{-i\omega y} dy} d\omega}\\
&=&\int{f(y)\left(\frac{1}{2\pi}\int{\frac{e^{i\omega(x-y)}}{1+h\omega^2}
    d\omega} \right) dy}.\\
\end{eqnarray*}

Using the method of residues, this can be simplified to give
\begin{equation}
\label{explicit_op_eqn}
u(x)=\left((I-h\Delta)^{-1}f\right)(x) = \frac{1}{2\sqrt{h}}\int{ f(y)
  e^{-|y-x|/\sqrt{h}} dy }.
\end{equation}
\end{calc}

\begin{calc}
\label{bound_calc}
Bounds on the $L^1$ and $L^\infty$ operator norms of
$(I-h\Delta)^{-1}$ are now computed.  First, let $f\in
L^\infty(\mathbb{R})$.  Then
\begin{eqnarray*}
|\left((I-h\Delta)^{-1}f\right)(x)| &=& \left|\frac{1}{2\sqrt{h}}\int{ f(y)
  e^{-|y-x|/\sqrt{h}} dy }\right|\\
&\le&\|f\|_\infty \frac{1}{2\sqrt{h}} \int{e^{-|y-x|/\sqrt{h}} dy }\\
&\le&\|f\|_\infty \frac{1}{\sqrt{h}} \int_0^\infty{e^{-s/\sqrt{h}}
  ds }\\
&\le&\|f\|_\infty,
\end{eqnarray*}
so $\|(I-h\Delta)^{-1}\|_\infty \le 1$.

Now, let $f\in L^1(\mathbb{R})$.  So then
\begin{eqnarray*}
\|(I-h\Delta)^{-1}f\|_1 &=& \int_{-\infty}^\infty{ \left|
  \frac{1}{2\sqrt{h}} \int_{-\infty}^\infty{f(y)e^{-|y-x|/\sqrt{h}} dy}
  \right| dx}\\ 
&\le& \frac{1}{2\sqrt{h}} \int_{-\infty}^\infty{
  \int_{-\infty}^\infty{ |f(y)|e^{-|y-x|/\sqrt{h}} dy} dx}\\
&\le& \frac{1}{\sqrt{h}} \int_{-\infty}^\infty{|f(y)|
  \int_{0}^\infty{ e^{-|y-x|/\sqrt{h}} dx} dy}\\
&\le& \int_{-\infty}^\infty{ |f(y)| dy} = \|f\|_1,
\end{eqnarray*}
which means $\|(I-h\Delta)^{-1}\|_1 \le 1$. 
\end{calc}

The third condition of Lemma \ref{excond_lem} is a control on the
slope error of the approximation.  A bound on this error may be
established for the implicit-explicit iteration as follows.

\begin{lem}
\label{hubbardbnd_lem}
Suppose $f_0 \in X$, $h>0$.  Put $f(x,t)=
f_0(x) + t D(x)$, where
\begin{equation*}
D=\frac{(I-h\Delta)^{-1}(f_0+hG(f_0)) - f_0}{h}
\end{equation*}
Then for every $0<t<h$, 
\begin{equation}
\label{hubbardbnd_eqn}
\|f'(t)-F(f(t))\|_\infty = O(h).
\end{equation}
\begin{proof}
Recall every function in $X$ will have bounded
partial derivatives up to fourth order from Definition \ref{spaces_df}.
\begin{eqnarray*}
\|f'(t)-F(f(t))\|_\infty &=& \|D-(\Delta (f_0+tD)+G(f_0+tD))\|_\infty\\
&=&\left\|D-\left(\Delta (f_0+tD) + \sum_{i=0}^\infty a_i (f_0+tD)^i\right)
\right\|_\infty\\
&\le&\left \|D-\Delta f_0 -t \Delta D - \sum_{i=0}^\infty a_i
 \left( \sum_{j=0}^i \binom{i}{j} f_0^j(tD)^{i-j} \right) \right\|_\infty\\
&\le&\left \|D-\Delta f_0 - t \Delta D - \sum_{i=0}^\infty a_i f_0^i
 \right\|_\infty + O(h)\\
&\le&\left \|\frac{(I-h\Delta)^{-1}-I}{h}f_0-\Delta f_0\right.\\
&&\left.+((I-h\Delta)^{-1}-I)G(f_0)\right \|_\infty+O(h)\\
\end{eqnarray*}
Now, using the fact that $(I-h\Delta)^{-1}-I=(I-h\Delta)^{-1}(h\Delta)$,
\begin{eqnarray*}
\|f'(t)-F(f(t))\|_\infty
&\le&\left \|(I-h\Delta)^{-1}\Delta f_0-\Delta
f_0\right.\\
&&\left.+(I-h\Delta)^{-1}(h\Delta)G(f_0)\right \|_\infty+O(h)\\
&\le&\|(I-h\Delta)^{-1}(h\Delta)(\Delta f_0 + G(f_0))\|_\infty+O(h)\\
&\le&h \|(I-h\Delta)^{-1}(\Delta F(f_0))\|_\infty+O(h)\\
&\le&h \|(I-h\Delta)^{-1}\|_\infty\|(\Delta F(f_0))\|_\infty+O(h)=O(h)\\
\end{eqnarray*}
\end{proof}
\end{lem}

\begin{lem}
\label{existunique_h_thm}
Suppose $0<h_i \rightarrow 0$.  Let $u_i$ be the implicit-explicit
iteration of size $h_i$ beginning at $f_0 \in X$ on $t \in [0,T]$.
Then provided there exist $A,B>0$ such that for each $i$ and $t \in
[0,T]$, $\|u_i(t)\|_1 \le A$ and $\|u_i(t)\|_\infty \le B$, then the
sequence $\{u_i(t)\}_{i=1}^\infty$ converges pointwise to a function
in $t$.  The limit function is piecewise differentiable in $t$.
\begin{proof}
Let $u_i$ be the implicit-explicit iteration of size $h_i$.  By Lemma
\ref{hubbardbnd_lem}, the slope error is bounded:
\begin{equation*}
\|u_i'(t)-F(u_i(t))\|_\infty = O(h_i) = \epsilon_i.
\end{equation*}
Notice that $\epsilon_i \rightarrow 0$.  Then, since $X \subset
C^2(\mathbb{R}^n)$, Lemma \ref{excond_lem} applies, giving a pointwise
limit function $u(t)$.  Finally, since the slope error uniformly
vanishes, Lemma \ref{derivative_lem} implies that the solution is
piecewise differentiable.
\end{proof}
\end{lem}

\section{``{\it A priori} estimates'' for the approximate solutions}

Now we demonstrate that the implicit-explicit method converges for all
initial conditions in $X$.  Specifically, for each $f_0 \in X$, there
exist $A,B>0$ such that for each $i$ and $t \in [0,T]$, $\|u_i(t)\|_1
\le A$ and $\|u_i(t)\|_\infty \le B$, given sufficiently small $T$.
We begin by recalling that from Calculation \ref{bound_calc}, the
$L^\infty$-norm of $(I-h\Delta)^{-1}$ is less than one.  This means
that for the implicit-explicit iteration,
\begin{eqnarray*}
\|f_{n+1}\|_\infty &\le& \|f_n+hG(f_n)\|_\infty\\
&\le& \|f_n\|_\infty + h \left \|\sum_{i=0}^\infty a_i f_n^i \right\|_\infty\\
&\le& \|f_n\|_\infty + h \sum_{i=0}^\infty \|a_i\|_\infty
\|f_n^i\|_\infty\\
&\le& \|f_n\|_\infty + h \sum_{i=0}^\infty \|a_i\|_\infty
\|f_n\|_\infty^i\\
&\le& \|f_n\|_\infty + h g_\infty(\|f_n\|_\infty)
\end{eqnarray*}
Hence the norm of each step of the implicit-explicit iteration will be
controlled by the behavior of the recursion
\begin{equation}
\label{rec_eqn}
f_{n+1}=f_n+hg_\infty(f_n),
\end{equation}
for $f_n,h,a>0$.  Since we are only concerned with short-time
existence and uniqueness, we look specifically at $h=T/N$ and $0\le n
\le N$, for fixed $T>0$ and $N \in \mathbb{N}$.  

\begin{rem}
\label{clever_rem}
The recursion defined by \eqref{rec_eqn} is an Euler solver for 
\begin{equation}
\label{rec_diff_eqn}
\frac{dy}{dt}=g_\infty(y), \text{ with }y(0)=f_0.
\end{equation}
This equation is separable, and $g_\infty$ is analytic near $f_0$, so
there exists a unique solution for the initial value problem
\eqref{rec_diff_eqn} for sufficiently short time.
Also, whenever $y(t)>0$
\begin{equation*}
\frac{d^2 y }{dt^2}=g_\infty'(y(t)) > 0,
\end{equation*}
the function $y(t)$ is concave up.  As a result, the exact solution to
\eqref{rec_diff_eqn} provides an upper bound for the recursion
\eqref{rec_eqn}.  More precisely, we have the following result.
\end{rem}

\begin{lem}
\label{clever_bound}
Suppose $y(0)=f_0 > 0$ in \eqref{rec_diff_eqn}.  Let $T>0$ be
given so that $y$ is continuous on $[0,T]$, and let $N\in\mathbb{N}$.
Then for each $0\le n \le N$, $f_n \le y(T)$, where $f_n$ satisfies
\eqref{rec_eqn} with $h=T/N$.
\begin{proof}
Since the right side of \eqref{rec_diff_eqn} is strictly positive, the
maximum of $y$ is attained at $T$ on any interval $[0,T]$ where $y$ is
continuous.  Furthermore, since $y(0)>0$, it follows from Remark
\ref{clever_rem} that $y$ is concave up on all of $[0,T]$.  Therefore,
$y$ is a convex function on $[0,T]$.  Hence Euler's method,
\eqref{rec_eqn}, will always underestimate the true value of $y$.
Another way of stating this is that
\begin{equation*}
f_n \le y(nh) \le y(T).
\end{equation*}
\end{proof}
\end{lem}


Using Lemma \ref{clever_bound}, the growth of iterates to
\eqref{rec_eqn} may be controlled independently of the step size.
This provides a uniform bound on the sequence of implicit-explicit
approximations.

\begin{lem}
\label{inf_bound}
Suppose $0<h_i=T/i$ for $i \in \mathbb{N}$.  Let $u_i$ be the
implicit-explicit iteration of size $h_i$ beginning at $f_0 \in X$ on
$t \in [0,T]$.  Then there exists a $B>0$ such that for each $i$ and
$t \in [0,T]$, we have $\|u_i(t)\|_\infty \le B$ for sufficiently
small $T>0$.
\begin{proof}
Suppose $f_{in}$ is the $n$-th step of the implicit-explicit iteration
of size $h_i$.  If we let $y(0)=\|f_0\|_\infty$, Lemma
\ref{clever_bound} implies that for any $i$ and any $0 \le n \le i$
\begin{equation*}
\|f_{in}\|_\infty \le y(T)
\end{equation*}
for sufficiently small T.  Hence by \eqref{hubbard_df_fcn} and the
triangle inequality, $\|u_i(t)\|_\infty \le B$ for all $i$ and $t \in
[0,T]$.
\end{proof}
\end{lem}

With the bound on the suprema of the approximations, we can obtain a
bound on the 1-norms.  

\begin{lem}
\label{one_bound}
Suppose $0<h_i=T/i$ for $i \in \mathbb{N}$.  Let $u_i$ be the
implicit-explicit iteration of size $h_i$ beginning at $f_0 \in X$ on
$t \in [0,T]$.  Then there exists an $A>0$ such that for each $i$ and
$t \in [0,T]$, we have $\|u_i(t)\|_1 \le A$ for sufficiently small
$T>0$.
\begin{proof}
First, notice that Lemma \ref{inf_bound} implies that there is a $B>0$
such that for each $i$ and $t \in [0,T]$, we have $\|u_i(t)\|_\infty
\le A$ for sufficiently small $T>0$.  Again suppose $f_{in}$ is the
$n$-th step of the implicit-explicit iteration of size $h_i$.  Then we
compute
\begin{eqnarray*}
\|f_{i,n+1}\|_1 &\le& \|f_{in}\|_1 + h_i \|G(f_{in})\|_1\\
&\le&\|f_{in}\|_1 + h_i \sum_{k=0}^\infty \|a_k f_{in}^k\|_1\\
&\le&\|f_{in}\|_1 + h_i \sum_{k=0}^\infty \int{|a_k f_{in}^k| dx}\\
&\le&\|f_{in}\|_1 + h_i \sum_{k=1}^\infty \|f_{in}\|_\infty^{k-1}
\|a_k\|_\infty \|f_{in}\|_1+h_i\|a_0\|_1\\
&\le&\|f_{in}\|_1\left(1 + h_i \sum_{k=1}^\infty \|a_k\|_\infty
B^{k-1}\right) + h_i\|a_0\|_1\\ 
&\le&\|f_{in}\|_1\left(1 + \frac{h_i}{B} g_\infty(B) -
\frac{h_i}{B}\|a_0\|_\infty \right) + h_i\|a_0\|_1\\ 
&\le&\|f_{in}\|_1\left(1 + h_i C \right) + h_i\|a_0\|_1\\ 
\end{eqnarray*}
This recurence leads to
\begin{eqnarray*}
\|f_{in}\|_1 &\le& \|f_0\|_1 ( 1 + h_i C)^n + h_i \|a_0\|_1
\sum_{m=0}^{n-1} {(1+h_i C)^m}\\
&\le& \|f_0\|_1 ( 1 + h_i C)^n + h_i \|a_0\|_1
\frac{(1+h_i C)^n-1}{h_i C}\\
&\le& \left(\|f_0\|_1 + \frac{1}{C}\|a_0\|_1\right) ( 1 + h_i C)^n -
\frac{1}{C}\|a_0\|_1\\ 
&\le& \left(\|f_0\|_1 + \frac{1}{C}\|a_0\|_1\right) \left( 1 +
\frac{CT}{i} \right)^n - \frac{1}{C} \|a_0\|_1\\
&\le& \left (\|f_0\|_1 + \frac{1}{C}\|a_0\|_1\right) \left( 1 +
\frac{CT}{i} \right)^i - \frac{1}{C} \|a_0\|_1\\ 
&\le&\left( \|f_0\|_1 + \frac{1}{C}\|a_0\|_1 \right) e^{CT} -
\frac{1}{C} \|a_0\|_1 = A. 
\end{eqnarray*}
Once again, by referring to \eqref{hubbard_df_fcn} and using the
triangle inequality, it follows that $\|u_i(t)\|_1 \le B$ for all $i$
and $t \in [0,T]$.
\end{proof}
\end{lem}

\begin{thm}
\label{full_existunique_h_thm}
Suppose $0<h_i=T/i$ for $i\in\mathbb{N}$.  Let $u_i$ be the
implicit-explicit iteration of size $h_i$ beginning at $f_0 \in X$ on
$t \in [0,T]$.  Then, for sufficiently small $T>0$, the sequence
$\{u_i(t)\}_{i=1}^\infty$ converges pointwise to a function in $t$.
The limit function is piecewise differentiable in $t$.
\begin{proof}
This compiles the results of Lemma \ref{existunique_h_thm}, Lemma
\ref{inf_bound}, and Lemma \ref{one_bound}.
\end{proof}
\end{thm}

\begin{rem}
These proofs can be generalized further to handle all equations of
the form

\begin{equation*}
\frac{\partial u(t)}{\partial t} = L(u(t)) + G(u),
\end{equation*}
where $G$ is as in \eqref{def_G}.  If the operator $L$ satisfies
\begin{itemize}
\item $L:L^1(\mathbb{R})\cap L^\infty(\mathbb{R}) \cap
  C^\infty(\mathbb{R}) \rightarrow L^\infty(\mathbb{R}) \cap
  C^\infty(\mathbb{R})$ is a sectorial linear operator \cite{Henry},
\item $\|(I-hL)^{-1}\|_1 \le 1$ and $\|(I-hL)^{-1}\|_\infty \le 1$,
\end{itemize}
then the implicit-explicit iteration
\begin{equation*}
f_{n+1} = (I-hL)^{-1}( f_n + h G(f_n) )
\end{equation*}
converges whenever $f \in X$.
\end{rem}

\begin{rem}
Additionally, the techniques can be easily extended to handle the
initial boundary value problem
\begin{equation*}
\frac{\partial u(x,t)}{\partial t}=\Delta u(x,t) + \sum_{i=0}^\infty
a_i(x) u^i(x,t),\text{ for } x \in K \subset \mathbb{R}^n,t>0
\end{equation*}
with $u(x,t)=v(x,t)$ a given Lipschitz function along $\partial K
\times [0,\infty)$, for $K$ compact with smooth boundary.  In this
case, a boundary term appears in the estimate for $\eta'(t)$ in Lemma
\ref{excond_lem}, which depends on the Lipschitz constant of $v$.
Additionally, in Definition \ref{hubbard_df}, one defines
$f_{n+1}$ to be the unique solution to the linear elliptic boundary
value problem
\begin{equation*}
(I-h\Delta)f_{n+1}=f_n+hG(f_n)
\end{equation*}
with $f_{n+1}(x)=v(x,nh)$ for $x \in \partial K$.
\end{rem}

\section{Conclusions}

The convergence proof for the implicit-explicit method presented here
has a number of advantages.  First of all, like all IMEX
methods, each approximation to the solution is computed explicitly.
As a result, a fully discretized version (as is standard in the
literature) is easy to program on a computer.  Theorem \ref{full_existunique_h_thm}
therefore assures the convergence of these fully discrete methods.

However, since the implicit-explicit method presented here is
discretized only in time, the convergence proof actually shows the
existence of a semigroup of solutions.  As a result, the convergence
proof forms a bridge between the functional-analytic viewpoint of
differential equations, namely that of semigroups, and the numerical
methods used to approximate solutions.  While the existence and
uniqueness of solutions for \eqref{pde2} has been known via semigroup
methods, the proof provided here gives a more elementary explanation
of how this occurs.  In particular, it approximates the semigroup
action directly.

\chapter{Classification of heteroclines}
\label{classify_ch}
\section{Introduction}

(This chapter is available on the arXiv as \cite{RobinsonClassify}.)

In this chapter, the global behavior of smooth solutions to the
semilinear parabolic equation \eqref{pde}
\begin{equation}
\label{pde_classify}
\frac{\partial u(x,t)}{\partial t}=\Delta u(x,t) - u^N + \sum_{i=0}^{N-1}
a_i(x) u^i(t,x)=\Delta u + P(u),
\end{equation}
for $(t,x) \in \mathbb{R}\times\mathbb{R}^n = \mathbb{R}^{n+1}$ is
considered, where $N\ge 2$ and $a_i \in L^\infty(\mathbb{R}^n)$ are
smooth with all derivatives of all orders bounded.  

The main result is that solutions to \eqref{pde_classify} which are
heteroclinic orbits connecting two sufficiently regular equilibrium
solutions of \eqref{pde_classify} are characterized by finite energy
(Definition \ref{energy_df}).  That this characterization is necessary
at all comes from the fact that the spatial domain of \eqref{pde_classify} is
unbounded.  For bounded spatial domains, all bounded global solutions
converge to equilibria \cite{Jost_2007}. The strength of our result
comes from the fact that the finite energy constraint makes solutions
behave rather well.  Therefore, this result is much sharper than what
has typically been obtained in the past, and it applies to more
complicated nonlinear terms.  

The disadvantage is that in doing so, we cannot treat some of the more
complicated aspects of the dynamics.  In particular, traveling wave
solutions do not have finite energy.  Even though a traveling wave
will often converge locally to equilibria, at least one of those
equilibria will not be admissible in our analysis.  On the other hand,
we can exclude traveling waves if we require that all the
coefficients $a_i$ decay fast enough and only consider one spatial
dimension.  Then our result establishes an equivalence between the
heteroclinic orbits and the finite energy solutions.

For somewhat more restricted nonlinearities, Du and Ma were able to
use squeezing methods to obtain similar results to what we obtain
here.  In particular, they also show that certain kinds of solutions
approach equilibria \cite{DuMa2001}. In a somewhat different setting,
Floer used a finite energy constraint for solutions and a regularity
constraint on equilibria to characterize heteroclinic orbits of an
elliptic problem \cite{Floer_gradient}. The techniques of Floer were
subsequently used by Salamon to provide a new characterization of
solutions to gradient flows on finite-dimensional
manifolds \cite{Salamon_1990}. In this chapter, we recast some of
Salamon's work into a parabolic setting, and of course work within an
infinite-dimensional space.

\section{Finite energy constraints}

From Chapter \ref{conv_ch}, we have that solutions to \eqref{pde_classify} exist along strips of
the form $(t,x)\in I\times \mathbb{R}^n$ for sufficiently small
$t$-intervals $I$.  One might hope to extend such solutions to all of
$\mathbb{R}^{n+1}$, but for certain choices of initial conditions such
eternal solutions may fail to exist.  Fujita's classic paper
\cite{Fujita} gives examples of this ``blow-up'' pathology.  We will
specifically avoid it by considering only eternal solutions to
\eqref{pde_classify}.  By eternal solutions, we mean those which are defined for
all $\mathbb{R}^{n+1}$, have one continous partial derivative in time,
and two continous partial derivatives in space.  It should be noted
that eternal solutions to \eqref{pde_classify} are quite rare: the
backwards-time Cauchy problem contains a heat operator, and so most
solutions will not extend to all of $\mathbb{R}^{n+1}$.

\begin{df}
\label{action_df}
Our analysis of \eqref{pde_classify} will make considerable use of the fact
that it is a gradient differential equation.  Recall that the right
side of \eqref{pde_classify} is the $L^2(\mathbb{R}^n)$ gradient of the
{\it action functional} \eqref{action_eqn}, defined for all $f \in
C^1(\mathbb{R}^n)$ by
\begin{equation*}
A(f)=\int_{\mathbb{R}^n} \frac{1}{2} \|\nabla f(x)\|^2 - \frac{f^{N+1}(x)}{N+1}+ \sum_{i=0}^{N-1}
\frac{a_i(x)}{i+1} f^{i+1}(x) dx.
\end{equation*}
\end{df}
It is then evident that along a solution $u(t)$ to \eqref{pde_classify},
\begin{eqnarray*}
\frac{dA(u(t))}{dt} &=& dA|_{u(t)}\left(\frac{\partial u}{\partial t}\right)\\
&=& \left < \nabla A(u(t)), \frac{\partial
  u}{\partial t}  \right >\\
&=&\left < \Delta u + P(u), \frac{\partial
  u}{\partial t}  \right >\\
&=&\left \|\frac{\partial
  u}{\partial t}  \right \|_2^2 \ge 0,\\
\end{eqnarray*}
so $A(u(t))$ is a monotone function.  As an immediate consequence,
nonconstant $t$-periodic solutions to \eqref{pde_classify} do not exist.  

\begin{df}
\label{energy_df}
The {\it energy functional} is the following quantity defined on the
space of functions $\mathbb{R}^{n+1}\to \mathbb{R}$ with one
continuous partial derivative in the first variable ($t$), and two
continuous partial derviatives in the rest ($x$):
\begin{equation}
E(u)=\frac{1}{2}\int_{-\infty}^\infty \int \left | \frac{\partial u}{\partial t}
\right |^2 + \left | \Delta u + P(u) \right |^2 dx\, dt.
\end{equation}
\end{df}

\begin{calc}
\label{action_energy_calc}
Suppose $u$ is in the domain of definition for the energy functional, then
\begin{eqnarray*}
E(u)&=&\frac{1}{2}\int_{-\infty}^\infty \int \left | \frac{\partial u}{\partial t}
\right |^2 + \left | \Delta u + P(u) \right |^2 dx\, dt\\
&=&\frac{1}{2}\int_{-\infty}^\infty \int \left ( \frac{\partial u}{\partial t}
 - \Delta u - P(u) \right )^2 + 2\frac{\partial u}{\partial t}\left( \Delta u + P(u) \right) dx\, dt\\
&=&\frac{1}{2}\int_{-\infty}^\infty \int \left ( \frac{\partial u}{\partial t}
 - \Delta u - P(u) \right )^2 dx\, dt + \int_{-\infty}^\infty \left< \frac{\partial
  u}{\partial t}, \Delta u + P(u) \right > dt\\
&=&\frac{1}{2}\int_{-\infty}^\infty \int \left ( \frac{\partial u}{\partial t}
 - \Delta u - P(u) \right )^2 dx\, dt + \int_{-\infty}^\infty \left< \frac{\partial
  u}{\partial t}, \nabla A(u(t))\right > dt\\
&=&\frac{1}{2}\int_{-\infty}^\infty \int \left ( \frac{\partial u}{\partial t}
 - \Delta u - P(u) \right )^2 dx\, dt + \int_{-\infty}^\infty
\frac{d}{dt} A(u(t)) dt\\
&=&\left.\frac{1}{2}\int_{-\infty}^\infty \int \left ( \frac{\partial u}{\partial t}
 - \Delta u - P(u) \right )^2 dx\, dt +  A(u(T)) \right|_{T=-\infty}^\infty.\\
\end{eqnarray*}
This calculation shows that finite energy solutions to \eqref{pde_classify}
minimize the energy functional.  If a solution to \eqref{pde_classify} is a
heteroclinic connection between two equilibria, then the energy
functional measures the difference between the values of the action
functional evaluated at the two equilibria.  The main result of this
chapter is the converse, so that finite energy characterizes the
solutions which connect equilibria.
\end{calc}

\begin{rem}
Finite energy solutions to \eqref{pde_classify} are even more rare than eternal
solutions.  However, the set of finite energy solutions is not entirely
vacuous, as will be shown in Chapter \ref{global_ch}.
\end{rem}

It is well-known that when equations like \eqref{pde_classify} exhibit the
correct symmetry, they can support traveling wave
solutions \cite{FiedlerScheel}. A typical traveling wave solution $u$
has a symmetry like $u(t,x)=U(x-ct)$ for some $c \in \mathbb{R}$.  As
a result, it is immediate that traveling waves will have infinite
energy.  On the other hand, they also evidently connect equilibria.
As a result, Calculation \ref{action_energy_calc} shows that a
necessary condition for traveling waves is that there exists at least
one equilibrium whose action is infinite.  In this chapter, we will
consider only equilibria with finite action, and solutions with finite
energy.  As a result, we will not be working with traveling waves.

\section{Convergence to equilibria}

In this section, we show that finite energy solutions tend to
equilibria as $|t| \to \infty$.  In doing this, we follow Floer in
\cite{Floer_gradient} which leads us through an essentially standard
parabolic bootstrapping argument.

\begin{lem}
\label{polybound_lem}
Let $U \subseteq \mathbb{R}^n$ and $u \in W^{k,p}(U)$ satisfy $\|D^j
u\|_\infty \le C < \infty$ for $0 \le j \le k$ (in particular, $u$ is
bounded).  If $P(u) = \sum_{i=1}^N a_i u^i$ with $a_i \in L^\infty(U)$
then there exists a $C'$ such that $\|P(u)\|_{k,p} \le C' \|u\|_{k,p}$.  
\begin{proof}
First, using the definition of the Sobolev norm,
\begin{equation*}
\|P(u)\|_{k,p}=\sum_{j=0}^k \|D^j P(u) \|_p \le \sum_{j=0}^k
\sum_{i=1}^N \|D^j a_i u^i \|_p.
\end{equation*}
Now $|D^j a_i u^i| \le P_{i,j}(u,Du,...,D^j u)$ is a polynomial in $j$
variables with constant coefficients, which has no constant term.  (It
has constant coefficients because the derivatives of the $a_i$ are
bounded.)  Additionally,
\begin{eqnarray*}
\|(D^m u)^q D^j u\|_p &=& \left ( \int \left | (D^m u)^q D^j u \right
|^p \right)^{1/p}\\
&\le& \|D^m u\|_\infty^q \left ( \int \left | D^j u \right
|^p \right)^{1/p} \le C^q \|D^j u\|_p,\\
\end{eqnarray*}
so by collecting terms,
\begin{equation*}
\|P(u)\|_{k,p} \le  \sum_{j=0}^k
\sum_{i=1}^N \|D^j a_i u^i \|_p \le \sum_{j=0}^k A_j \|D^j u\|_p \le C' \|u\|_{k,p}.
\end{equation*}
\end{proof}
\end{lem}

The following result is a parabolic bootstrapping argument that does
most of the work.  In it, we follow Floer in \cite{Floer_gradient},
replacing ``elliptic'' with ``parabolic'' as necessary.

\begin{lem}
\label{parabolic_bootstrap_lem}
If $u$ is a finite energy solution to \eqref{pde_classify} with $\|D^j
u\|_{L^\infty((-\infty,\infty)\times V)} \le C < \infty$ for $0 \le j
\le k$ with $k\ge 1$ on each compact $V \subset \mathbb{R}^n$, then
each of $\lim_{t \to \pm \infty} u(t,x)$ exists, and converges with
$k$ of its first derivatives uniformly on compact subsets of
$\mathbb{R}^n$.  Further, the limits are equilibrium solutions to
\eqref{pde_classify}.
\begin{proof}
Define $u_m(t,x)=u(t+m,x)$ for $m=0,1,2...$.  Suppose $U \subset
\mathbb{R}^{n+1}$ is a bounded open set and $K \subset U$ is compact.  Let
$\beta$ be a bump function whose support is in $U$ and takes the value
1 on $K$.  We take $p>1$ such that $kp>n+1$.  Then we can consider
$u_m \in W^{k,p}(U)$ (recall that $u$ and its first $k$ derivatives of
$u$ are bounded on the closure of $U$), and we have
\begin{equation*}
\|u_m\|_{W^{k+1,p}(K)} \le \|\beta u_m \|_{W^{k+1,p}(U)}.
\end{equation*}
Then using the standard parabolic regularity for the heat operator,
\begin{equation*}
\|\beta u_m \|_{W^{k+1,p}(U)} \le C_1 \left \| \left (
\frac{\partial}{\partial t} - \Delta\right ) (\beta u_m)\right \|_{W^{k,p}(U)}.
\end{equation*}
Let $P'(u)=-u^N+\sum_{i=1}^{N-1} a_i u^i$, noting carefully that we have left
out the $a_0$ term.  The usual product rule, and a little
work, as suggested in \cite{Salamon_1990} yields the following sequence of
inequalities
\begin{eqnarray*}
\|u_m\|_{W^{k+1,p}(K)}&\le& C_1 \left \|\beta
\left(\frac{\partial}{\partial t} - \Delta \right) u_m \right
\|_{W^{k,p}(U)} + C_2 \|u_m\|_{W^{k,p}(U)}\\
&\le& C_1 \left \|\beta
\left(\frac{\partial}{\partial t} - \Delta \right) u_m + \beta P'(u_m)
- \beta P'(u_m) \right
\|_{W^{k,p}(U)} \\&&+ C_2 \|u_m\|_{W^{k,p}(U)}\\
&\le& C_1  \|\beta a_0\|_{W^{k,p}(U)} + C_1 \| \beta P'(u_m)
\|_{W^{k,p}(U)} + C_2 \|u_m\|_{W^{k,p}(U)}\\
&\le& C_1 \|\beta a_0\|_{W^{k,p}(U)} + C_3 \|u_m\|_{W^{k,p}(U)},\\
\end{eqnarray*}
where the last inequality is a consequence of Lemma
\ref{polybound_lem}.  By the hypotheses on $u$ and $a_0$, this implies
that there is a finite bound on $\|u_m\|_{W^{k+1,p}(K)}$, which is
independent of $m$.  Now by our choice of $p$, the general
Sobolev inequalities imply that $\|u_m\|_{C^{k+1-(n+1)/p}(K)}$ is
uniformly bounded.  By choosing $p$ large enough, there is a
subsequence $\{v_{m'}\} \subset \{u_m\}$ such that $v_{m'}$ and its
first $k$ derivatives converge uniformly on $K$, say to $v$.  For any
$T>0$, we observe
\begin{eqnarray*}
\int_{-T}^T \int \left | \frac{\partial v}{\partial t} \right |^2 dx
\, dt &=& \lim_{{m'} \to \infty}\int_{-T}^T \int \left | \frac{\partial v_{m'}}{\partial t} \right |^2 dx
\, dt\\
&=&\lim_{{m'} \to \infty} \int_{{m'}-T}^{{m'}+T} \int \left | \frac{\partial u}{\partial t} \right |^2 dx
\, dt=0,\\
\end{eqnarray*}
where the last equality is by the finite energy condition.  Hence
$\left|\frac{\partial v}{\partial t}\right| = 0$ almost everywhere,
which implies that $v$ is an equilibrium and that $\lim_{t\to\infty}
u(t,x) = v(x)$.  Similar reasoning works for
$t\to - \infty$.
\end{proof}
\end{lem}

Now we would like to relax the bounds on $u$ and its derivatives, by
showing that they are in fact consequences of the finite energy
condition.

\begin{lem}
\label{finite_energy_consequences_lem}
Suppose that either $n=1$ (one spatial dimension) or $N$ is odd, then
we have the following. If $u$ is a finite energy solution to
\eqref{pde_classify}, then the the limits $\lim_{t\to\pm\infty} u(t,x)$ exist
uniformly on compact subsets, and
additionally,
\begin{itemize}
\item $u$ is bounded,
\item the derivatives $Du$ are bounded,
\item and therefore the limits are continuous equilibrium solutions.
\end{itemize}
\begin{proof}
Note that since 
\begin{equation*}
E(u)=\frac{1}{2}\int_{-\infty}^\infty \int \left | \frac{\partial u}{\partial t}
\right |^2 + \left | \Delta u + P(u) \right |^2 dx\, dt < \infty,
\end{equation*}
we have that for any $\epsilon>0$,
\begin{equation*}
\lim_{T\to\infty}\frac{1}{2}\int_{T-\epsilon}^{T+\epsilon} \int \left | \frac{\partial u}{\partial t}
\right |^2 + \left | \Delta u + P(u) \right |^2 dx\, dt = 0,
\end{equation*}
whence $\lim_{t\to\infty} \left|\frac{\partial u}{\partial t}\right| =
0$ for almost all $x$.  So this gives that the limit is an equilibrium
almost everywhere.  Of course, this argument works for $t\to -\infty$.

Now in the case of $N$ being odd, a comparison principle shows that
solutions to \eqref{pde_classify} are always bounded.  So we need to consider
the case with $N$ even.  In that case, a comparison principle on
\eqref{pde_classify} shows that $u$ is bounded from {\it above}. On the other
hand, if $N$ is even we have assumed that $n=1$ in this case, and it
follows from an easy ODE phase-plane argument that unbounded
equilibria are bounded from {\it below}.  (Here we have used that the
coefficients $a_i$ are bounded.)  As a result, we must conclude that
if a solution to \eqref{pde_classify} tends to any equilibrium, that
equilibrium (and hence $u$ also) must be bounded.

Now observe that $\left|\frac{\partial u}{\partial t}\right|\to 0$ as
$t\to\infty$ on almost all of any compact $K \subset \mathbb{R}^n$,
and that $\left|\frac{\partial u}{\partial t}\right| \le a < \infty$
for some finite $a$ on $\{(t,x)|t=0,x\in K\}$ by the smoothness of
$u$.  By the compactness of $K$, this means that if
$\left\|\frac{\partial u}{\partial
    t}\right\|_{L^\infty((-\infty,\infty)\times K)} = \infty$, there
must be a $(t^*,x^*)$ such that $\lim_{(t,x)\to
  (t^*,x^*)}\left|\frac{\partial u}{\partial t}\right| = \infty$.
This contradicts smoothness of $u$, so we conclude
$\left|\frac{\partial u}{\partial t}\right|$ is bounded on the strip
$(-\infty,\infty)\times K$.  On the other hand, the finite energy
condition also implies that for each $v\in\mathbb{R}^n$,
\begin{equation*}
\lim_{s\to\infty} \int_{-\infty}^\infty \int_{K+sv}
\left|\frac{\partial u}{\partial t} \right |^2 dx\, dt = 0,
\end{equation*}
whence we must conclude that $\lim_{s\to\infty}\left| \frac{\partial
    u(t,x+sv)} {\partial t}\right|=0$ for almost every
$t\in\mathbb{R}$ and $x\in K$.  Thus the smoothness of $u$ implies
that $\left| \frac{\partial u}{\partial t}\right|$ is bounded on all
of $\mathbb{R}^{n+1}$.

Next, note that since $\left|\frac{\partial u}{\partial t}\right|$ and
$u$ are both bounded, then so is $\Delta u$.  (Use the boundedness of
the coefficients of $P$.)  Taken together, this implies that all the
spatial first derivatives of $u$ are also bounded.

As a result, we have on $K$ a bounded equicontinuous family of
functions, so Ascoli's theorem implies that they (after extracting a
suitable subsequence) converge uniformly on compact subsets of $K$ to
a continuous limit.
\end{proof}
\end{lem}

\begin{cor}
\label{limits_to_equilibria}
Suppose that either $n=1$ or $N$ is odd.  An eternal solution $u$
to \eqref{pde_classify} has finite energy if and only if each of the following hold:
\begin{itemize}
\item each of $U_{\pm}(x) =\lim_{t \to
\pm \infty} u(t,x)$ exists and converges with its first derivatives
uniformly on compact subsets of $\mathbb{R}^n$,
\item $U_{\pm}$ are bounded, continuous equilibrium solutions to \eqref{pde_classify},
\item and either $|A(U_+)-A(U_-)|<\infty$ or $U_+=U_-$.
\end{itemize}
\end{cor}

\begin{rem}
If we consider the more limited case of \eqref{limited_pde}, then the
asymptotic decay rate for equilibria indicates that all equilibria
have finite action (Corollary \ref{finite_action_eq_cor}).  In this
case, Corollary \ref{limits_to_equilibria} characterizes {\it all}
heteroclinic orbits, not just those whose action difference is finite.
\end{rem}

\begin{thm}
\label{unif_limits_to_equilibria}
Suppose that $n=1$ and that all equilibria have finite action.  If
$u$ is a finite energy solution then it converges uniformly to
equilibria as $|t|\to\infty$.
\begin{proof}
Suppose that $u$ tends to equilibrium solutions $f_\pm$ as $t\to
\pm\infty$.  Suppose that this convergence is not uniform, so that
there exists an $\epsilon>0$ such for each $T$, there is a $|t|>|T|$
with either $\|u(t)-f_-\|_\infty > \epsilon$ or $\|u(t)-f_+\|_\infty
> \epsilon$.  We therefore postulate the existence of a pair of
sequences $\{t_j\}$, $\{x_j\}$ such that $|t_j| \to \infty$ and
$\min\{|u(t_j,x_j)-f_-(x_j)|,|u(t_j,x_j)-f_+(x_j)|\}>\epsilon$ for all
$j$.  We assume that for each $j$, $x_j$ is chosen so that
$\min\{|u(t_j,x_j)-f_-(x_j)|,|u(t_j,x_j)-f_+(x_j)|\}$ is maximized.
Notice that since $u\to f_\pm$ uniformly on compact subsets, we must
have $|x_j|\to\infty$.

To simplify the discussion, we find an $R>0$ such that for all
$|y|>\epsilon$ and $|x|>R$,
\begin{equation}
\label{controlled_x}
\left|\sum_{i=0}^{N-1} a_i(x) y^i\right| < \frac{1}{2}|y|^N.
\end{equation}
We assume that $|x_j|>R$ for all $j$.  This condition ensures that the
leading nonlinear coefficient of \eqref{pde_classify} dominates.

We discern three cases, which we can consider without loss of
generality after extracting a suitable subsequence of
$\{(t_j,x_j)\}$.  In each of the cases, we shall perform a coordinate
transformation so that the equilibrium to which $u$ converges is the
zero function.  In particular, we start the sums at 1 rather than 0.

\begin{enumerate}
\item Suppose $t_j \to +\infty$ and $u(t_j,x_j)>\epsilon>0$.  Since
  $x_j$ is chosen at a maximum of $u(t_j)$ for each $j$, we have that
  $\frac{\partial^2 u(t_j,x_j)}{x^2}<0$ by the maximum principle.  As
  a result, 
\begin{eqnarray*}
\frac{\partial}{\partial t} u(t_j,x_j) &=& \frac{\partial^2}{\partial
  x^2}u(t_j,x_j) - u^N(t_j,x_j) +
  \sum_{i=1}^{N-1}a_i(x_j)u^i(t_j,x_j)\\
&\le&- u^N(t_j,x_j) +
  \sum_{i=1}^{N-1}a_i(x_j)u^i(t_j,x_j)\\
&\le& -\frac{\epsilon^N}{2}, \text{ by \eqref{controlled_x}}.\\
\end{eqnarray*}
Therefore we conclude that $\|u(t)\|_\infty\to 0$.

\item Suppose $t_j \to -\infty$ and $u(t_j,x_j)>\epsilon>0$.
Consider the time-reversed version of \eqref{pde_classify}, namely
\begin{equation}
\label{backwards_pde}
\frac{\partial u}{\partial t}=-\frac{\partial^2 u}{\partial x^2} + u^N
- \sum_{i=1}^{N-1} a_i u^i.
\end{equation}
The comparison principle works in reverse for this equation!  Suppose
that $v(t,x)=U(t)$ is a spatially constant solution to
\eqref{backwards_pde} with $U(t_j)=v(t_j,x_j)=u(t_j,x_j)>0$ for some
$j$.  Then, shortly thereafter, $\|u(t)\|_\infty > U(t)$, since 
\begin{eqnarray*}
\frac{\partial u(t_j,x_j)}{\partial t}&\ge& U^N(t_j) -
\sum_{i=1}^{N-1} a_i(x_j) U^i(t_j)\\
&\ge& \frac{1}{2} U^N(t_j) > 0.
\end{eqnarray*}
  On the other hand, this rate of growth indicates that $u$ blows up
in finite time.  This contradicts the fact that $u$ is an eternal
solution.

\item Suppose $u(t_j,x_j)<-\epsilon<0$ and that $t_j \to -\infty$ or
  $t_j \to +\infty$.  If $N$ is odd, then this case can be covered by
  the previous ones, {\it mutatis mutandis}.  Therefore, we assume $N$
  is even.  We assume that the limit as $t\to +\infty$ of $u$ is the
  zero function.  From Lemma \ref{finite_energy_consequences_lem}, we have a constant $A=\left
  \|\frac{\partial^2 u}{\partial x^2} \right\|_\infty < \infty$ which is
  independent of $t$.  Thus for each $t_j$, we have an upper bound for
  $u(t_j)$ which looks like
\begin{equation}
\label{u_upper_bound}
U_j(x) = \min\left\{\mu_j, u(t_j,x_j)+\frac{A}{2}(x-x_j)^2 \right\},
\end{equation}
for some $\mu_j>0$.  In particular, note that $\mu_j \to 0$ by the
previous cases.  

  We show that the forward Cauchy problem for \eqref{pde_classify} started with
  $U_j$ as an initial condition blows up for sufficiently large $j$.
  By the comparison principle, this implies that $u$ cannot be an
  eternal solution, which is a contradiction.  This can be shown using
  the method of Fujita, which we briefly sketch here.  

Apply the coordinate transformation $w=u-\mu$ for some $\mu >
  \mu_j>0$.  Therefore, the initial condition can be made entirely
  negative, and by the previous cases, the solution {\it stays}
  negative for arbitrarily long future time (by taking $j$ large).
  (Notice that it may not remain negative for {\it all} future time in
  the case where $t_j\to -\infty$.)  This transformation changes
  \eqref{pde_classify} into
\begin{equation*}
\frac{\partial w}{\partial t}=\frac{\partial^2 w}{\partial x^2} - w^N
+ \sum_{i=1}^{N-1}a_i w^i + \sum_{k=0}^{N-1} w^k \left( -
\begin{pmatrix} N\\k \end{pmatrix} \mu^{N-k} + \sum_{i=k+1}^{N-1}
\begin{pmatrix} i \\ k \end{pmatrix} a_i \mu^{i-k} \right).
\end{equation*}
As is usual for the Fujita method, we choose a solution $v$ to
$\frac{\partial v}{\partial t}=-\frac{\partial^2
  v}{\partial x^2}$.  In particular, fix $T>0$ and choose
$\epsilon>0$ to define
\begin{equation*}
v_\epsilon(s,x)=\frac{1}{\sqrt{4\pi(T-s+\epsilon)}} e^{-\frac{1}{4(T-s+\epsilon)}(x-x_j)^2}.
\end{equation*}
Then we define $J(s)=\int v_\epsilon(s,x) w(s-t_j,x) dx$ and compute
\begin{eqnarray*}
\frac{dJ}{ds}&=&\int \frac{\partial v_\epsilon}{\partial s} w +
v_\epsilon \frac{\partial w}{\partial s} dx \\
&=&- \int v_\epsilon w^N dx
+ \sum_{i=1}^{N-1} \int a_i v_\epsilon w^i dx \\&&+ \sum_{k=0}^{N-1}
\int v_\epsilon w^k \left( -
\begin{pmatrix} N\\k \end{pmatrix} \mu^{N-k} + \sum_{i=k+1}^{N-1}
\begin{pmatrix} i \\ k \end{pmatrix} a_i \mu^{i-k} \right) dx\\
&\le&-J^N + \sum_{i=1}^{N-1} \|w\|_\infty^i \int a_i v_\epsilon dx \\&&+ \sum_{k=0}^{N-1}
\|w\|_\infty^k \left( \begin{pmatrix} N\\k \end{pmatrix} \mu^{N-k} + \sum_{i=k+1}^{N-1}
\begin{pmatrix} i \\ k \end{pmatrix} \|a_i\|_\infty \mu^{i-k} \right)\\
\end{eqnarray*}
where we have used Lemma \ref{finite_energy_consequences_lem} to bound $w$, and we have
used the assumption that $N$ is even in the last step.  Since $a_i$
decays to zero, the second term can be made arbitrarily small for an
arbitrarily large $s$ by taking $j$ large as well (for fixed $T$ and
$\epsilon$).  (The second term may eventually grow larger.)  The last
term is a constant, independent of $T,\epsilon$ and can be made
arbitrarily small by taking $j$ large.  By \eqref{u_upper_bound}, for
sufficiently large $j$, $J(0)<0$, and for larger $j$, $J(0)$ becomes
more negative.  Therefore, for a certain $T$ and sufficiently large
$j$, $J(s)$ tends to $-\infty$ for some $0<s<T$.  However, this
contradicts the fact that $w$ is bounded.
\end{enumerate}
\end{proof}
\end{thm}

\begin{cor}
Suppose that $n=1$ and that all equilibria have finite action.  If
$u$ is a finite energy solution then it converges uniformly to
equilibria as $|x|\to\infty$.  
\begin{proof}
Really, the only thing that must be noticed is that Theorem
\ref{unif_limits_to_equilibria} shows that there is uniform
convergence in the time direction.  For a given $\epsilon$, there is a
$T>0$ such that $|u(t,x)-f_\pm(x)|< \epsilon$ for all $|t|>T$.
However, this means that for $t \in [-T,T]$, this does not hold.
However, $[-T,T]$ is compact, and the proof of Lemma
\ref{parabolic_bootstrap_lem} indicates that there is uniform
convergence to equilibria as $\|x\|\to \infty$ on compact subsets.
\end{proof}
\end{cor}

\begin{cor}
\label{heteroclines_in_L1}
The above Corollary implies that the asymptotic spatial behavior of
heteroclinic orbits is determined entirely by the asymptotic spatial
behavior of equilibria.  In particular, in Chapter \ref{nonauto_ch},
it is shown that the equilibria for the case of \eqref{limited_pde}, 
\begin{equation*}
\frac{\partial u}{\partial t} = \frac{\partial^2 u}{\partial x^2} -
u^2 + \phi
\end{equation*}
with $\phi$ decaying to zero, all lie in $L^1(\mathbb{R})$.  As a
result, each timeslice of a heterocline lies in $L^1(\mathbb{R})$ as
well.
\end{cor}

\section{Discussion}

The point of employing the bootstrapping argument of Lemma
\ref{parabolic_bootstrap_lem} is only to extract uniform convergence
of the derivatives of the solution.  As can be seen from the proof of
Lemma \ref{finite_energy_consequences_lem}, such regularity arguments
are unneeded to obtain good convergence of the solution only.

While Corollary \ref{limits_to_equilibria} is probably true for
all spatial dimensions, the proof given here cannot be generalized to
higher dimensions.  In particular, V\'{e}ron in \cite{Veron_1996}
shows that in the case of $P(u)=-u^N$, there are solutions to the
equilibrium equation $\Delta u - u^N=0$ which are {\it unbounded
below} and {\it bounded above} when the spatial dimension is greater
than one.  This breaks the proof of Lemma
\ref{finite_energy_consequences_lem}, that the limiting equilibria of
finite energy solutions are bounded for $N$ even, since the proof
requires exactly the opposite.

On the other hand, the case of $P(u)=-u|u|^{N-1}+ \sum_{i=0}^{N-1}a_i
u^i$ is considerably easier than what we have considered here.  In
particular, all solutions to \eqref{pde_classify} are then bounded.  In that
case, the proof of Lemma \ref{finite_energy_consequences_lem} works
for all spatial dimensions.

\chapter{Equilibrium analysis}
\label{nonauto_ch}
\section{Introduction}

(This chapter is available on the arXiv as \cite{RobinsonNonauto}, and
has been accepted for publication in {\it Ergodic Theory and Dynamical
Systems}.)

Since the dynamics of solutions to the semilinear parabolic equation
\eqref{pde} depend strongly on the equilibrium solutions, it is important to
understand the number and structure of equilibrium solutions.  As will
be shown, this is a somewhat ill-defined and rather delicate goal.
Therefore, to fix ideas and techniques, we shall focus on the specific
case of the equilibria of \eqref{limited_pde}
\begin{equation}
\label{limited_pde_nonauto}
\frac{\partial u(t,x)}{\partial t}=\frac{\partial^2 u(t,x)}{\partial x^2} - u^2(t,x)
+ \phi(x),
\end{equation}
where $\phi$ tends to zero as $|x|\to\infty$.  The resulting questions
and techniques we encounter have obvious generalizations to the more
general equation.  Therefore, we are faced with the task of analyzing
a nonlinear {\it ordinary} differential equation, and finding its
{\it global} solutions.  Additionally, the asymptotic properties of
such solutions will be crucial in Chapters \ref{classify_ch},
\ref{global_ch}, and \ref{unstable_ch}.

Finding global solutions to nonlinear ordinary differential equations
on an infinite interval can be rather difficult.  Numerical
approximations can be particularly misleading, because they examine
only a finite-dimensional portion of the infinite-dimensional space in
which solutions lie.  Additionally, the conditions for global
existence can be rather delicate, which a numerical solver may have
difficulty rigorously checking.  In situations where there is
well-defined asymptotic behavior for global solutions, it is possible
to exploit the asymptotic information to answer questions about global
existence and uniqueness of solutions directly.  Additionally, more
detailed information may be provided by using the asymptotic behavior
to install artificial boundary conditions for use in a numerical
solver.  The numerical solver can then be used on the remaining
(bounded) interval with boundary conditions that match the numerical
approximation to an asymptotic expansion valid on the rest of the
solution interval.

In this chapter, we consider the behavior of global solutions
satisfying the equilibrium equation for \eqref{limited_pde_nonauto},
namely
\begin{equation}
\label{long_ode1}
0=f''(x)-f^2(x)+\phi(x),\text{ for all }x\in\mathbb{R}.
\end{equation}
In particular, we wish to know how many solutions there are for a
given $\phi$.  (There may be uncountably many solutions, as in the
case where $\phi \equiv \text{const}>0$.)  This problem depends rather
strongly on the asymptotic behavior of solutions to \eqref{long_ode1}
as $|x| \to \infty$, so it is useful to study instead the pair of
initial value problems
\begin{equation}
\label{ode1}
\begin{cases}
0=f''(x) - f^2(x) + \phi(x) \text{ for }x>0\\
(f(0),f'(0))\in Z,
\end{cases}
\end{equation}
and
\begin{equation}
\begin{cases}
\label{ode1_backwards}
0=f''(x) - f^2(x) + \phi(x) \text{ for }x<0\\
(f(0),f'(0))\in Z',
\end{cases}
\end{equation}
where $\phi \in C^\infty(\mathbb{R}).$ The sets $Z$,$Z'$ supply the
  initial conditions for which solutions exist to \eqref{ode1} for all
  $x>0$ and to \eqref{ode1_backwards} for all $x<0$, respectively.
  Solutions to \eqref{long_ode1} will occur exactly when $Z \cap Z'$
  is nonempty.  Indeed, the theorem on existence and uniqueness for
  ODE gives a bijection between points in $Z\cap Z'$ and solutions to
  \eqref{long_ode1} \cite{LeeSmooth}. Since \eqref{ode1} and
  \eqref{ode1_backwards} are related by reflection across $x=0$, it is
  sufficient to study \eqref{ode1} only.

Due to the asympotic behavior of solutions to \eqref{ode1}, the
methods we employ here will be most effective in the specific cases
where $\phi$ is nonnegative and monotonically decreasing to zero.  (We
denote the space of smooth functions that decay to zero as
$C^\infty_0(\mathbb{R})$.)  The decay condition on $\phi$ allows the
differential operators in \eqref{long_ode1} through
\eqref{ode1_backwards} to be examined with a perturbative approach as
$x$ becomes large, and makes sense if one is looking for smooth
solutions in $L^p(\mathbb{R})$ with bounded derivatives.  

When $\phi$ is strictly negative, it happens that no solutions exist
to \eqref{ode1} for all $x>0$.  The monotonicity restriction on $\phi$
provides some technical simplifications and sharpens the results that
we obtain.  This leads us to restrict $\phi$ to a class of functions
that captures this monotonicity restriction but allows some
flexibility, which we shall call the {\it M-shaped functions}.

It is unlikely that we will be able to solve \eqref{ode1} explicitly
for arbitrary $\phi$, so one might think that numerical approximations
might be helpful.  However, most numerical approximations will not be
able to count the number of global solutions accurately.  For
instance, finite-difference methods are typically only useful for
finding solutions valid on finite intervals of $\mathbb{R}$.  However,
one cannot easily infer a solution's behavior for large values of
$|x|$ when it is only known on a finite interval.  In particular,
global solutions to \eqref{ode1} must tend to zero (Theorem
\ref{limzero_iff_bounded_lem}).  All other solutions fail to exist for
all of $\mathbb{R}$.  Worse, the space of initial conditions which
give rise to global solutions is at best a 1-dimensional submanifold
of the 2-dimensional space of initial conditions (Theorem
\ref{properties_of_Z}).  Therefore, a typical finite-difference
solution that {\it appears} to tend to zero may in fact not, and as a
result fails to be a solution over all $x>0$.

Because of this failure, we need to understand the asymptotic behavior
of solutions to \eqref{ode1} as we take $x \to \infty$.  Equivalently,
since $\phi \in C_0^\infty(\mathbb{R})$, this means that we should
examine solutions with $\phi$ small.  The driving motivation for this
discussion is that solutions to $0=f''(x)-f^2(x)+\phi(x)$ for $\phi$
small behave much like solutions to $0=f''(x)-f^2(x)$.  In the latter
case, we can completely characterize the solutions which exist on
intervals like $[x_0,\infty)$.

In Section \ref{review_sec} we review what is known about the much
simpler case where $\phi$ is a constant.  Of course, then \eqref{ode1}
is autonomous, and the results are standard.  In Section
\ref{asymp_exist_sec}, we establish the existence of solutions which
are asymptotic to zero.  Some of these solutions are computed
explicitly using perturbation methods in Section \ref{series_sec},
where low order approximations are used to gather qualitative
information about the initial condition sets $Z$ and $Z'$.  In
Sections \ref{extension_sec} and \ref{geom_props_Z}, these qualitative
observations are made precise.  Section \ref{long_sec} applies these
observations about $Z$ and $Z'$ to give existence and uniqueness
results for \eqref{long_ode1}.  Finally, in Section \ref{numer_sec},
we use the information gathered about $Z$ and $Z'$ to provide
artificial boundary conditions to a numerical solver on a bounded
interval, which sharpens the results from Section \ref{long_sec}.  We
exhibit the numerical results for a typical family of $\phi$, showing
bifurcations in the global solutions to \eqref{long_ode1}.

\section{Review of behavior of solutions to $0=f''(x)-f^2(x)+P$}
\label{review_sec}

It will be helpful to review the behavior of 
\begin{equation}
\label{IVP2}
\begin{cases}
0=f''(x) - f^2(x) + P\\
f(0),f'(0)\text{ given},
\end{cases}
\end{equation}
where $P$ is a constant, since varying $\phi$ can be viewed as a
perturbation on the case $\phi(x)=P$.  In particular, we need to
compute some estimates for later use.  We shall typically take $P>0$,
as there do not exist solutions for all $x$ if $P<0$.  




\begin{figure}
\begin{center}
\includegraphics[height=3in]{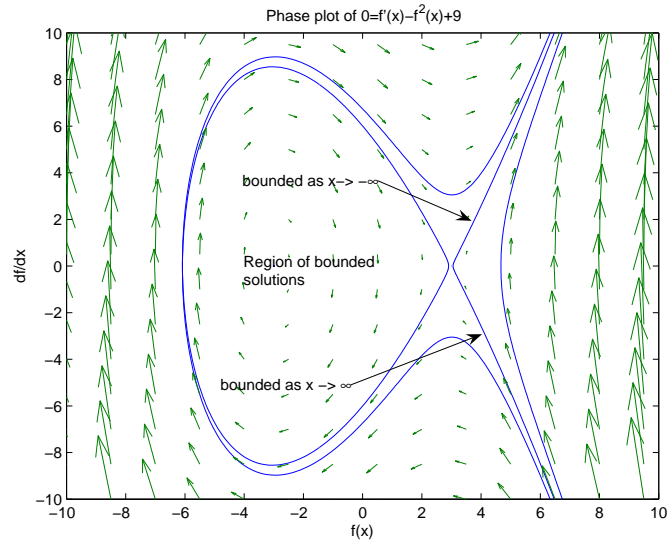}
\end{center}
\caption{The phase plot of $f''-f^2+9=0$.  Bounded solutions live in a
small region, the rest are unbounded.}
\label{P_eq_p9}
\end{figure}

\begin{lem}
\label{no_upper_f_bnd_lem}
Suppose $f$ is a solution to the initial value problem
\eqref{IVP2} with $f(0)>\sqrt{P}$ and $f'(0)>0$.  Then there does
not exist an upper bound on $f(x)$, when $x>0$.  Additionally, if
$P<0$, there does not exist an upper bound on $f(x)$.
\begin{proof}
Observe that for $f>\sqrt{P}$ or if $P<0$
\begin{equation*}
f''=f^2-P>0.
\end{equation*}
Hence, since $f'(0)>0$, and $f'$ is monotonic increasing, $f(x)$ is
monotonic increasing at an increasing rate.  Thus it must be
unbounded from above.
\end{proof} 
\end{lem}

\begin{df}
The differential equation \eqref{IVP2} comes from a Hamiltonian, namely
\begin{equation*}
H(f,f')=\frac{1}{3} f^3 - \frac{1}{2}f'^2-fP+\frac{2}{3}P^{3/2}.
\end{equation*}
\end{df}

\begin{df}
A useful tool in the study of smooth dynamical systems is the {\it
funnel}.  Suppose $\Phi$ is a local flow on a manifold $M$.  A funnel
$F$ is a set such that if $x \in F$, then $\Phi_x(t) \in F$ for all
$t>0$.  A funnel $F$ with an oriented, piecewise $C^1$ boundary is
characterized by having the vector field $\frac{d}{dt} \Phi_x$ being
inward-pointing for all $x \in \partial F$.
\end{df}

\begin{lem}
\label{bounded_in_funnel_lem}
Suppose $f$ is a solution to the equation \eqref{IVP2} on
$\mathbb{R}$.  All bounded solutions lie in the funnel
\begin{equation}
\label{M}
M=\{(f,f')|H(f,f')\ge 0 \text{ and } f \le \sqrt{P}\}.
\end{equation}
Any solution which includes a point outside the closure of $M$ is
unbounded, either for $x>0$ or $x<0$.  (Note that $M$ is the
teardrop-shaped region in Figure \ref{P_eq_p9}.)
\begin{proof}
\begin{itemize}
\item $M$ is a bounded set.  Notice that $H(f,0) \ge H(f,f')$, or in other
  words within $M$,
\begin{equation*}
0<\frac{1}{3}f^3 - \frac{1}{2}f'^2-f P + \frac{2}{3}P^{3/2}\le
 \frac{1}{3} f^3 - fP + \frac{2}{3}P^{3/2}.
\end{equation*}
Elementary calculus reveals that this inequality establishes a
lower bound on $f$, namely that 
\begin{equation}
\label{f_bound}
-\sqrt{3P} \le f \le \sqrt{P}
\end{equation} 
On the other hand,
\begin{equation}
\label{fp_bound}
|f'|<\sqrt{\frac{4}{3}P^{3/2}+\frac{2}{3}f^3-2fP} \le \sqrt{\frac{8}{3}}P^{3/4}
\end{equation}
immediately establishes a bound on $f'$.

\item $M$ is a funnel, from which solutions neither enter nor leave.
This is immediate from the fact that $H$ is the Hamiltonian, and the
definition of $M$ simply says that $H(f,f')\ge 0$.
This suffices since solutions to \eqref{IVP2} are tangent to level
curves of $H$.

\item If $(f(0),f'(0)) \notin M$ then $f$ is unbounded.  Evidently if
  $f(0)> \sqrt{P}$ and $f'(0)>0$, then Lemma \ref{no_upper_f_bnd_lem}
  applies to give that $f$ is unbounded.  For the remainder, discern
  two cases.  First, suppose $f(0)>\sqrt{P}$ and $f'(0)<0$.
  Evidently, $H(f(0),f'(0))=H(f(0),-f'(0))$, so it's just a matter of
  verifying that a solution curve transports our solution to the first
  quadrant.  But this is immediately clear from the formula for
\begin{equation*}
f'=\pm \sqrt{\frac{2}{3}f^3-2fP-2H(f(0),f'(0))},
\end{equation*}
which gives $f'=\pm f'(0)$ when $f=f(0)$.  The other case is when
   $H(f(0),f'(0)) < 0$.  Then we show that there
   is a point $(\sqrt{P},g)$ on the same solution curve, and then
   Lemma \ref{no_upper_f_bnd_lem} applies.  So we try to satisfy
\begin{eqnarray*}
\frac{1}{3}P^{3/2} - \frac{1}{2}g^2-P^{3/2} + \frac{2}{3} P^{3/2} &=&
H(f(0),f'(0)) < 0 \\
g^2 &=& -2 H(f(0),f'(0)) > 0,
\end{eqnarray*}
which clearly has a solution in $g$.
\end{itemize}
\end{proof}
\end{lem}

\begin{lem}
\label{asymptote_lem}
If $f$ is a solution to \eqref{IVP2} with $f(0)>\sqrt{P}$, and
$f'(0)>0$ then there exists a $C$ such that $\lim_{x \to C} f(x)
= \infty$. 
\begin{proof}
From Lemma \ref{bounded_in_funnel_lem}, we have that $f$
is unbounded, and goes to $+\infty$.  Using the Hamiltonian, we
can solve for
\begin{equation*}
\frac{df}{dx} = \pm \sqrt{\frac{2}{3} f^3 - 2 f P - 2 H(f(0),f'(0))},
\end{equation*}
or viewing $f$ as the independent variable,
\begin{eqnarray*}
\frac{dx}{df} &=& \frac{1}{\pm \sqrt{\frac{2}{3} f^3 - 2 f P - 2
    H(f(0),f'(0))}}\\ 
&\sim& \sqrt{\frac{3}{2}} f^{-3/2},
\end{eqnarray*}
as $f$ becomes large.  Solving this asymptotic differential equation
is easy, and leads to
\begin{eqnarray*}
x &\sim& -\frac{1}{2}\sqrt{\frac{3}{2 f}}+C,\\
f &\sim& \frac{3}{8(x-C)^2}, \text{ (for $|x-C|$ small)}
\end{eqnarray*}
which has an asymptote at $x=C$.
\end{proof}
\end{lem}

\section{Existence of asymptotic solutions for $\phi \in
  C_0^\infty(\mathbb{R})$} 
\label{asymp_exist_sec}

The first collection of results we obtain will make the assumption
that $\phi$ tends to zero.  From this, a number of useful asymptotic
results follow.  Working in the phase plane will be useful for
understanding \eqref{ode1}.  Of course \eqref{ode1} is not autonomous,
but by adding an additional variable, it becomes so.

\begin{df}
We think of \eqref{ode1} as a vector field $V$ on $\mathbb{R}^3$, defined
by the formula
\begin{equation}
\label{ode1_system}
V(f,f',x)=\begin{pmatrix}
f' \\ f^2 - \phi(x) \\ 1
\end{pmatrix}.
\end{equation}
Notice that the first coordinate of an integral curve for this vector
field solves \eqref{ode1}.
\end{df}

\begin{df}
Define $H(f,f',x)=\frac{1}{3}f^3-\frac{1}{2}f'^2-f \phi(x) +
\frac{2}{3} \phi^{3/2}(x).$  Notice that for constant $\phi = P$, this
reduces to a Hamiltonian for \eqref{IVP2}.  
\end{df}

\begin{thm}
\label{limzero_iff_bounded_lem}
Suppose $f$ is a solution to the problem \eqref{ode1} where $\phi \in
C_0^\infty(\mathbb{R})$.  If $f$ does not tend to zero as
$x\to \infty$, then there exists a $z$ such that $\lim_{x\to z} f(x) =
  \infty$.  Stated another way, if $f$ solves \eqref{ode1} for all
  $x>0$, then $\lim_{x \to \infty} f(x) = 0$.
\begin{proof}
If $f$ does not tend to zero, this means that there is an $R>0$ such
that for each $x_0>0$, there is an $x>x_0$ so that $|f(x)|>R$.  But
since $\phi$ tends to zero as $x \to \infty$, for any $P>0$ we can
find an $x_1>0$ such that for all $x>x_1$, $|\phi(x)|<P$.  Choose such
a $P$ so that the set $M$ in Lemma \ref{bounded_in_funnel_lem}
associated to \eqref{IVP2} is contained entirely within the strip
$-R<f<R$.  We can do this since the set $M$ is bounded, and its radius
decreases with decreasing $P$, as shown in \eqref{f_bound} and
\eqref{fp_bound}.  But this means that there is an $x_2>x_1$ such that
$|f(x_2)|>R$.

Construct the following regions (See Figure \ref{limzero_fig}):
\begin{equation*}
I=\{(f,f',x)| f \ge R \text{ and } f' \le 0\},
\end{equation*}
\begin{equation*}
II=\{(f,f',x)| f \ge R \text{ and } f' \ge 0\},
\end{equation*}
\begin{equation*}
III=\{(f,f',x)| f \le -R \},
\end{equation*}
and
\begin{equation*}
\begin{split}
IV=\left\{(f,f',x)| f' \ge 0 \text{ and } f\ge -R \text{ and } \left(
\frac{1}{3}f^3 - \frac{1}{2} f'^2 - f P + \frac{2}{3} P^{3/2} \le 0 \right. \right.
\\ \left. \left. \text{ if } f \le \sqrt{P}\right) \right\}.
\end{split}
\end{equation*}

\begin{figure}
\begin{center}
\includegraphics[height=3in]{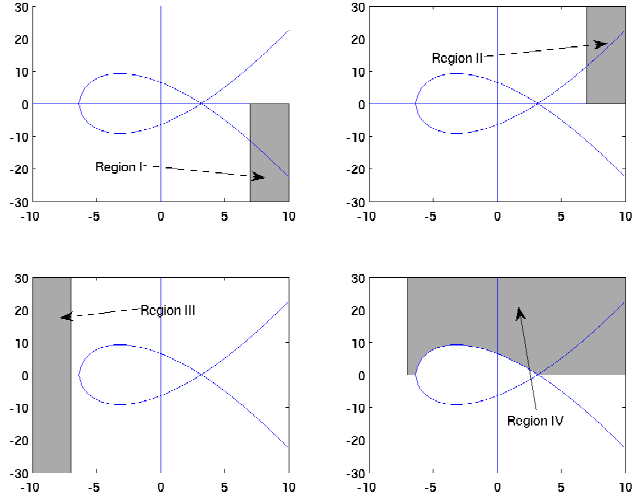}
\end{center}
\caption{The Regions $I$, $II$, $III$, and $IV$ of Theorem \ref{limzero_iff_bounded_lem}}
\label{limzero_fig}
\end{figure}

The following statements hold:
\begin{itemize}
\item Region $I$ is an antifunnel.  Along $f=R$ and $f'=0$, solutions
  must exit.  Once a solution exits Region $I$, it cannot reenter.
  Also, because $f > \sqrt{P}$, $f''=f^2-\phi> f^2 - P > 0$, solutions
  must exit Region $I$ in finite $x$.
\item Region $II$ is a funnel.  Along $f=R$ and $f'=0$, solutions
  enter.  Now $f''=f^2 - \phi > f^2 - P \ge 0$ and $f' \ge 0$, so
  solutions will increase at an increasing rate and so, they
  are unbounded. 
\item Solutions remain in Region $III$ for only finite $x$, after
  which they must enter Region $IV$.  This occurs since $f \le
  -\sqrt{P} < 0$, and so $f'$ always increases.  Note that for $f'<0$,
  solutions will enter Region $III$ along $f=-R$, and for $f'>0$,
  solutions exit along $f=-R$.
\item Region $IV$ is a funnel.  Solutions enter along $f=-R$ and along
  $f'=0$ (note that $|f| \ge \sqrt{P}$ in both cases).  Along the
  curve boundary of Region $IV$, we have that
\begin{eqnarray*}
\nabla \left(\frac{1}{3} f^3 - \frac{1}{2} f'^2 - fP +
\frac{2}{3}P^{3/2} \right) \cdot  V(f,f',x)&=&
\begin{pmatrix} f^2 - P \\ -f' \\ 0 \end{pmatrix}^T \begin{pmatrix} f' \\ f^2 - \phi \\ 1
\end{pmatrix}\\
&=&f'(P-\phi)< 0,
\end{eqnarray*}
so that solutions enter.
\end{itemize}

Now suppose $(f(x_2),f'(x_2),x_2)\in I$.  After finite $x$, say at
$x=x_2'$, the solution through that point must exit Region $I$, never
to return.  Then, there is an $x_3 > x_2'$ such that $|f(x_3)|>R$.  So
this solution has either $(f(x_3),f'(x_3),x_3) \in II$ or $\in III$.
The former gives the conclusion we want, so consider the latter case.
The solution will only remain in Region $III$ for finite $x$, after
which it enters Region $IV$, say at $x=x_3'$.  Then there is an
$x_4>x_3'$ such that $|f(x_4)|>R$.  Now the only possible location for
$(f(x_4),f'(x_4),x_4)$ to be is within Region $II$, since it must also
remain in Region $IV$.  As a result, the solution is unbounded by an
easy extension of Lemma \ref{no_upper_f_bnd_lem}.  As $x$ becomes
large, $\phi$ tends to zero, so the solution will be asymptotic to an
unbounded solution of $0=f''-f^2$.  But Lemma \ref{no_upper_f_bnd_lem}
above assures us that such a solution is unbounded from above, and
Lemma \ref{asymptote_lem} gives that it has an asymptote.  Hence, our
solution must blow up at a finite $x$.
\end{proof}
\end{thm}

This result indicates that solutions to \eqref{ode1} which exist for
all $x>0$ are rather rare.  Those which exist for all $x>0$ must tend
to zero, and it seems difficult to ``pin them down.''  We now apply
topological methods, similar to those employed in \cite{HubbardWest},
to ``capture'' the solutions we seek.  The methods we use are due to
Wa\.{z}ewski \cite{Wazewski}.

We begin by extending the usual definition of a flow slightly to the
case of a manifold with boundary.

\begin{df}
Suppose $M$ is a manifold with boundary.  A {\it flow domain} $J$ is a
  subset of $\mathbb{R} \times M$ such that if $x \in M$ then
  $J_x=\text{pr}_1 (J \cap \mathbb{R} \times \{x\})$ is an interval
  containing 0, and if $x$ is in the interior of $M$ then 0 is
  in the interior of $J_x$.  ($\text{pr}_1:\mathbb{R}\times M \to
  \mathbb{R}$ is projection onto the first factor)
\end{df}

\begin{df}
A {\it (smooth) flow} is a smooth map $\Phi$ from a flow domain $J$
to a manifold with boundary $M$, satisfying
\begin{itemize}
\item $\Phi(0,x)=x$ for all $x \in M$ and
\item $\Phi(t_1+t_2,x)=\Phi(t_1,\Phi(t_2,x))$ whenever both
  sides are well-defined.
\end{itemize}
Additionally, we assume that flows are {\it maximal} in the sense that
they cannot be written as a restriction of a map from a larger flow
domain which satisfies the above axioms.  We call the curve $\Phi_x:J_x
\to M$ defined by $\Phi_x(t)=\Phi(t,x)$ the {\it integral curve
through $x$} for $\Phi$.
\end{df}

\begin{df}
Suppose $\Phi:J \to M$ is a flow on $M$ and $x \in \partial M$.  Then
the flow at $x$ is said to be {\it inward-going} (or simply {\it
inward}) if $J_x$ is an interval of the form $[0,a)$ or $[0,a]$ for
some $0<a\le \infty$.  Likewise, the flow at $x$ is {\it outward-going} if
$J_x$ is of the form $(a,0]$ or $[a,0]$ for $-\infty \le a < 0$.
\end{df}

\begin{thm} (Wa\.{z}ewski's antifunnel theorem)
\label{wazewski_thm}
Suppose $\Phi:J \to M$ is a flow on $M$ and that $\{A,B\}$ forms a
partition of the boundary of $M$ such that the flow of $\Phi$ is
inward along A and outward along B.  If every integral curve of $\Phi$
intersects $B$ in finite time (ie. $J_x$ is bounded for each $x$),
then $A$ is diffeomorphic to $B$.
\begin{proof}
For each $x \in A$, $J_x=[0,t_x]$, where $t_x$ is the time
which the integral curve through $x$ intersects $B$.  (We have that
$\Phi(t_x,x)$ is outward-going, since $J_x$ is closed, so it is in
$B$.)

Using this, we can define a map $F:A \to B$ by $F(x)=\Phi(t_x,x)$.
$F$ takes $A$ smoothly and injectively into $B$.  The
smoothness follows from the smoothness of $\Phi$ and that $\partial M$
is a smooth submanifold.  To see the injectivity, suppose $F(x)=F(y)$ for
some $x,y \in A$, so $\Phi(t_x,x)=\Phi(t_y,y)$.  Without loss of
generality, suppose $0 < t_x \le t_y$.  Then we have that
\begin{eqnarray*}
F(x)&=&F(y)\\
\Phi(-t_x,F(x))&=&\Phi(-t_x,F(y))\\
\Phi(-t_x,\Phi(t_x,x))&=&\Phi(-t_x,\Phi(t_y,y))\\
\Phi(t_x-t_x,x)&=&\Phi(t_y-t_x,y)\\
x&=&\Phi(t_y-t_x,y).
\end{eqnarray*}
But the flow is inward at $x$, so it is also inward at
$\Phi(t_y-t_x,y)$.  This means that $(t_y - t_x - \epsilon,y) \notin J$
for every $\epsilon > 0$.  But this contradicts the fact that $(t_y,y)\in J$
unless we have $t_y \le t_x$.  As a result, $t_y = t_x$, so $x=y$.

In just the same way as for $F$, we construct a map $G:B \to A$ so
that $G$ takes $B$ smoothly and injectively into $A$.  Namely, we
suppose $J_y=[s_y,0]$ for some $s_y$, and put $G(y)=\Phi(s_y, y)$.
Notice that by maximality, if there were to be an $x \in A$ such that
$F(x)=y$, $s_y = - t_x$.

Now we claim that $G$ is the inverse of $F$. We have that 
\begin{eqnarray*}
(G\circ F)(x)&=&\Phi(s_{F(x)},F(x)) \\
&=&\Phi(s_{F(x)},\Phi(t_x,x))\\
&=&\Phi(s_{F(x)}+t_x,x)\\
&=&\Phi(-t_x+t_x,x)=x,\\
\end{eqnarray*}
where we employ the remark about $s_y$ above.
\end{proof}
\end{thm}

\begin{rem}
We can extend the Antifunnel theorem to a topological space $X$ on
which a flow $\Phi:J \to X$ acts in the obvious way.  In that case,
there is no reasonable definition of the boundary of $X$.  However,
the notion of inward- and outward-going points still makes sense.  If
we let $A$ be the set of inward-going points and $B$ be the set of
outward-going points in $X$, then the conclusion is that $A$ is
homeomorphic to $B$.
\end{rem}

Now we employ the Antifunnel theorem to deduce the existence of
a bounded solution to $0=f''-f^2+\phi$ for $x>x_0$ for some $x_0\ge 0$.

\begin{figure}
\begin{center}
\includegraphics[height=3in]{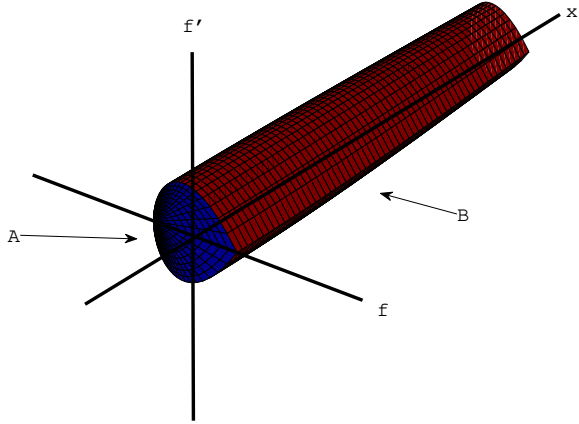}
\end{center}
\caption{Schematic of the region $R_1$, showing the boundary partition
$A$ and $B$.}
\label{region_r1_fig}
\end{figure}

\begin{thm}
\label{region_r1_thm}
Suppose $0 \le \phi(x) \le K$ for all $x \ge x_0$ for some $x_0>0$ and
$0 < K < \infty$, and that there exists an $x_1 \ge x_0$ such that for all
$x>x_1$, $\phi(x)>0$.  Then the region $R_1$ given by
$R_1=\{(f,f',x)|H(f,f',x) \ge 0, x \ge 0, f\le \sqrt{\phi(x)}\}$
contains a bounded solution to $0=f''(x)-f^2(x)+\phi(x)$, which exists
for all $x$ greater than some nonnegative $x_2$.
\begin{proof}
Without loss of generality, we may take $x_0=0$, because otherwise
solutions must exit the portions of $R_1$ in $\{(f,f',x)|x<x_0\}$
since the $x$-component of $V(f,f',x)$ is equal to 1.

If $\phi(0)>0$, partition the boundary of $R_1$ into two pieces:
$A=\{(f,f',x)|x = 0\}$ and $B=\{(f,f',x)|H(f,f',x)=0\}$ (See Figure
\ref{region_r1_fig}).  The flow of $V$ is evidently inward along $A$.
As for $B$, notice that $\nabla H$ is an inward-pointing vector field
normal to $B$.  We compute
\begin{eqnarray*}
\nabla H \cdot V &=& 
\begin{pmatrix} f^2 - \phi(x) \\ -f' \\ (-f + \sqrt{\phi(x)})\phi'(x) 
  \end{pmatrix}^T
\begin{pmatrix} f' \\ f^2 - \phi(x) \\ 1 \end{pmatrix}\\
&=& ( - f + \sqrt{\phi(x)} ) \phi'(x),
\end{eqnarray*}
which has the same sign as $\phi'(x)$ when $f < \sqrt{\phi(x)}$ in
$R_1$.  Finally, we must deal with the case where $f=\sqrt{\phi(x)}
\in B$.  But in this case, $f'=0$ from the equation for $H$, so we see
that $V(\sqrt{\phi(x)},0,x)=(0,0,1)^T$, so the flow is inward when
$\phi'(x)<0$ and outward when $\phi'(x)>0$.  This means that
the portion of the boundary of $R_1$ on which the flow is outward is a
disjoint union of annuli.  On the other hand, the portion of the
boundary of $R_1$ on which the flow is inward is the disjoint union of
a disk (namely $R_1 \cap \{x=0\}$) and some annuli.  

We now consider the case of $\phi(0)=0$, in which case the set $A$
above is just a point.  Assume without loss of generality that
$\phi(x)$ is strictly positive for all $x>0$, so we let $x_1=0$.  Let
\begin{equation*}
x'=\begin{cases}
\text{inf } \{x \in (0,\infty) | \phi'(x)=0\} \text{ or } \\
\infty \text{ if } \phi'(x) > 0 \text{ for all }x>0.
\end{cases}
\end{equation*}
In this case, the set $\{(f,f',x)|0 \le x \le x'\} \cap \partial R_1$,
is a contractible (it may be a point if $\phi$ oscillates rapidly as
$x\to 0$), connected component of the inflow portion of the boundary
of $R_1$.  It is obvious that the remainder of the inflow portion of
the boundary is homeomorphic to a disjoint union of annuli, since
$\phi$ is smooth and strictly positive.

We can apply the Antifunnel theorem to conclude that there is a
solution which does not intersect either the inflow or outflow
portions of the boundary.  There is a lower bound on the
$x$-coordinate of such a solution, since the $x$-component of
$V(f,f',x)$ is equal to 1, and the Region $R_1$ lies within the
half-space $x>0$.  Therefore, there must exist a solution which
enters $R_1$, and remains inside the interior of $R_1$ for all larger
$x$.  That such a solution is bounded follows from the fact that each
constant $x$ cross section of $R_1$ has a radius bounded by the
inequalities \eqref{f_bound} and \eqref{fp_bound}, and the fact that
$\phi(x) \le K < \infty$.  
\end{proof}
\end{thm}

\section{Asymptotic series solution}
\label{series_sec}

Theorem \ref{region_r1_thm} ensures the existence of solutions to
$0=f''-f^2+\phi$ for $x$ sufficiently large.  However, it does not
give any description of the initial condition set $Z$ which leads to
such solutions, nor does it give a description of the maximal
intervals of existence.  Fortunately, it is relatively easy to
construct an asymptotic series for solutions to \eqref{ode1}, which
will provide a partial answer to this concern.  In doing so, we
essentially follow standard procedure, as outlined in \cite{Holmes},
for example.  However, our case is better than the standard situation,
because under relatively mild restrictions this series {\it converges}
to a true solution.

We begin by supposing that our solution has the form
\begin{equation}
\label{f_format}
f=\sum_{k=0}^\infty f_k,
\end{equation}
where we temporarily assume $f_{k+1} \ll f_k$ and $f_0 \gg \phi$, as
$x \to +\infty$.  (This assumption will be verified in Lemma
\ref{series_bound_lem}.)  Substituting \eqref{f_format} into
\eqref{ode1}, we get

\begin{eqnarray*}
0&=&\sum_{k=0}^\infty \left[ f_{k}'' - \sum_{m=0}^{k}f_m
  f_{k-m}\right]+ \phi \\
0&=&f_0''-f_0^2+(f_1''-2f_0f_1+\phi)+\sum_{k=2}^\infty
  \left[f_{k}''-2f_0f_{k}-\sum_{m=1}^{k-1}f_mf_{k-m}\right].
\end{eqnarray*}
We solve this equation by setting different orders to zero.
Namely,
\begin{eqnarray*}
0&=&f_0''-f_0^2\\
0&=&f_1''-2f_0f_1+\phi\\
0&=&f_k''-2f_0f_k-\sum_{m=1}^{k-1}f_mf_{k-m}.
\end{eqnarray*}
The equation for $f_0$ is integrable, and therefore easy to solve.
(There are two families of solutions for $f_0$.  We select the
nontrivial one, because the other one simply results in $f(x) \sim
-\int_x^\infty \int_t^\infty \phi(s) ds\, dt$.)  The equations for
$f_k$ are linear and can be solved by a reduction of order.  Thus
formally, the solutions are
\begin{equation}
\label{series_soln_coef}
\begin{cases}
f_0=\frac{6}{(x-d)^2}\\
f_1=\frac{1}{(x-d)^3}\left[ K + \int^x(t-d)^6\int_t^\infty
  \frac{\phi(s)}{(s-d)^3}ds\,dt\right]\\
f_k=-\frac{1}{(x-d)^3}\int^x (t-d)^6 \int_t^\infty \frac{ \sum_{m=1}^{k-1}
  f_m(s) f_{k-m}(s)}{(s-d)^3}ds\,dt,
\end{cases}
\end{equation}
for $d,K$ constants.  Notice that these constants parametrize the set of
initial conditions $Z$.  

\begin{lem}
\label{series_bound_lem}
Suppose $f(x)=\sum_{k=0}^\infty f_k(x)$ where the $f_k$ are given by
\eqref{series_soln_coef}.  If there
exists an $M>0$, an $R>0$, and an $\alpha>5$ such that 
\begin{equation}
\label{phi_cond}
|\phi(x)|<\frac{M}{(x-d)^\alpha} \text{ for all } |x-d|>R>0,
\end{equation}
then $f(x)$ is bounded above by the power series
\begin{equation}
\label{f_power_series_bound}
|f(x)| \le \frac{1}{(x-d)^2} \sum_{k=0}^\infty \left|\frac{A_k}{x-d}\right|^k.
\end{equation}
\begin{proof}
We proceed by induction, and begin by showing that the $f_1$ term is
appropriately bounded:

\begin{eqnarray*}
|f_1(x)|&=&\left|\frac{1}{(x-d)^3}\left[K+\int^x(t-d)^6\int_t^\infty
  \frac{\phi(s)}{(s-d)^3}ds\,dt\right]\right|\\
&\le&\left|\frac{1}{(x-d)^3}\left[K+\int^x(t-d)^6\int_t^\infty
  \frac{M}{(s-d)^{3+\alpha}}ds\,dt\right]\right|\\
&\le&\left|\frac{1}{(x-d)^3}\left[K+
  \frac{M}{(2+\alpha)(5-\alpha)(x-d)^{\alpha-5}} \right]\right|\\
\end{eqnarray*}
Now since $|x-d|>R$ and $\alpha > 5$, we have that
\begin{eqnarray*}
|f_1(x)|&\le& \frac{1}{|x-d|^3} \left[|K|+
  \frac{M}{(2+\alpha)|5-\alpha|R^{\alpha-5}} \right].\\
&\le&\frac{A_1}{|x-d|^3}.
\end{eqnarray*}
with
\begin{equation}
\label{A_1_def}
A_1 = |K|+\frac{M}{(2+\alpha)(\alpha-5) R^{\alpha-5}}.
\end{equation}

For the induction hypothesis, we assume that
$|f_i|\le\frac{A_i}{|x-d|^{2+i}}$ with $A_i \ge 0$ and for all $i \le
k-1$.  We have that

\begin{eqnarray*}
\sum_{m=1}^{k-1} f_m f_{k-m} &\le& \sum_{m=1}^{k-1}
\frac{A_m}{|x-d|^{2+m}} \frac{A_{k-m}}{|x-d|^{2+k-m}} \\
&\le& \frac{1}{|x-d|^{k+4}} \sum_{m=1}^{k-1} A_m A_{k-m}, \\
\end{eqnarray*}
so by the same calculation as for $f_1$, we obtain
\begin{equation*}
f_k \le \frac{\sum_{m=1}^{k-1} A_m A_{k-m}}{(k+6)(k-1)}\frac{1}{|x-d|^{k+2}}.
\end{equation*}

Hence we should take
\begin{equation}
\label{A_k_def}
A_k = \frac{\sum_{m=1}^{k-1} A_m A_{k-m}}{(k+6)(k-1)}.
\end{equation}

Hence we have that 

\begin{equation*}
|f(x)| \le \sum_{k=0}^\infty |f_k(x)|\le
 \frac{1}{|x-d|^2}\sum_{k=0}^\infty A_k \left|\frac{1}{x-d}\right|^k.
\end{equation*}
\end{proof}
\end{lem}

\begin{lem}
\label{series_conv_lem}
The power series given by 
\begin{equation*}
\sum_{k=0}^\infty \frac{A_k}{|x-d|^k}, 
\end{equation*}
with $A_0, A_1 \ge 0$ given, and 
\begin{equation*}
A_{k}=\frac{\sum_{m=1}^{k-1} A_m A_{k-m}}{(k+6)(k-1)} =
\frac{\sum_{m=1}^{k-1} A_m A_{k-m}}{k^2+5k-6}
\end{equation*}
converges for $|x-d|>R$ if $A_1 \le 8R$.
\begin{proof}
We show that under the conditions given, the series passes the usual
ratio test.  That is, we wish to show
\begin{equation*}
\lim_{k \rightarrow \infty} \left|\frac{A_{k+1}}{A_k}\right| \le R.
\end{equation*}
Proceed by induction.  Take as the base case, $k=1$: by the formula for
$A_{k}$,
\begin{equation*}
A_2=\frac{A_1^2}{8}, \text{ so }
\frac{A_2}{A_1}=\frac{A_1^2}{8A_1}=\frac{A_1}{8} \le R.
\end{equation*}

Then for the induction step, 
\begin{eqnarray*}
\frac{A_{k+1}}{A_k}&=&\frac{\sum_{m=1}^k A_mA_{k-m+1}}{A_k (k^2+7k)}\\
&=&\frac{\sum_{m=2}^k A_mA_{k-m+1}+A_1 A_k}{A_k (k^2+7k)}\\
&=&\frac{\sum_{m=1}^{k-1} A_{m+1}A_{k-m}+A_1 A_k}{A_k (k^2+7k)}\\
&\le&\frac{R\sum_{m=1}^{k-1} A_mA_{k-m}+A_1 A_k}{A_k (k^2+7k)}\\
&\le&\frac{RA_k(k^2+5k-6)+A_1 A_k}{A_k (k^2+7k)}\\
&\le&\frac{R(k^2+5k-6)+A_1}{(k^2+7k)}\\
&\le&\frac{R(k^2+5k-2)}{(k^2+7k)}\\
&\le&R,\\
\end{eqnarray*}
since $A_1 \le 8 R$.  Thus $\left| \frac{A_{k+1}}{A_k} \right| \le R$
for all $k$, so the power series converges.
\end{proof}
\end{lem}

Lemma \ref{series_conv_lem} provides conditions for the convergence of the
bounding series found in Lemma \ref{series_bound_lem}.  Hence we have actually proven the
following:

\begin{thm}
\label{series_conv_thm}
Suppose $f(x)=\sum_{k=0}^\infty f_k(x)$ where the $f_k$ are given by
\eqref{series_soln_coef}.  If there exists an $M>0$, an $R>0$, an
$\alpha>5$ such that \eqref{phi_cond} holds, and furthermore
\begin{equation}
\label{series_conv_cond}
M < 8 (\alpha+2)(\alpha-5)R^{\alpha-4},
\end{equation}
then the series for $f(x)$ converges for all $x$ such that $|x-d|>R$.
\begin{proof}
Combining Lemmas \ref{series_bound_lem} and \ref{series_conv_lem}, we
find that the key condition is that $A_1 \le 8R$, which by
substitution into \eqref{A_1_def} yields
\begin{equation*}
0<|K|+\frac{M}{(\alpha+2)(\alpha-5)R^{\alpha-5}}\le 8R.
\end{equation*}
But in order to have $|K| \ge 0$, this gives
\begin{equation*}
0<8R-\frac{M}{(\alpha+2)(\alpha-5)R^{\alpha-5}},
\end{equation*}
which leads immediately to the condition stated.
\end{proof}
\end{thm}

\begin{cor}
\label{finite_action_eq_cor}
As an immediate consequence of Theorem \ref{series_conv_thm}, the
action functional $A$ given by \eqref{action_eqn} evaluated at each
equilibrium is finite provided the conditions \eqref{phi_cond} and 
\eqref{series_conv_cond} on $\phi$ hold.
\end{cor}

\begin{rem}
It is worth noting that if the spatial dimension $n>1$
(\eqref{long_ode1} is now an elliptic {\it partial} differential
equation), then the asymptotic decay rate will typically be slower
than that of the series solution given here.  As a result, Corollary
\ref{finite_action_eq_cor} will not hold for higher spatial
dimensions.  Indeed, whether anything like Corollary \ref{finite_action_eq_cor}
holds in higher spatial dimensions is an open question.
\end{rem}

\begin{figure}
\begin{center}
\includegraphics[height=3in]{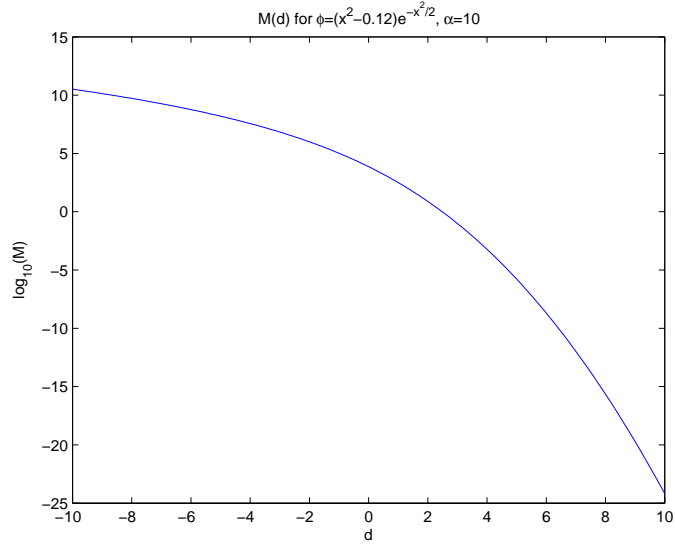}
\end{center}
\caption{A typical $M(d)$ function}
\label{m_of_d}
\end{figure}

\begin{figure}
\begin{center}
\includegraphics[height=3in]{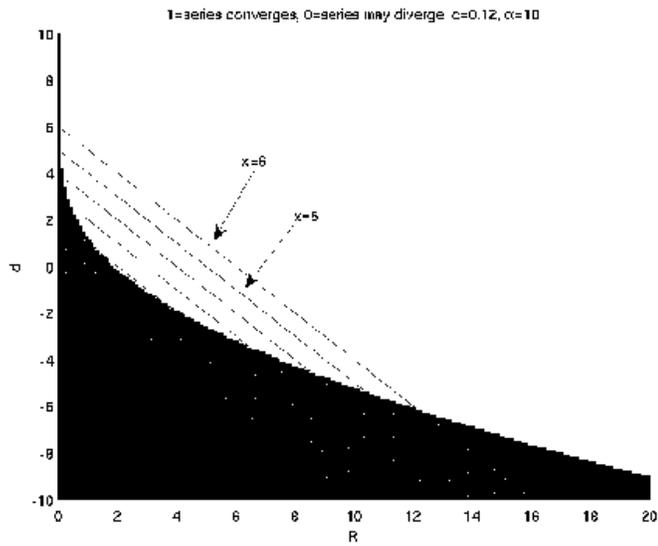}
\end{center}
\caption{Series convergence test, for $\phi(x)=(x^2-0.12)e^{-x^2/2}$:
  white = series converges, black = series may diverge }
\label{series_conv}
\end{figure}

\begin{eg}
It is important to notice that the $M$ defined above in Lemma
\ref{series_bound_lem} can depend crucially upon the value of $d$ and
the shape of the curve $\phi(x)$.  For the case of
$\phi(x)=(x^2-c)e^{-x^2/2}$, a typical plot of $M(d)$ is shown in Figure
\ref{m_of_d}.  It should be noted that for various values of $c$, the
$M(d)$ function is numerically very similar.

This also means that the condition \eqref{series_conv_cond} defines a
somewhat complicated region over which parameters $d, K$ and $R$ yield
convergent series solutions.  An example with our given $\phi(x)$
function is shown in Figure \ref{series_conv}.  Thus it appears that
our series solution converges if one goes out far enough, and
specifies small enough initial conditions.
\end{eg}

\begin{rem}
The convergence of the series solution is controlled by the
convergence of a well-behaved power series.  It follows that as
the $\phi$ function becomes smaller, fewer terms in the series are
needed to accurately approximate the solution.  Indeed, each term in
the series solution is asymptotically smaller than the previous one.
Thus, we can gain some qualitative information from the leading two
terms of the series, which are
\begin{equation*}
f(x) \sim \frac{6}{(x-d)^2}+\frac{1}{(x-d)^3}\left[ K + \int^x(t-d)^6\int_t^\infty
  \frac{\phi(s)}{(s-d)^3}ds\,dt\right].
\end{equation*}
Taking a derivative by $x$ gives
\begin{equation*}
\begin{split}
f'(x) \sim \frac{-12}{(x-d)^3}+\frac{-3}{(x-d)^4}\left[ K + \int^x(t-d)^6\int_t^\infty
  \frac{\phi(s)}{(s-d)^3}ds\,dt\right]+ \\(x-d)^3\int_x^\infty
  \frac{\phi(s)}{(s-d)^3}ds.
\end{split}
\end{equation*}
On the other hand, using the standard expansion for $(a+b)^{3/2}$, one
obtains
\begin{eqnarray*}
f^{3/2}(x) &\sim& \left( \frac{6^3}{(x-d)^6}+\frac{3\cdot
  6^2}{(x-d)^7}\left[ K + \int^x(t-d)^6\int_t^\infty
  \frac{\phi(s)}{(s-d)^3}ds\,dt\right] \right)^{1/2}\\
&\sim&\frac{6^{3/2}}{(x-d)^3}+\frac{(x-d)^3}{2\cdot 6^{3/2}}\frac{3\cdot
  6^2}{(x-d)^7}\left[ K + \int^x(t-d)^6\int_t^\infty
  \frac{\phi(s)}{(s-d)^3}ds\,dt\right]\\
&\sim&\frac{6^{3/2}}{(x-d)^3}+\frac{3\cdot
  18}{6^{3/2}(x-d)^4}\left[ K + \int^x(t-d)^6\int_t^\infty
  \frac{\phi(s)}{(s-d)^3}ds\,dt\right]\\
\end{eqnarray*}
which leads to 
\begin{equation}
\label{handy_asymp_exp}
f'(x) \sim -\sqrt{\frac{2}{3}} f^{3/2} + (x-d)^3\int_x^\infty
  \frac{\phi(s)}{(s-d)^3}ds.
\end{equation}
Notice that this equation depends only on $d$, not $K$.  So from this
we should expect that the initial data for solutions to be confined to
a thin region in the plane $x=0$.  This will be confirmed in Theorem
\ref{properties_of_Z}

Additionally, the relation $f'=-\sqrt{2/3} f^{3/2}$ holds exactly for
the bounded solutions of $0=f''-f^2$.  Indeed, in that case, the set
$Z$ is $\{(f,f')|3f'^2=2f^3,f'<0\}$.  So \eqref{handy_asymp_exp}
indicates that the presence of $\phi \neq 0$ will deflect the set $Z$
largely in the $f'$ direction.  This is exactly what we show in Section
\ref{geom_props_Z}.
\end{rem}

\section{Restriction to $\phi$ nonnegative and monotonically
  decreasing}
\label{extension_sec}

We now examine what stronger results can be obtained by requiring
$\phi(x)\ge 0$ and $\phi'(x)<0$ for all $x>0$.  This can be expected
to provide stronger results, in particular because the region $R_1$
employed in Theorem \ref{region_r1_thm} acquires a simpler inflow and
outflow structure on the boundary, and in particular, solutions will
exist for all $x>0$.  A collection of four results indicate that all
bounded solutions to \eqref{ode1} lie within a narrow region.

\begin{lem}
\label{region_r1_lem}
Suppose $\phi(x)\ge 0$ and $\phi'(x) < 0$ for all $x \ge 0$.  Then the
region given by $R_1=\{(f,f',x)|H(f,f',x) \ge 0, x \ge 0, f\le
\sqrt{\phi(x)}\}$ contains a bounded solution to \eqref{ode1}.
\begin{proof}
Following the proof of Theorem \ref{region_r1_thm}, we partition the
boundary of $R_1$ into two pieces: $A=\{(f,f',x)|x = 0\}$ and
$B=\{(f,f',x)|H(f,f',x)=0\}$, noting that the flow of $V$ is inward
along $A$.  Reviewing the computation in Theorem \ref{region_r1_thm}, 
the flow is outward along all of $B$.

Now we employ the Antifunnel theorem, noting that while $A$ is
simply-connected, $B$ is not.  Hence they cannot be homeomorphic, and
so there must be a solution that remains inside $R_1$ (which evidently
starts on $A$).  But the first coordinate of such an integral curve
must obviously be bounded, since the $x$ cross-sections of $R_1$ form
a decreasing sequence of sets, ordered by inclusion, and the
cross-section for $x=0$ is a bounded set.
\end{proof}
\end{lem}

\begin{lem}
\label{region_r2_lem}
Suppose $\phi(x)\ge 0$ and $\phi'(x) < 0$ for all $x \ge 0$.  Then the
region given by $R_2=\{(f,f',x)|H(f,f',x) \le 0, \frac{1}{3}f^3 -
\frac{1}{2} f'^2 \ge 0, x \ge 0,f'\le 0\}$ contains a bounded solution
to \eqref{ode1}.
\begin{proof}
Partition the boundary of $R_2$ into two pieces:
\begin{equation*}
A=\{(f,f',0)|f' \le 0\}\cup\{(f,f',x)|H(f,f',x)=0,f\le
\sqrt{\phi(x)},f'\le 0\},
\end{equation*}
and
\begin{equation*}
\begin{split}
B=\{(f,f',x)|H(f,f',x)=0,f\ge \sqrt{\phi(x)},f'\le 0 \}\cup \\ \{(f,f',x)|
\frac{1}{3}f^3 - \frac{1}{2} f'^2=0, f'\le 0 \}.
\end{split}
\end{equation*}
By the calculation in Theorem \ref{region_r1_thm}, the flow along $A$ is
inward-going.  Additionally, the flow along the first connected
component of $B$ is outward-going.  Finally, we put
$S(f,f',x)=\frac{1}{3}f^3 - \frac{1}{2} f'^2$ and observe that $\nabla
S$ is an inward pointing normal vector field to $B$.  We compute
\begin{eqnarray*}
\nabla S \cdot V &=& \begin{pmatrix} f^2 \\ - f' \\ 0 \end{pmatrix}^T 
\begin{pmatrix} f' \\ f^2 - \phi(x) \\ 1 \end{pmatrix} \\
&=& f' \phi(x) \le 0,
\end{eqnarray*}
so the flow along this component of $B$ is outward-going.  As a
result, we can apply the Antifunnel theorem, noting that $A$ is
connected, while $B$ is not.  Therefore, there exists a solution to
\eqref{ode1} that remains in $R_2$.  Note that there is a lower bound
on the $x$-coordinate of this solution, since the $x$-component of
$V(f,f',x)$ is equal to 1, and the Region $R_2$ lies within the
half-space $x>0$.  So this solution must enter $R_2$ through $A$, and
then never intersect $B$.  Additionally, notice that such a solution
will have $f' \le 0$ and $f \ge 0$, so it must be bounded.
\end{proof}
\end{lem}

\begin{lem}
\label{unbounded_funnel_lem}
Suppose $\phi(x)\ge 0$ and $\phi'(x) < 0$ for all $x \ge 0$.  The
complement of the set $A=R_1 \cup R_2$ consists of solutions which are
unbounded, and blow up in finite $x$. 
\begin{proof}
Let the complement of the set $A$ be called $C$, namely
$C=\{(f,f',x)|x>0\} - A$.  Now the calculations in Lemmas
\ref{region_r1_lem} and \ref{region_r2_lem} show that $C$ is a funnel,
in that the flow through the entire boundary of $C$ is inward.  If
$\phi$ does not tend to zero, then the argument in the proof of
Theorem \ref{limzero_iff_bounded_lem} completes the proof, as there is
a tubular neighborhood about $\{f=f'=0\}$ with strictly positive
radius in which solutions in $C$ cannot remain.  So without loss of
generality, we assume $\phi \to 0$.

\begin{figure}
\begin{center}
\includegraphics[height=3in]{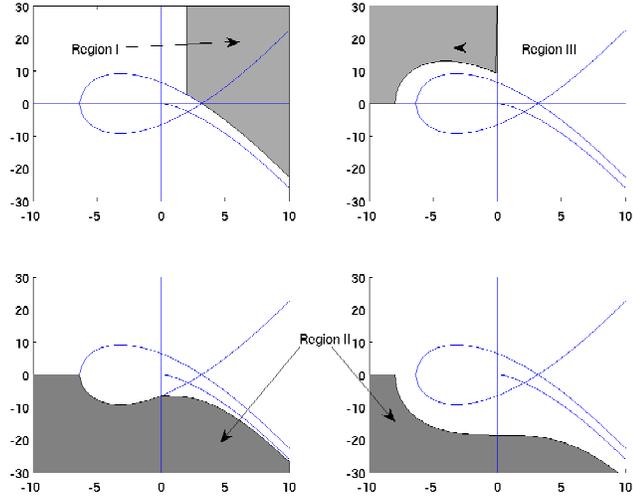}
\end{center}
\caption{The Regions $I$, $II$, and $III$ of Lemma \ref{unbounded_funnel_lem}}
\label{unboundedfunnel_fig}
\end{figure}

Define the Region $I$ by
\begin{equation*}
I=\left\{(f,f',x)| f > \sqrt{\phi(x)} \text{ and } \left(f'>0 \text{ or
  } H(f,f',x)>0\right) \right\}.
\end{equation*}
There are two bounding faces of Region $I$, along which the flow is
inward.  The first is $S_1 = f-\sqrt{\phi(x)} = 0$, along which  
\begin{eqnarray*}
\nabla S_1 \cdot V(f,f',x) &=& \begin{pmatrix}1\\ 0 \\ -
  \frac{\phi'(x)}{2 \sqrt{\phi(x)}} \end{pmatrix}^T
  \begin{pmatrix}f'\\ f^2 - \phi(x) \\ 1\end{pmatrix}\\
&=& f' - \frac{\phi'(x)}{2 \sqrt{\phi(x)}} > 0.
\end{eqnarray*}
The second was computed already in the proof of Theorem
\ref{region_r1_thm}.  Notice that $f''=f^2 - \phi(x) > 0 $ in Region
$I$, so $f(x)$ is concave-up, so solutions which enter Region $I$ are
unbounded.  Using similar reasoning to that of Theorem
\ref{limzero_iff_bounded_lem}, such solutions blow up in finite $x$.

Now suppose we have a point $(a,a',x_0) \in C$ with $a' < 0$.  We claim
that for some $x_1 > x_0$, the integral curve through this point will
cross the $f'=0$ plane.  To see this, construct Region $II$ by
\begin{equation*}
II=\left\{(f,f',x)|\frac{1}{3}f^3 - \frac{1}{2}f'^2 - \frac{1}{3}a^3 +
\frac{1}{2} a'^2 \le 0 \text{ and } f' \le 0\right\} \cap C.
\end{equation*}
Note that
\begin{equation*}
\nabla\left(\frac{1}{3}f^3 - \frac{1}{2} f'^2\right)\cdot V(f,f',x) =
f'\phi(x) \le 0,
\end{equation*}
so the flow is inward along Region $II$ except along $f'=0$ (along
which it is outward).  Also note that Region $II$ excludes a tubular
neighborhood of the line $f=f'=0$ with strictly positive radius.  As a
result of this, the integral curve through $(a,a',x_0)$ proceeds at
least as far as to allow $f < - \sqrt{\phi(x)}$, at which point, a
finite amount of distance in $x$ takes it to $f'=0$.

So at that point, the integral curve has entered Region $III$, say at
$x=x_1$, where
\begin{equation*}
III=\{(f,f',x)| H(f,f',x_1) \le 0 \text{ and } f \le 0 \text{ and } f'
\ge 0\}.
\end{equation*}
The flow is evidently inward along $f'=0$ and the curved portion by
previous calculations, and outward along $f=0$.  Again, note that the
line $f=f'=0$ is excluded from Region $III$ by a tubular neighborhood
of strictly positive radius, so there is an $x_2 > x_1$ where the
integral curve exits Region $III$ through $f=0$.

Now, consider a point $(0,c',x_2)$ along this integral curve with
$c'>0$.  In this case, the flow moves such a point rightward.  On the
other hand, the left boundary of Region $I$ moves leftward,
approaching $f=0$.  So there must be an $x_3 > x_2$ such that the
integral curve through $(0,c',x_2)$ enters the Region $I$.  Collecting
our findings, we see that every point in $C$ has an integral curve
which passes to Region $I$, and therefore corresponds to a solution
which is unbounded, and blows up for some finite $x$.  
\end{proof}
\end{lem}

\begin{thm}
\label{properties_of_Z}
Suppose $\phi(x)\ge 0$ and $\phi'(x) < 0$ for all $x \ge 0$.  The set $Z$
of initial conditions to \eqref{ode1} that lead to bounded solutions
\begin{enumerate}
\item lies within $A=R_1 \cup R_2$ and is
\item nonempty,
\item closed,
\item unbounded,
\item connected, and
\item simply connected.
\item Additionally, the portion of $Z$ corresponding to solutions that enter
the interior of $R_2$ is a 1-dimensional submanifold of $\{(f,f',x)|x=0\}$.
\end{enumerate}
\begin{proof}
\begin{enumerate}
\item From Lemma \ref{unbounded_funnel_lem}, all bounded solutions
  must lie in $A$. 
\item That there exist bounded solutions in $A$ is the content of
  Lemmas \ref{region_r1_lem} and \ref{region_r2_lem}.

\item Now, put $A_0 = A \cap \{(f,f',0)\}$ and $B_0 = \partial A -
  A_0$.  Observe that from the proofs of the previous theorems, the
  flow of $V$ along $A_0$ is inward, and the flow along $B_0$ is
  outward.  Since the last component of $V$ does not vanish, the flow
  of $V$ causes each point of $B_0$ to lie on an integral curve
  starting on $A_0$.  This establishes a homeomorphism $\Omega$ from
  $B_0$ into a subset of $A_0$.  In particular, $\Omega$ is an open
  map.  Now every solution passing through $B_0$ is of course
  unbounded, so $Z= A_0 - \Omega(B_0)$ is evidently closed (it is the
  complement of an open set).

\item $B_0$ clearly has the topology of $\mathbb{R} \times
  [0,\infty)$, so $\pi_1(B_0) = 0$.  Hence, $\pi_1(\Omega(B_0))=0$
  also, but notice that $\Omega(B_0)$ contains $\partial A_0$.  Suppose
  $Z$ were a bounded set.  Then it is contained in some disk $D$.  But
  $\partial D$ is homotopic to a loop in $A_0 - Z$, which either lies
  in $\text{int}(A_0 - Z)$ (in which case the homotopy need not move
  it) or in $\partial A_0$.  But this means that the loop encloses all
  of $Z$, and so cannot be contractible in $\Omega(B_0)$, which
  contradicts the triviality of $\pi_1(\Omega(B_0))$.   Hence $Z$ is
  unbounded.

\item We first show that the portion of $Z$ lying in the region $R_2$
satisfies the horizontal line test. First, note that a solution
starting in $Z\cap R_2$ cannot exit $R_2$.  For one, it cannot enter
$R_1$, since $R_1$ is an antifunnel.  Secondly, it cannot exit into
$\mathbb{R}^3 - (R_1 \cup R_2)$ since solutions there are all
nonglobal.  Suppose that $f_1(0) \ge f_2(0) \ge 0$ and
$f_1'(0)=f_2'(0)$ with $(f_1(0),f_1'(0))$ and $(f_2(0),f_2'(0))$ both
in $Z \cap R_2$.  But then
\begin{eqnarray*}
\frac{d}{dx}(f_1'(x)-f_2'(x))&=&f_1''(x)-f_2''(x)\\
&=&f_1^2(x) - f_2^2(x) \ge 0,\\
\end{eqnarray*}
with equality only if $f_1(0)=f_2(0)$.  Hence,
$\frac{d}{dx}(f_1(x)-f_2(x)) \ge 0$ for $x>0$, again with equality
only if $f_1(0)=f_2(0)$.  Now all solutions which remain in $R_2$ are
monotonic decreasing and bounded from below, so they must have
limits.  On the other hand, the only possible limit is
$(0,\lim_{x\to\infty}\sqrt{\phi(x)})$, so therefore all bounded
solutions in $R_2$ must have a common limit.  Therefore, we must
have that $f_1(0)=f_2(0)$.  Now this means that the portion of $Z$ in
the region $R_2$ can be realized as the graph of a function from the $f'$
coordinate to the $f$ coordinate.  Therefore, if $Z$ were not
connected, at least one component of $Z$ would be a bounded subset, which is
a contradiction.

\item Finally, if $Z$ were not simply connected, the Jordan curve
  theorem gives that there are two (or more) path components to
  $\Omega(B_0)=A_0 - Z$, which contradicts the continuity of $\Omega$.

\item By the connectedness of $Z$ and the horizontal line test in
  $R_2$, the function from the $f'$ coordinate to the $f$ coordinate
  whose graph is $Z\cap\text{int }R_2$ must be continuous.
  Additionally, by the connectedness of $Z$ and the uniquenss of
  solutions to ODE, this implies that the rest of $Z$ whose solutions
  enter the interior of $R_2$ is also a 1-manifold.
\end{enumerate}
\end{proof}
\end{thm}

\begin{df}
It is convenient to define, in addition to the initial condition set
$Z$, other sets $Z_{x_0} \subset \{(f,f',x)|x=x_0\}$ such that any
integral curve passing through a point in $Z_{x_0}$ exists for all
$x>0$.  Similarly, one can define $Z_{x_0}'$.
\end{df}

\begin{rem}
\label{approx_Z_rem}
If $\phi\to 0$ as $x\to\infty$, we conjecture that $Z$ acquires the
structure of a 1-manifold with boundary.  The series solution
\eqref{series_soln_coef} is not valid at such a boundary of $Z$, since
such a solution must remain in $R_1$ and therefore decays quicker than
the leading coefficient of \eqref{series_soln_coef}.  Indeed, by
analogy with the case where $\phi \equiv 0$, the leading term $f_0$ of
the series solution would vanish, and the solution is then asymptotic
to $-\int_x^\infty \int_t^\infty \phi(s) ds\, dt$.

All solutions in the form of the series solution
\eqref{series_soln_coef} enter $R_2$, so a result of this theorem is
that one of the two parameters $d$ or $K$ in the series solution is
superfluous.  Since $d$ parametrizes solutions when $\phi\equiv 0$, we
conventionally take $K=0$.  Using this, \eqref{handy_asymp_exp}
indicates that a good approximation (as $x_0 \to \infty$, locally near
$f=f'=0$) to the set $Z_{x_0}$ is the set
\begin{equation*}
\{H(f,f')=0\}=\{(f,f')|\frac{1}{3}f^3=\frac{1}{2}f'^2\}.
\end{equation*}
\end{rem}

\begin{rem}
If $\phi \to P >0$ as $x\to\infty$, then it is not true that $Z$ is a
1-manifold (with boundary).  Indeed, $Z$ has the structure of a
1-manifold attached to the teardrop-shaped set $M$ from Lemma
\ref{bounded_in_funnel_lem}.
\end{rem}

\section{Geometric properties of the initial condition set $Z$}
\label{geom_props_Z}

\begin{lem}
\label{Z_intersects_f_axis}
Suppose $\phi(x)>0$, $\phi'(x)<0$ for all $x>0$ and $\phi \to 0$ as
$x \to \infty$.  Then the set $Z$ intersects $\{(f,f',x)|f'=0\}$.
\begin{proof}
First, observe that $Z$ intersects the boundary of $R_1$ in $x=0$, since we
have by Lemmas \ref{region_r1_lem} and \ref{region_r2_lem} solutions
entirely within $R_1$ and its complement.  Using the fact that $Z$ is
connected and the Jordan curve theorem, $Z$ must intersect the
boundary of $R_1$ in the plane $x=0$.  This reasoning also applies for
each $Z_{x_0}$ with $x_0\ge 0$, so that we can find points in the
intersections $Z_{x_0} \cap \partial R_1$ for each $x_0 \ge 0$.  Also
note that for the backwards flow associated to our equation (ie. the
flow of $-V$), solutions which enter $R_1$ must exit through the plane
$x=0$.  Hence there exists a sequence of points $\{F_n\} \subset Z$
with $F_n=(f_n,f_n',0)$ such that the integral curve through $F_n$
passes through $G_n=(g_n,g_n',n) \in Z_n \cap \partial R_1$ for each
integer $n \ge 0$.

Discern three cases:
\begin{enumerate}
\item If any $F_n$ are in Quadrants I or II, then since $Z$ is
  connected, it must intersect $\{f'=0\}$.
\item If any $F_n$ are in Quadrant III, observe that the flow across
  the surface $S=\left \{(f,f',x)|\frac{1}{3}f^3=\frac{1}{2}f'^2, f' \le
  0\right \}$ is right-to-left.  Thus the integral curve must cross
  into Quadrant II on its way to $G_n$.  Therefore, the set $Z$
  cannot intersect the surface $S$, and so it must intersect
  $\{f'=0\}$.
\item Assume all the $F_n$ lie in Quadrant IV.  Observe that $\{F_n\}$
  is a closed subset of $R_1 \cap \{x=0\}$, which is compact.
  Hence some subsequence of $\{F_n\}$ must have a limit, say $F$.
  Since $Z$ is closed, $F \in Z$.  But in the portion of $R_1$ lying
  in the $x=0$ plane and in Quadrant IV, we have that 
\begin{equation*}
\frac{d}{dx}f'=f^2 - \phi < 0
\end{equation*}
and
\begin{equation*}
\frac{d}{dx}f=f' < 0.
\end{equation*}
Hence $f_n' \ge g_n'$.  But since $\phi \to 0$, $g_n' \to 0$, so $F$
lies on $\{f'=0\}$.
\end{enumerate}
\end{proof}
\end{lem}

\begin{lem}
\label{Z_intersects_fp_axis_weak}
Under the same hypotheses as Lemma \ref{Z_intersects_f_axis}, $Z$ also
intersects the half plane $\{f=0,f'>0\}$.
\begin{proof}
Using Lemma \ref{Z_intersects_f_axis}, we form a sequence $\{F_n\}
\subset Z$ such that the integral curve through $F_n$ passes through
$\{f'=0,f \ge 0, x=n\}$ for each integer $n$.  (This can be done
without loss of generality, because if any integral curves pass
through $\{f'=0,f < 0\}$, then the proof is complete by connectedness
of $Z$.)  Note that this sequence is entirely contained within $R_1$
by Lemma \ref{unbounded_funnel_lem}.

Discern three cases:
\begin{enumerate}
\item There exists an $F_n$ in either of Quadrants II or III.  The
  result follows by the connectedness of $Z$.
\item There exists $F_n$ in Quadrant IV.  This cannot occur unless the
  integral curve through $F_n$ passes through Quadrant III since the
  flow along $\{f'=0\}$ points inward into the portion of Quadrant IV
  inside $R_1$.
\item Otherwise, we assume $\{F_n\}$ is entirely contained within
  Quadrant I.  In this case, note that 
\begin{equation*}
\frac{d}{dx} f' = f^2 - \phi < 0.
\end{equation*}
Hence the $f'$-coordinate of the integral curve through each $F_n$ is
positive on the interior of Quadrant I.  Hence 
\begin{equation*}
\frac{d}{dx} f = f' > 0,
\end{equation*}
so $f_n \le g_n$.  But $g_n \to 0$ since $\phi \to 0$, so any limit
point of $\{F_n\}$ will have $f$-coordinate equal to zero.  By the
compactness of $R_1 \cap \{x=0\}$ and the closedness of $Z$, this
implies that $Z$ intersects $\{f=0,f'>0\}$.
\end{enumerate}
\end{proof}
\end{lem}

\begin{lem}
\label{Z_intersects_fp_axis}
Suppose $\phi(x) >0$ for all $x>0$, $\phi \to 0$ as $x \to \infty$,
and that there exists an $x_0 \ge 0$ such that for all $x>x_0$,
$\phi'(x)<0$.  Then the set $Z$ intersects $\{f=0, f'>0\}$.
\begin{proof}
We follow the pattern of proving the existence of an intersection for
an open interval in $x$ containing $x_0$, and then constructing an
{\it a priori} estimate for the $f'$-coordinate of this intersection.  

Apply Lemma \ref{Z_intersects_fp_axis_weak} to $x_0$, we have that
$Z_{x_0}$ intersects $\{f=0,f'>0\}$.  Let $(0,f_0',x_0)$ lie in this
intersection.  Note that 
\begin{equation*}
\frac{d}{dx} f = f'>0
\end{equation*}
and
\begin{equation*}
\frac{d}{dx}f' = f^2 - \phi = -\phi < 0
\end{equation*}
when evaluated there.  As a result, the integral curve passing through
$(0,f_0',x_0)$ must pass through Quadrant II first, say for $x \in
(x_1,x_0)$.  Then evidently, $Z_{x_1}$ must intersect $\{f=0,f'>0\}$.

Now since $\phi(x)>0$ between $x_1$ and $x_0$, and $[x_1,x_0]$ is
compact, there is an open set in $\mathbb{R}^3$ containing the
intersection of each $Z_x$ with $\{f=0,f'>0\}$ for each $x\in
[x_1,x_0]$, such that in this open set $\frac{d}{dx} f' \le K < 0$.  As a
result, $f_1' \ge f_0'$.  Hence the $f'$-coordinate of the
intersection point of $Z_x$ with $\{f=0,f'>0\}$ is decreasing with
increasing $x$.  (Since we have $f^2 - \phi > -\phi$, it is decreasing
at a rate no faster than $\phi$.  This implies that this intersection point
has $f'$-coordinate no larger than $\int_0^{x_0} \phi(x) dx + f_0'$ at
$x=0$.)  Now since solutions through $Z_{x_1}$ exist for all $x>0$ by
definition, this suffices to show that $Z$ intersects $\{f=0,f'>0\}$.
\end{proof}
\end{lem}

\begin{rem}
\label{delicateness_rem}
The line of reasoning used in the third case of each of Lemmas
\ref{Z_intersects_f_axis} and \ref{Z_intersects_fp_axis_weak} (and
also in \ref{Z_intersects_fp_axis}) fails if we try to continue $Z$ much
farther.  This is due to the nonmonotonicity of $df'/dx$ in Quadrants
II and III.  More delicate control of $\phi$ must be exercised to say
more.
\end{rem}

\begin{calc}
Towards the end of the more delicate results mentioned in Remark
\ref{delicateness_rem}, it is useful to know the maximum speed along
integral curves on points in the region $R_1$ in the $f$- and
$f'$-directions.  By this we mean to compute for fixed $x$ the maximum
values of 
\begin{equation}
\label{directional_speeds}
\begin{cases}
|f'| \text{ for the }f\text{-direction}\\
|f^2-\phi(x)| \text{ for the }f'\text{-direction}\\
\end{cases}
\end{equation}
in $R_1$.  The first is easy to maximize: we simply look for the
maximum value of $f'$ in $R_1$, which is a maximum of 
\begin{equation*}
f'= \sqrt{\frac{2}{3} f^3 - 2f\phi(x) + \frac{4}{3} \phi^{3/2}(x)},
\end{equation*}
for $-2\sqrt{\phi(x)} \le f \le \sqrt{\phi(x)}$.  This occurs at
$f=-\sqrt{\phi(x)}$, and has the value of $\sqrt{8/3}\phi^{3/4}$.
For the second part of \eqref{directional_speeds}, it is easy to see
that the maximum is $3 \phi(x)$. In summary,
\begin{equation}
\label{max_speed_calc}
\begin{cases}
|f'| \le \sqrt{\frac{8}{3}}\phi^{3/4}(x) \text{ for the }f\text{-direction}\\
|f^2-\phi(x)| \le 3\phi(x) \text{ for the }f'\text{-direction}\\
\end{cases}
\end{equation}
on $R_1$.
\end{calc}

Using this calculation, we can impose a stronger bound on the decay of
$\phi(x)$, and constrain the set $Z$ further.

\begin{lem}
\label{constrain_Z_lem_weak}
Suppose $\phi(x)>0$, $\phi'(x) < -D
\frac{4\sqrt{2}}{k\sqrt{3}}\phi^{5/4}(x)$ for all $x>0$ for some
$0<k<1$ and $D>1$.  Then the set $Z$ is contained within $\{f \ge
-k\sqrt{\phi(0)}\}$ and intersects each vertical and horizontal line
in $\{f \ge 0\}$ exactly once, and intersects $\{f'=0\}$ only once.
\begin{proof}
That $Z$ intersects $\{f=0,f' \ge 0\}$ and $\{f'=0,f\ge 0\}$ at all
follows from Lemmas \ref{Z_intersects_f_axis} and
\ref{Z_intersects_fp_axis}.  Now consider the region $A \subset R_1$
shown in Figure \ref{constrain_Z_fig} and defined by 
\begin{equation*}
A =R_1 \cap \left(\{f'\ge 0, f \le k \sqrt{\phi(x)}\} \cup \{f' \le
0, 2 f^3 \le 3 f'^2\} \right).
\end{equation*}

\begin{figure}
\begin{center}
\includegraphics[height=3in]{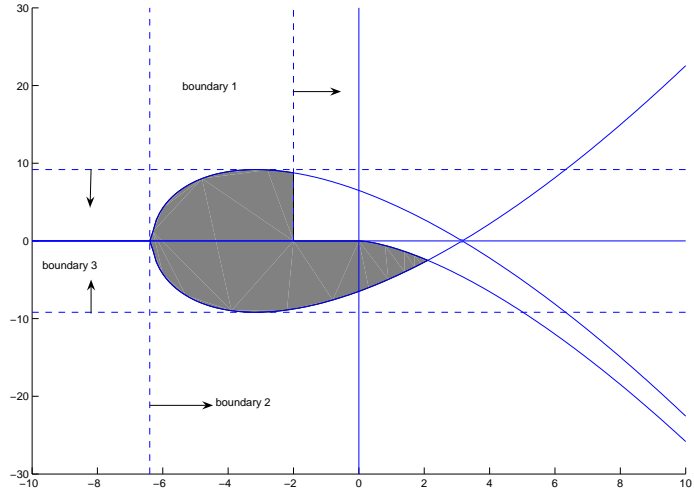}
\end{center}
\caption{The region $A$ of Lemma \ref{constrain_Z_lem_weak}}
\label{constrain_Z_fig}
\end{figure}

The boundary segments strictly to the right of the boundary labelled 1 in
Figure \ref{constrain_Z_fig} are evidently inflow, so long as $\phi >
0$.  The boundary labelled as 1 in the figure moves with speed
\begin{equation*}
\frac{d}{dx} ( - k \sqrt{\phi(x)}) = \frac{-k}{2\sqrt{\phi(x)}}
\phi'(x) > D \frac{2\sqrt{2}}{\sqrt{3}} \phi^{3/4}(x)
\end{equation*}
which is greater than maximum speed in the $f$-direction given in
\eqref{max_speed_calc}.  This implies that the boundary moves faster
than any solution inside $R_1$.  Hence it is an inflow portion of the
boundary.  On the other hand, the curved segment of the boundary to
the left has been shown to be outflow, in Lemma \ref{region_r1_lem}.

We observe that the boundary marked 2 in Figure
\ref{constrain_Z_fig} moves with speed
\begin{equation*}
\frac{d}{dx} ( - 2 k \sqrt{\phi(x)} ) = \frac{ -k}{\sqrt{\phi(x)}} \phi'(x),
\end{equation*}
which is strictly faster than the boundary marked 1 in Figure
\ref{constrain_Z_fig}, and the boundary marked 3 in Figure
\ref{constrain_Z_fig} moves with speed
\begin{equation*}
\frac{d}{dx} \left( \pm \sqrt{\frac{8}{3}} \phi^{3/4}(x) \right ) =
\pm \frac{\sqrt{3}}{\sqrt{2} \phi^{1/4}(x)} \phi'(x),
\end{equation*}
noting that $f'(-\sqrt{\phi(x)})=\pm\sqrt{8/3}\phi^{3/4}(x)$ is the
value of the maximum $f'$-coordinate of $R_1$ at a given $x$ value.
This last speed is greater than the maximum speed in the $f'$-direction given by
\eqref{max_speed_calc} since $\phi'(x) < - D 
\sqrt{6}\phi^{5/4}(x)$.  (Notice that
$\sqrt{6} < \frac{4\sqrt{2}}{k\sqrt{3}}$, since $0<k<1$.)  

Since $D>1$, this means that both the boundaries marked 2 and 3 in
Figure \ref{constrain_Z_fig} overtake any solution constrained to be
within $R_1$.  As a result, every solution within the region $A$ must
leave it within finite $x$.  But the only way to leave $A$ causes a
solution to enter $\mathbb{R}^3 - (R_1 \cup R_2)$, so every solution
which contains a point in $A$ cannot exist for all $x>0$ by Lemma
\ref{unbounded_funnel_lem}.  Therefore, $Z$ is contained within $(R_1
\cup R_2) - A$.

Now consider the region $B$ which is defined by
\begin{equation*}
B =R_1 \cap \left(\{f'\ge 0, f \le 0\} \cup \{f' \le
0, 2 f^3 \le 3 f'^2\} \right),
\end{equation*}
which is simply the region $A$, with $k$ taken to be zero.  The
portion of the boundary of $B$ lying in the $\{f=0\}$ plane is inflow.
We can therefore apply the reasoning of the vertical line test:
Suppose $(f_1,f_1',0),(f_2,f_2',0) \in Z$ with $f_1=f_2 > 0$ and $f_1'
\ge f_2'$.  Then we have both (at $x=0$)
\begin{equation*}
\frac{d}{dx}(f_1'-f_2')=f_1^2-f_2^2 =0 
\end{equation*}
and
\begin{equation*}
\frac{d}{dx}(f_1-f_2)=f_1'-f_2' \ge 0,
\end{equation*}
which gives that $f_1^2-f_2^2 \ge 0$ for some open interval about
$x=0$.  Then, $\frac{d}{dx}(f_1'-f_2')\ge 0$, which implies that
in fact $\frac{d}{dx}(f_1-f_2) \ge 0$.  However, since all bounded
solutions tend to the common limit of zero, we have that this implies
$f_1'=f_2'$ at $x=0$.  (Note that since each solution starts in $Z \cap
(R_1 - B)$, we have that neither solution can become negative, since
that would involve entering $B \subset A$ or leaving $R_1 \cup R_2$.)
This implies that there is a unique intersection of $Z$ with each vertical
line.  The same reasoning applies in the case of the horizontal line
test, as in Theorem \ref{properties_of_Z}.
\end{proof}
\end{lem}

\begin{lem}
\label{constrain_Z_lem}
Suppose $\phi(x)>0$ for all $x>0$, and that $\phi'(x) < -D 
\frac{4\sqrt{2}}{k\sqrt{3}}\phi^{5/4}(x)$ for all $x>x_0\ge 0$ for some $0<k<1$ and $D>1$.
Additionally, suppose that for all $x \in [0,x_0]$, 
\begin{equation}
\label{constrain_Z_eqn}
x_0 - x < \frac{\sqrt{\phi(x)}-k\sqrt{\phi(x_0)}}{
\sqrt{\frac{8}{3}} P^{3/4}},
\end{equation}
where $P=\max_{x \in [0,x_0]} \phi(x)$. Then the set $Z$ is contained
within $\{f \ge -\sqrt{\phi(0)}\}$ and intersects each vertical and
horizontal line in $\{f \ge 0\}$ exactly once, and intersects
$\{f'=0\}$ only once.
\begin{proof}
The set $Z_{x_0}$ is constrained to lie within the set $\{f' \ge
-k\sqrt{\phi(x_0)}\}$, by Lemma \ref{constrain_Z_lem_weak} (replacing
$x_0$ by zero).  Now using the $f$-direction part of
\eqref{max_speed_calc}, the smallest $f$-value attained in $Z_x$ is
\begin{equation*}
\int_{x_0}^x \sqrt{\frac{8}{3}}\phi^{3/4}(x) dx -
k \sqrt{\phi(x_0)}.
\end{equation*}
If $x<x_0$, we have
\begin{eqnarray*}
\int_{x_0}^x \sqrt{\frac{8}{3}}\phi^{3/4}(x) dx -
k \sqrt{\phi(x_0)} &\ge& \sqrt{\frac{8}{3}} P^{3/4} (x-x_0) - k\sqrt{\phi(x_0)}\\
&>&-\sqrt{\phi(x)},\\
\end{eqnarray*}
by \eqref{constrain_Z_eqn}.  As a result, $Z_x \subset \{f \ge
-\sqrt{\phi(x)}\}$ for each $x<x_0$.  This additionally means that in
the backwards flow, the entire portion of $Z_x$ contained in $\{f \le
0\}$ is moving away from the plane $\{f'=0\}$, which completes the
proof.
\end{proof}
\end{lem}

\begin{rem}
The condition that $\phi'(x) < C \phi^{5/4}(x)$
implies 
\begin{eqnarray*}
\phi^{-5/4} \phi'(x) &<& C\\
-\frac{1}{4} \phi^{-1/4}(x) &<& Cx + C'\\
\phi(x) &<&\frac{C'''}{\left(C''-x\right)^4}
\end{eqnarray*}
for some $C''$ and $C'''$.  Notice that this condition is satisfied
when the series solution converges by Theorem \ref{series_conv_thm}.
\end{rem}

\section{Solutions on the entire real line}
\label{long_sec}
We now combine the results for \eqref{ode1} and \eqref{ode1_backwards}
to discuss properties of the solutions to \eqref{long_ode1}.
When $\phi(x)$ is monotonically decreasing for $x \ge 0$, we have by Lemma
\ref{unbounded_funnel_lem} that the initial condition set stays within
$R_1 \cup R_2$.  In particular, $Z \subset \{f' \le \sqrt{\frac{8}{3}}
\phi^{3/4}(0)\}$.  If we relax the restriction of monotonicity, we
obtain a similar result.

\begin{lem}
\label{bound_on_fp}
If $f=f(x)$ is a bounded solution to the initial value problem
\eqref{ode1} with $\phi \in C^\infty\cap
  L^\infty(\mathbb{R})$ then $f'(0) < \sqrt{8/3}
    \|\phi\|_\infty^{3/4}$.
\begin{proof}
Since $f$ is a solution to \eqref{ode1}, then it must satisfy
\begin{equation*}
f''=f^2-\phi(x) \ge f^2 - \|\phi\|_\infty.
\end{equation*}
Now Lemma \ref{bounded_in_funnel_lem} shows that all bounded solutions
to $g''=g^2 - \|\phi\|_\infty$ lie in the closure of the set $M$ given
by
\begin{equation*}
M=\left \{(g,g')|\frac{1}{3} g^3 - \frac{1}{2} g'^2 - g
\|\phi\|_\infty + \frac{2}{3} \|\phi\|_\infty^{3/2}>0,g <
\sqrt{\|\phi\|_\infty} \right\}.
\end{equation*}

Since this set $M$ is bounded, we can find the maximum value of $f'$,
which is $f'_{max} = \sqrt{8/3}\|\phi\|_\infty^{3/4}.$
\end{proof}
\end{lem}

\begin{lem}
\label{bound_on_phi_integral_1}
Consider solutions to \eqref{long_ode1} on the real line, with $\phi
\in C_0^\infty \cap L^\infty(\mathbb{R})$.  If for some $-\infty < A <
B < \infty$,
\begin{equation*}
-\int_A^B \phi(x) dx > \sqrt{\frac{8}{3}} \left (
 (\sup_{x\in(-\infty,A]} |\phi(x)|)^{3/4} +
 (\sup_{x\in[B,\infty)} |\phi(x)|)^{3/4}  \right)
\end{equation*}
then no bounded solutions exist.
\begin{proof}
For a solution $f$, we have that $f''=f^2-\phi(x) \ge -\phi(x).$
Integrating both sides we have
\begin{equation*}
f'(B)-f'(A) \ge - \int_A^B \phi(x) dx.
\end{equation*}
By Lemma \ref{bound_on_fp}, bounded solutions on
\begin{itemize}
\item $x>B$ have $f'(B) < \sqrt{\frac{8}{3}} (\sup_{x\in(-\infty,A]}
  |\phi(x)|)^{3/4}$, and
\item on $x>A$, they have $f'(A) < (\sup_{x\in[B,\infty)}
  |\phi(x)|)^{3/4}$,
\end{itemize}
so a necessary condition for there to be a bounded solution is that
\begin{equation*}
-\int_A^B \phi(x) dx \le \sqrt{\frac{8}{3}} \left (
 (\sup_{x\in(-\infty,A]} |\phi(x)|)^{3/4} +
 (\sup_{x\in[B,\infty)} |\phi(x)|)^{3/4}  \right).
\end{equation*}
\end{proof}
\end{lem}

\begin{cor}
\label{phi_integral_bound}
A necessary condition for bounded solutions to \eqref{long_ode1} to exist if $\phi \in
C_0^\infty \cap L^\infty(\mathbb{R})$ is $\int_{-\infty}^\infty \phi(x) dx > 0$.
\begin{proof}
Suppose bounded solutions exist.  By the proof of Lemma
\ref{bound_on_phi_integral_1}, if we let
\begin{equation*}
g_n=-\int_{-n}^n \phi(x) dx,
\end{equation*}
and
\begin{equation*}
h_n=\sqrt{\frac{8}{3}}\left (
 (\sup_{x\in(-\infty,-n]} |\phi(x)|)^{3/4} +
 (\sup_{x\in[n,\infty)} |\phi(x)|)^{3/4}  \right),
\end{equation*}
then $g_n < h_n$ for each positive integer $n$.  But the continuity of
limits gives 
\begin{equation*}
- \int_{-\infty}^\infty \phi(x) dx = \lim_{n\to \infty} g_n < \lim_{n
  \to \infty} h_n = 0.
\end{equation*}
\end{proof}
\end{cor}

\begin{df}
A function $\phi \in C_0^\infty \cap L^\infty(\mathbb{R})$ will be
  called {\it M-shaped} if there exists an $x_0>0$ such that for all
  $|x|>x_0$, $\phi(x) > 0$ and
\begin{itemize}
\item $\phi$ is monotonic increasing for $x < -x_0$ and 
\item $\phi$ is monotonic decreasing for $x > x_0$.
\end{itemize}
\end{df}

\begin{thm}
\label{long_soln_exists}
Suppose $\phi$ is a positive M-shaped function, then solutions exist
to \eqref{long_ode1}. 
\begin{proof}
Observe that by Lemma \ref{Z_intersects_fp_axis}, we have that the set
$Z$ intersects $\{f=0,f'>0\}$.  Additionally, by Theorem
\ref{properties_of_Z}, we have that $Z$ also lies in $R_2$, which is
unbounded in Quadrant IV.  Likewise, the set $Z'$ (for
\eqref{ode1_backwards}) intersects $\{f=0,f'<0\}$, and becomes
unbounded in Quadrant I, so $Z \cap Z'$ must be nonempty,
and at least one point in this intersection is in the half-plane $\{x=0,f>0\}$.
\end{proof}
\end{thm}

\begin{thm}
\label{long_soln_unique}
Suppose $\phi$ is a positive M-shaped function which additionally
satisfies the decay constraints of Lemma \ref{constrain_Z_lem} for
$x>0$ and $x<0$ separately, then a unique positive solution exists
to \eqref{long_ode1}.  (Note that for $x<0$, the inequalities and
signs in Lemma \ref{constrain_Z_lem} must be reversed, {\it mutatis
  mutandis}.)
\begin{proof}
By the Theorem \ref{long_soln_exists}, there exist solutions to
\eqref{long_ode1}, one of which comes from the intersection of $Z \cap
Z'$ in the half-plane $\{x=0,f>0\}$.  The vertical-line test in Lemma
\ref{constrain_Z_lem} allows one to conclude that the solution which
passes through that half-plane must continue directly to the region
$R_2$ of Lemma \ref{region_r2_lem}, without crossing the plane
$\{f=0\}$.  Thus this solution is strictly positive.

On the other hand, Lemma \ref{constrain_Z_lem} indicates that $Z$ may
lie only in Quadrants I, II, and IV, while the set $Z'$ must lie in
Quadrants I, III, and IV.  On the other hand, the vertical- and
horizontal-line tests ensure a unique intersection of $Z$ and $Z'$ in
Quadrants I and IV, so the solution is unique.
\end{proof}
\end{thm}

\begin{eg}
\label{example_1}
We examine the family $\phi_c(x)=c e^{-x^2/2}$, which is M-shaped when
$c\ge 0$.  Notice that when $c<0$, then the necessary condition of
Corollary \ref{phi_integral_bound} is not met, so solutions do not
exist for all $x \in \mathbb{R}$.  When $c=0$, then the trivial
solution $f=0$ is the only solution.  For $c>0$, we examine
$\phi'_c(x) = -xce^{-x^2/2}$.  Figure \ref{Z_example_1} shows the sets
$Z$ and $Z'$ for the case when $c=0.05$.  In particular, one notes
that there appears to be a unique point of intersection.

\begin{figure}
\begin{center}
\includegraphics[height=3in]{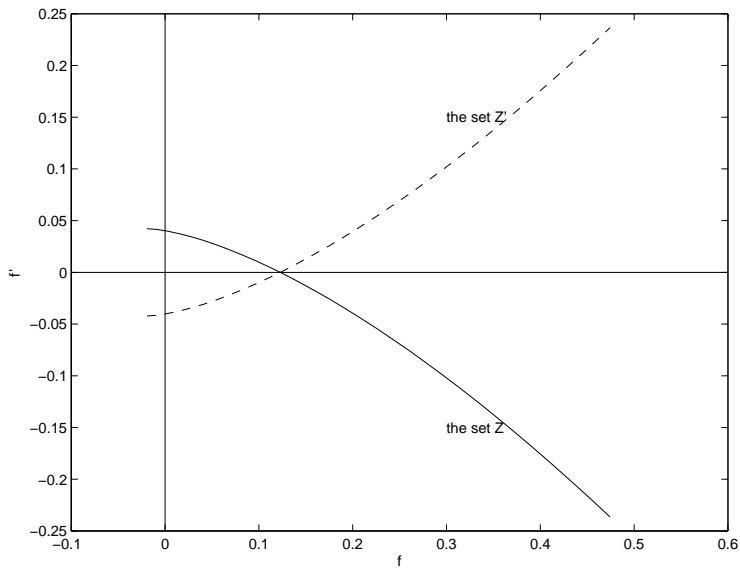}
\end{center}
\caption{The sets $Z$ and $Z'$ in Example \ref{example_1}}
\label{Z_example_1}
\end{figure}

We find the $x_0$ for which larger $x$ satisfy $\phi'(x) <
-4\sqrt{2}\phi^{5/4}(x)/(k\sqrt{3})$:
\begin{eqnarray*}
-xce^{-x^2/2} &<& -\frac{4\sqrt{2}}{k\sqrt{3}}c^{5/4}e^{-5x^2/8}\\
xe^{x^2/8}&>&\frac{4\sqrt{2}}{k\sqrt{3}}c^{1/4},
\end{eqnarray*}
which occurs if $x>\frac{4\sqrt{2}}{k\sqrt{3}}c^{1/4}$, so we may take
$x_0=\frac{4\sqrt{2}}{k\sqrt{3}}c^{1/4}$.  

By way of example, if we fix $x_0 = 4/3$, then $k=\sqrt{6}c^{1/4}$.
(We enforce $0<k<1$ by taking $c$ small.)  Now we must check to see if
\eqref{constrain_Z_eqn} holds.  In this case, we need to see if $c$
can be chosen so that $x_0-x=4/3-x$ is bounded above by
\begin{eqnarray*}
\frac{\sqrt{\phi(x)}-k\sqrt{\phi(x_0)}}{
\sqrt{\frac{8}{3}} P^{3/4}} &=&
\frac{\sqrt{c}e^{-x^2/4}-\sqrt{6} c^{3/4}
  e^{-16/36}}{\sqrt{8/3}c^{3/4}}\\
&=&\frac{e^{-x^2/4}-\sqrt{6}c^{1/4}e^{-16/36}}{\sqrt{8/3}c^{1/4}}\\
&\ge&\frac{e^{-16/36}-\sqrt{6}c^{1/4}e^{-16/36}}{\sqrt{8/3}c^{1/4}},
\end{eqnarray*} 
which can be made as large as one likes by taking $c$ sufficiently
small.  Noting that this last line is a constant in $x$ completes the
bound.  Therefore, there is a unique positive solution for
$0=f''-f^2+ce^{-x^2/2}$ with $c \in [0,\epsilon)$ for some
  $\epsilon>0$.

\end{eg}

\begin{rem}
\label{asymp_summary_rem}
Taken together, the results of Corollary \ref{phi_integral_bound} and
Theorems \ref{long_soln_exists} and \ref{long_soln_unique} for
M-shaped $\phi$ provide the following story about solutions to the
equation $0=f''-f^2+\phi$ on the real line:
\begin{itemize}
\item If the portion of $\phi$ where it is allowed to be negative is
  sufficiently negative, then no solutions exist,
\item If $\phi$ is positive, then a solution will exist.  There is no
  particular reason to believe that this solution will be strictly
  positive or unique.  
\item If the decay in the monotonic portions of $\phi$ is fast enough,
  there is exactly one solution, which is strictly positive.
\end{itemize}
\end{rem}

\section{Numerical examination}
\label{numer_sec}
\subsection{Computational framework}

Notice that the results of Remark \ref{asymp_summary_rem} are not
sharp: nothing is said if $\phi$ has a portion which is negative, but
still satisfies the necessary condition of Corollary
\ref{phi_integral_bound}.  Further, if $\phi$ is positive, but does
not satisfy the decay rate conditions, nothing is said about the
number of global solutions that exist.  Answers to these questions can
be obtained by combining the asymptotic information we have collected
about the sets $Z$ and $Z'$ with a numerical solver.  In particular,
we can obtain information about the number of global solutions to
\eqref{long_ode1} for any M-shaped $\phi$.

Implicit in the use of a numerical solver is the following Conjecture:

\begin{conj}
\label{finitely_many_equilib_conj}
There are only finitely many smooth global solutions to \eqref{long_ode1}.
\end{conj}

Suppose that $\phi$ is an M-shaped function, and that $x_0$ is such
that $\phi(x)$ is monotonic decreasing for all $x>x_0$ and is
monotonic increasing for all $x<-x_0$.  (If $\phi$ decreases fast
enough, we can choose $x_0$ so that the series solution converges on
the complement of $(-x_0,x_0)$ for sufficiently small initial
conditions.)  Then we have the sets $Z'_{-x_0}$ and $Z_{x_0}$ of
initial conditions to ensure existence of solutions on
$(-\infty,-x_0]$ and $[x_0,\infty)$ respectively.  Then any solution
to the boundary value problem
\begin{equation}
\label{ode1_bvp}
\begin{cases}
0=f''(x)-f^2(x)+\phi(x) \text{ for } -x_0<x<x_0\\
(f(-x_0),f'(-x_0))\in Z'_{-x_0}\\
(f(x_0),f'(x_0))\in Z_{x_0}\\
\end{cases}
\end{equation}
extends to a global solution of \eqref{long_ode1}.  So all one must do
is solve \eqref{ode1_bvp} numerically.  An easy way to do this is to
numerically extend the sets $Z'_{-x_0}$ and $Z_{x_0}$ to $Z'$ and $Z$
respectively (ie. extend them to the plane $x=0$) and compute $Z' \cap
Z$.  

\begin{figure}
\begin{center}
\includegraphics[height=3in]{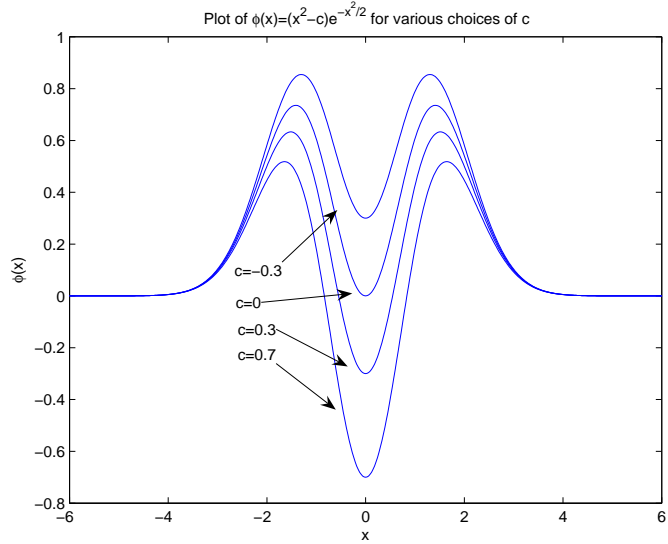}
\end{center}
\caption{The function $\phi(x;c)$ for various $c$ values}
\label{sample_phi_fig}
\end{figure}

In order to analyze \eqref{long_ode1} numerically, it is
necessary to make a choice of $\phi$.  Evidently, the numerical
results for that particular choice of $\phi$ cannot be expected to
apply in general.  However, a good choice of $\phi$ will suggest
features in the solutions that are common to a larger class of
$\phi$.  We shall use
\begin{equation}
\label{sample_phi}
\phi(x;c)=(x^2-c)e^{-x^2/2},
\end{equation}
where $c$ is taken to be a fixed parameter.  (See Figure
\ref{sample_phi_fig})  This choice of $\phi$ has the following
features which make for interesting behavior in solutions to
\eqref{long_ode1}:
\begin{itemize}
\item $\phi(x;c)>0$ for $c<0$.  In this case, there are solutions to
  \eqref{long_ode1}, by Theorem \ref{long_soln_exists}.  On the other
  hand, the decay rate conditions are not met over all of $\mathbb{R}$
  so the uniqueness result of Theorem \ref{long_soln_unique} does not
  apply.  Inded, the decay rate conditions are met only for
  sufficiently large $|x|$, but not for $|x|$ small.  
\item If $c>0$ is large enough, it should happen that no solutions to
  \eqref{long_ode1} exist, since the necessary condition of Corollary
  \ref{phi_integral_bound} is not met.  Indeed, the integral of $\phi$
  vanishes when $c=1$.
\end{itemize}

\subsection{Bifurcations in the global solutions}

\begin{figure}
\begin{center}
\includegraphics[height=2in]{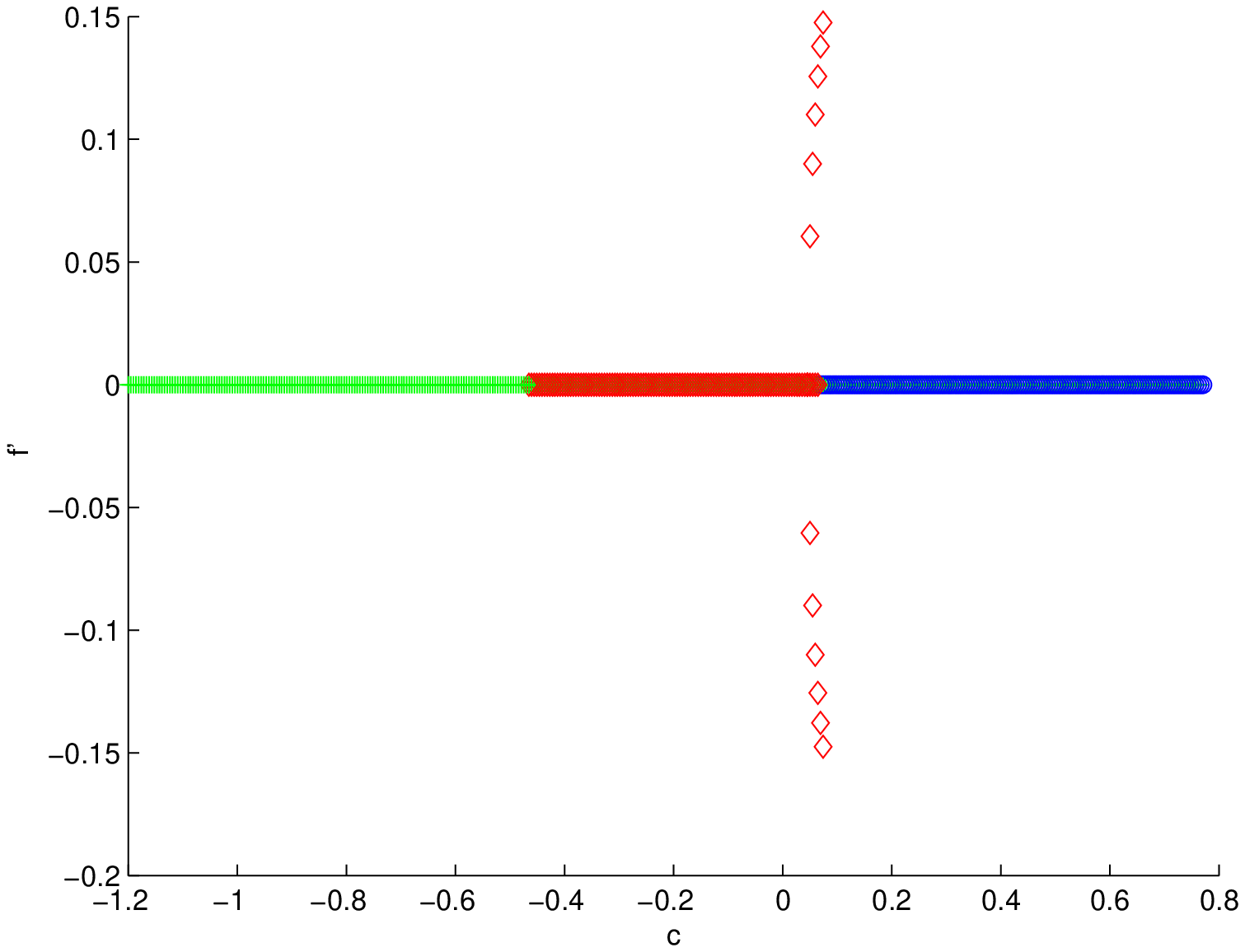}
\includegraphics[height=2in]{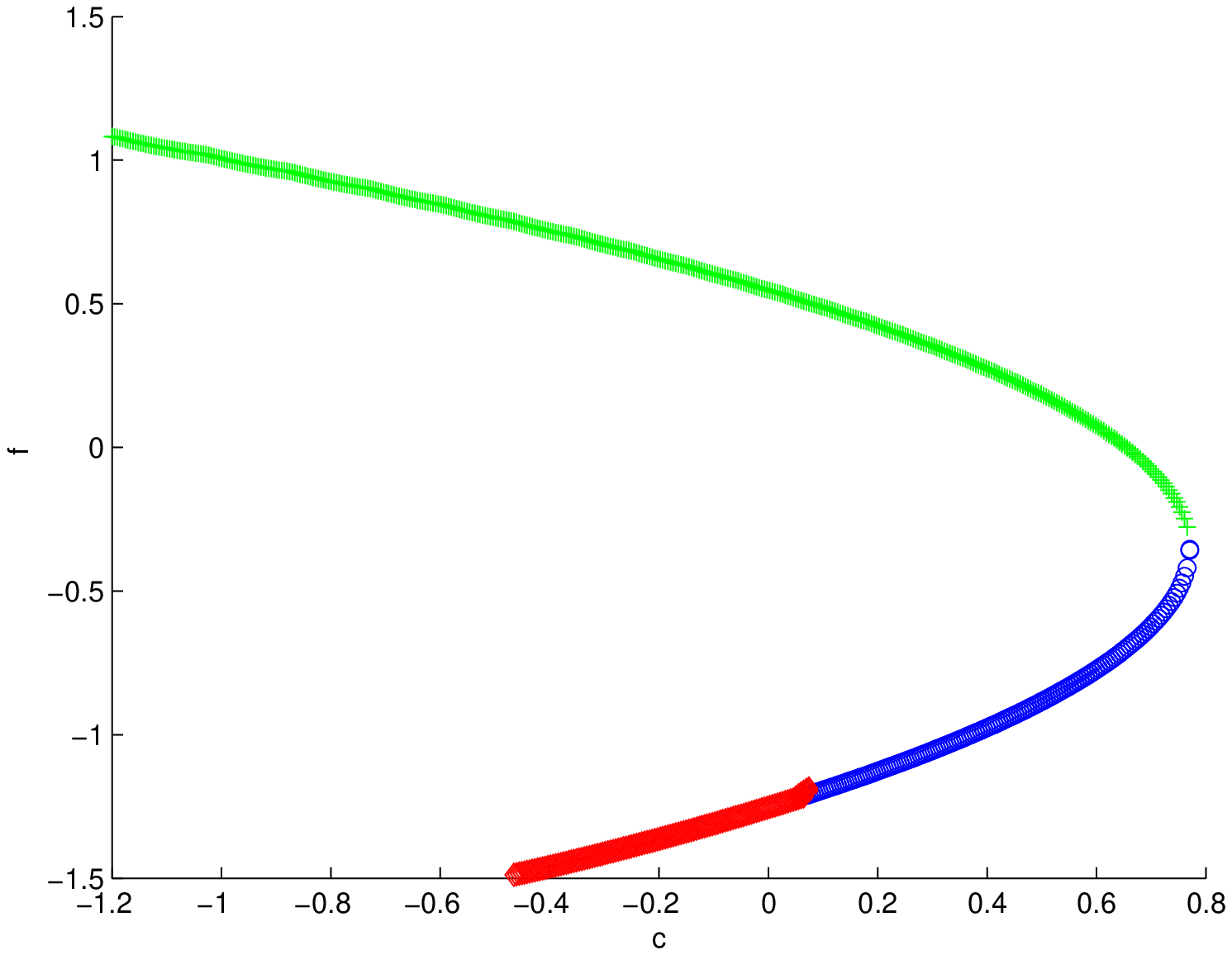}
\includegraphics[height=2in]{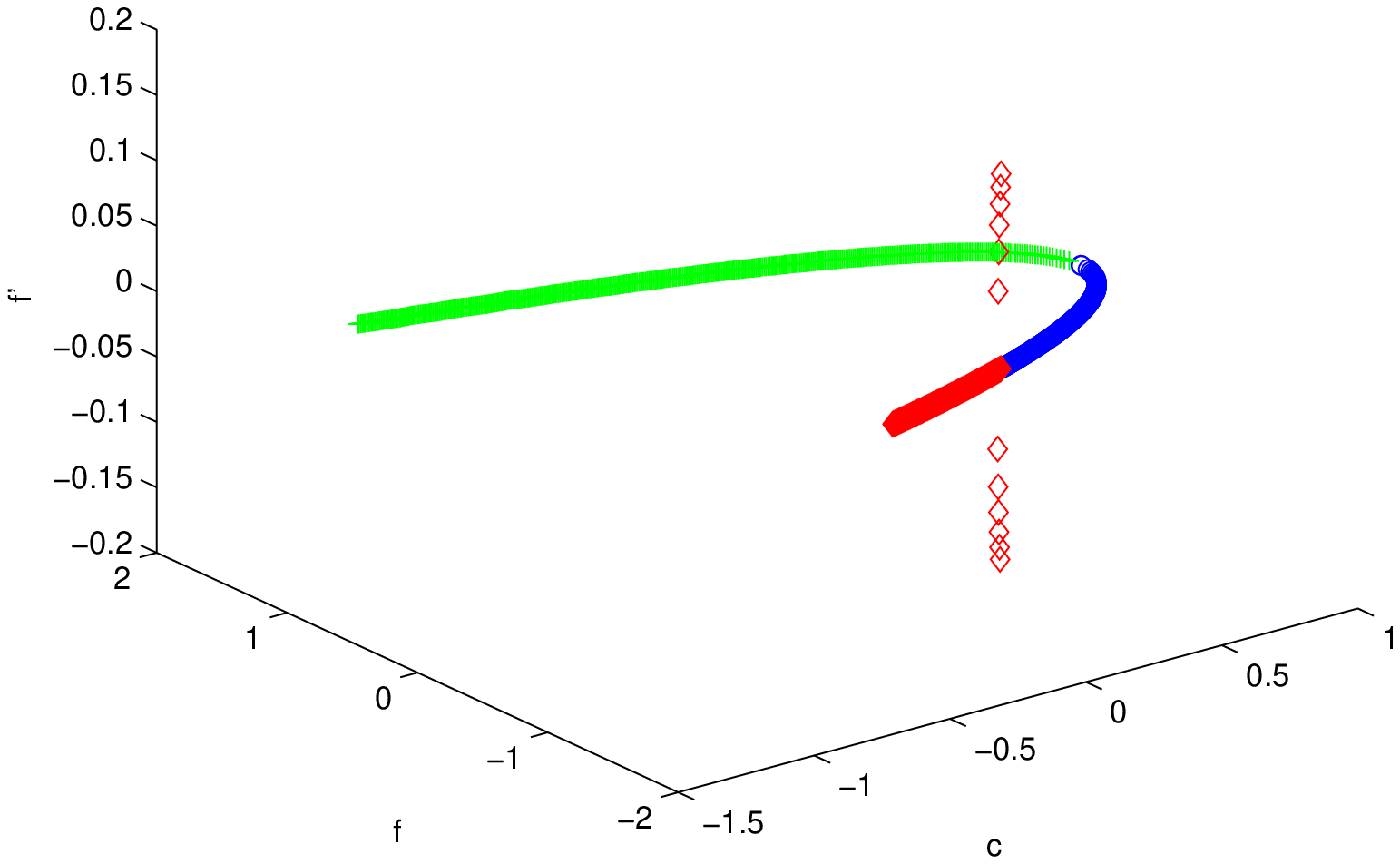}
\end{center}
\caption{Bifurcation diagram, coded by spectrum of
$\frac{d^2}{dx^2}-2f$: green = nonpositive spectrum, blue = one positive
eigenvalue, red = two positive eigenvalues}
\label{bif_diag}
\end{figure}

Once computed, the numerical solutions can then be tabulated
conveniently in a bifurcation diagram.  That is, consider the set in
$\mathbb{R}^3$ given by $(c,f(0),f'(0))$ for each solution $f$.
Evidently, by existence and uniqueness for ordinary differential
equations, each solution can be uniquely represented by such a point.
The results of such a computation are shown in Figure \ref{bif_diag}.
In this diagram, the solutions are color-coded by the number of
positive eigenvalues of $\frac{d^2}{dx^2}-2f$ as an operator
$C^2(\mathbb{R})\to C^0(\mathbb{R})$, which will be shown in Chapter
\ref{unstable_ch} to be the dimension of the unstable manifold of an
equilibrium solution $f$. (It should be noted that the green curve
continues for $c<-1.2$, but was stopped for display reasons.)

\begin{figure}
\begin{center}
\includegraphics[height=3in]{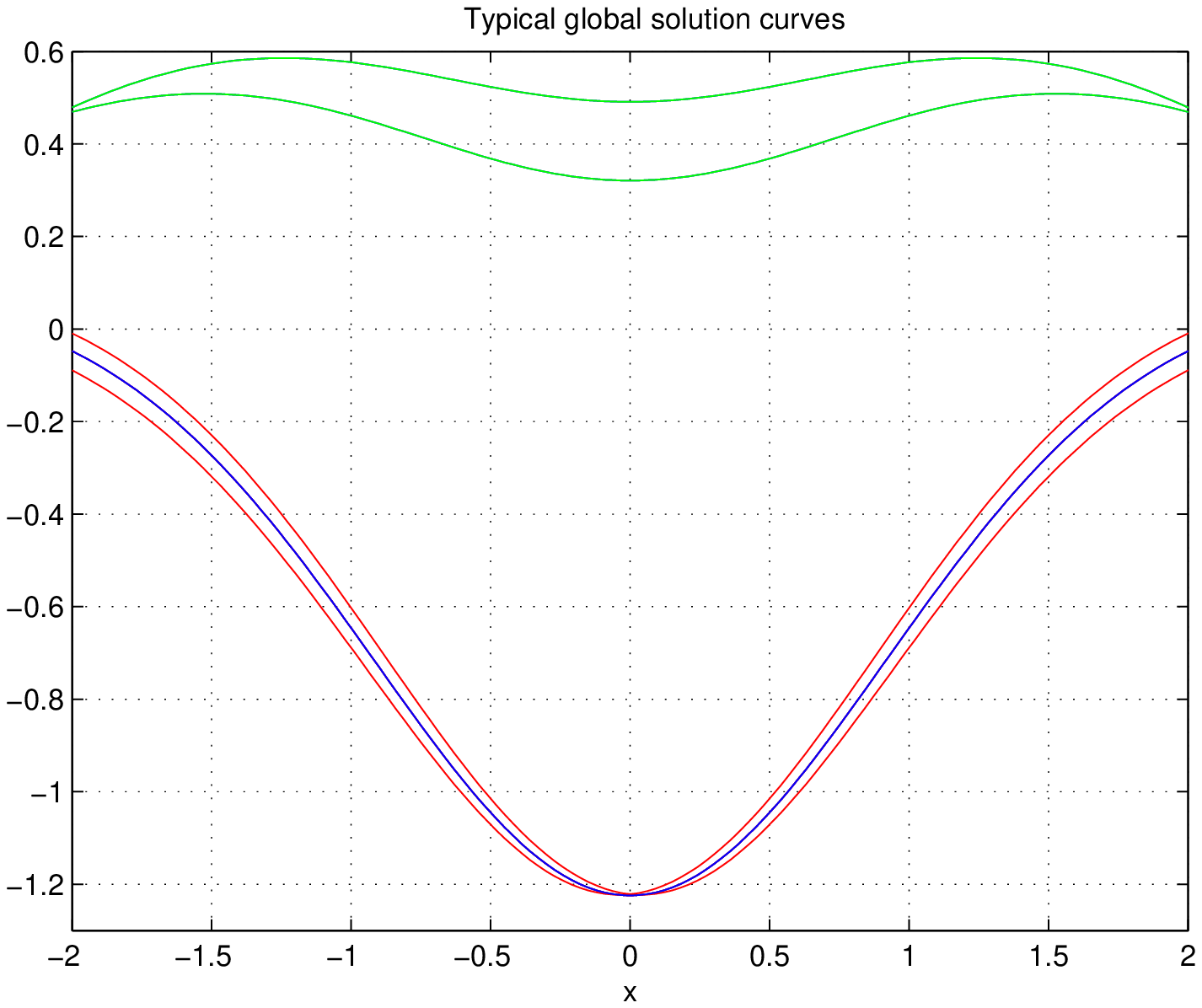}
\end{center}
\caption{Typical global solutions: green are from the positive branch,
  the blue one is taken from the lower branch with $f'(0)=0$, and the red
  ones are from the fork arms}
\end{figure}

Considering the bifurcation diagram, it appears to indicate that
\eqref{ode1} undergoes a saddle node bifurcation at approximately
$c=0.7706$, and a subcritical pitchfork bifurcation at $c=0.0501$.
The results agree with Theorem \ref{long_soln_exists}, in that
solutions do exist when $c<0$.  The saddle node bifurcation was
anticipated by the general shape of $\phi$.  For $c>0.7706$, global
solutions do not exist, which was qualitatively predicted
by Corollary \ref{phi_integral_bound}.

\begin{figure}
\begin{center}
\includegraphics[height=3in]{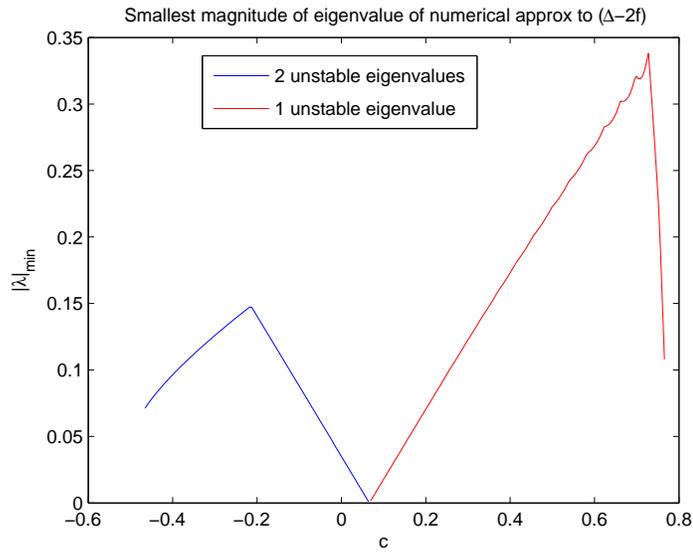}
\end{center}
\caption{Smallest-magnitude eigenvalue measured along the lower branch with $f'(0)=0$}
\label{small_eig}
\end{figure}

However, there are some stranger features of the bifurcation diagram.
Most prominently, the bifurcation diagram appears simply to {\it end}
near $c=-0.4652$, and at each branch of the pitchfork at $c=0.0740$.
It is important to verify that these are not numerical or
discretization errors.  If these {\it ends} are to be thought of as
valid bifurcations, very likely, $\frac{d^2}{dx^2} - 2f$ acquires a
zero eigenvalue there.  Plotting the smallest magnitude eigenvalue
gives some credence to this possibility.  (See Figure \ref{small_eig})

\begin{figure}
\begin{center}
\includegraphics[height=3in]{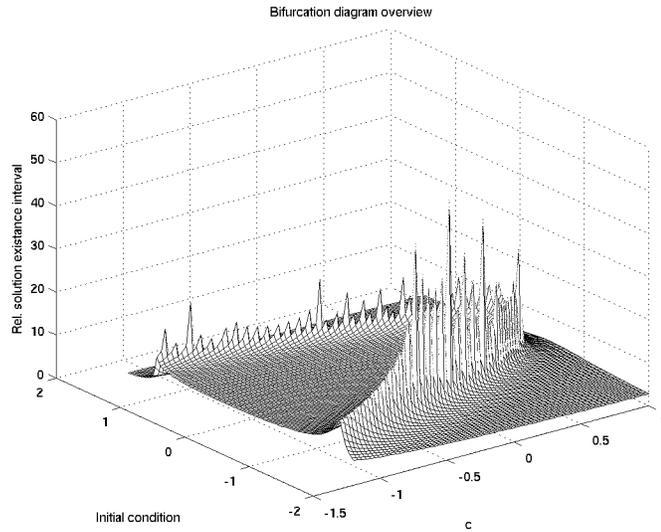}
\end{center}
\caption{Estimate of existence interval length}
\label{bif_exist}
\end{figure}

As another check, one can measure the size of the existence interval
for solutions to \eqref{long_ode1}, centered at $x=0$.  Looking in the
$(c,f(0))$-plane (taking $f'(0)=0$), one can find the first $x$ such
that the solution exceeds a particular value.  This is shown in Figure
\ref{bif_exist}, in which one sees the same general shape as in the
bifurcation diagram.  (The jagged nature of the graph along the actual
bifurcation diagram is due to aliasing.)  However, for $c<-0.4652$,
the lower branch clearly continues into solutions that exist for only
finite $x$.  So the end bifurcation indicates a failure of the
solutions to \eqref{long_ode1} to exist for all $x$.

From the point of view of \eqref{limited_pde_nonauto} (the parabolic
problem), the end bifurcations indicate that the equilibria are
degenerate in the sense of Morse.  It is easy to construct a
1-parameter family of 2-dimensional flows for which end bifurcations
occur.  The resulting equilibrium solutions in that case always
acquire a center manifold -- essentially a zero eigenvalue as noted
above.  In an infinite-dimensional flow, however, equilibria can be
degenerate without having a center manifold.  This is a manifestation
of the fact that infinite-dimensional spaces are not locally compact.
In Chapter \ref{instability_ch}, we show that in fact all equilibrium
solutions are asymptotically unstable.  However, in Chapter
\ref{unstable_ch}, we find that all equilibria have finite-dimensional
unstable manifolds whose dimension is determined by the dimension of
the positive eigenspace of $\frac{d^2}{dx^2} - 2f$.  In particular,
there are equilibria whose unstable manifolds are empty, yet they are
unstable.

\section{Conclusions}

In this chapter, an approach for counting and approximating global
solutions to a nonlinear, nonautonomous differential equation was
described that combines asymptotic and numerical information.  The
asymptotic information alone is enough to give necessary and
sufficient (but not sharp) conditions for solutions to exist, and
provides a fairly weak uniqueness condition.  More importantly, the
asymptotic approximation can be used to supply enough information to
pose a boundary value problem on a bounded interval containing a
smaller interval where asymptotic approximation is not valid.  This
boundary value problem is well-suited for numerical examination, and
the combined approach yields much more detailed results
than either method alone.

The techniques and results of this chapter should apply to more
general kinds of differential equations.  Indeed, there should be no
particular obstruction to extending any of the analysis to
equations of the form
\begin{equation*}
0=f''(x)-f^N(x)+\phi(x),
\end{equation*}
where $N>2$.  Regions $R_1$ and $R_2$ are then relatively easy to
construct, and similar results hold for them as are shown here.
Additionally, the asymptotic series for large $x$ can be obtained
since $0=f''-f^N$ has easily-found explicit solutions.  

There is considerably more difficulty in trying to understand
solutions to equations like
\begin{equation}
\label{more_general_ode}
0=f''(x)+G(f)+\phi(x),
\end{equation}
where $G$ is some polynomial with $G(0)=0$.  In this more general
setting, explicit solutions to $0=f''+G(f)$ are significantly harder
to find and work with.  There is also no reason to expect that Theorem
\ref{limzero_iff_bounded_lem} will hold, since $G$ may have several
zeros.  As a result, the asymptotic analysis becomes essentially
unavailable.  Therefore, the only remaining tool is Wa\.{z}ewski's
antifunnel theorem, for which one still needs a good description of
the asymptotic behavior of solutions.

\chapter{Existence of nontrivial eternal solutions}
\label{global_ch}
\section{Introduction}

The existence of eternal solutions poses a potentially difficult
problem, because the backward-time Cauchy problem is well known to be
ill-posed.  Obviously, equilibrium solutions are trivial examples of
such eternal solutions, and in Chapter \ref{nonauto_ch} we showed that
they exist.  It is not at all clear that there are other eternal
solutions, and indeed there may not be.  In this chapter we use a pair
of nonintersecting equilibrium solutions to construct a
heteroclinic orbit which connects them.  Therefore, the set of
heteroclinic orbits is generally nonempty.

As has been done in previous chapters, we will work with the more
limited equation \eqref{limited_pde}
\begin{equation}
\label{pde1}
\frac{\partial u(t,x)}{\partial t} = \frac{\partial^2 u(t,x)}{\partial
  x^2} - u^2(t,x) + \phi(x),
\end{equation}
where $\phi$ is a certain smooth function which decays to zero.  In
particular, when there are at least two equilibria whose difference is
never zero, there exist nonequilibrium eternal solutions to
\eqref{pde1}.  By an eternal solution, we mean a classical solution that is defined for
all $t\in\mathbb{R}$ and $x\in\mathbb{R}$.  We follow the general
technique for constructing ``ancient solutions,'' which was used in a
different context by Perelman.

\section{Equilibrium solutions}

We choose $\phi(x)=(x^2-0.4)e^{-x^2/2}$.  It has been shown in Chapter
\ref{nonauto_ch} (see Figure \ref{bif_diag}), that in this situation, there exists a pair of equilibrium solutions $f_+,f_-$ with the following properties:

\begin{enumerate}
\item $f_+$ and $f_-$ are smooth and bounded,
\item $f_+$ and $f_-$ have bounded first and second derivatives,
\item $f_+$ and $f_-$ are asymptotic to $6/x^2$ for large $x$, and so
  both belong to $L^1(\mathbb{R})$,
\item $f_+(x)>f_-(x)$ for all $x$,
\item there is no equilibrium solution $f_2$ with
  $f_+(x)>f_2(x)>f_-(x)$ for all $x$,
\end{enumerate}
and additionally, there exists a one-parameter family $g_c$ of
solutions to 
\begin{equation}
\label{gc_eqn}
0=g_c''(x) - g_c^2(x) + \phi_c(x)
\end{equation}
with 
\begin{enumerate}
\item $c\in[0,1)$,
\item $g_0=f_-$ and $\phi_0=\phi$,
\item $\phi_a(x) < \phi_b(x)$ and $g_a(x) > g_b(x)$ for all $x$ if $a>b$.
\end{enumerate}

The latter set of properties can occur as a consequence of the
specific structure of $f_-$.  For instance, consider the following result.

\begin{prop}
  Suppose $f_- \in C^{2,\alpha}(\mathbb{R})$ satisfies the above
  conditions and additionally, there is a compact $K \subset
  \mathbb{R}$ with nonempty interior such that $f_-$ is negative on
  the interior of $K$ and is nonnegative on the complement of $K$.
  Then such a family $g_c$ above exists.
\begin{proof}(Sketch)
Work in $T_{f_-}C^{2,\alpha}(\mathbb{R})$, the tangent space at
$f_-$.  Then \eqref{gc_eqn} becomes its linearization (for $h_c$, say), namely
\begin{equation}
\label{gc_lin_eqn}
0=h_c''(x) - 2 f_-(x) h_c(x) + (\phi_c - \phi).
\end{equation}
Consider the slighly different problem, 
\begin{equation}
\label{gc_lin_eg_eqn}
0=y''(x)-2f_-(x) y(x) + v(x) y(x),
\end{equation}
where $v$ is a smooth function to be determined.  If we can find a $v
\le 0$ such that $y>0$ and $y \to 0$ as $|x| \to \infty$, then we are
done, because we simply let $vy = \phi_c - \phi$ in
\eqref{gc_lin_eqn}.  In that case, $h_c=y$ has the required
properties.  
We sketch why such a $v$ exists:
\begin{itemize}
\item If $v \equiv 0$, then $y \equiv 0$ is a solution, giving
  $g_c=f_-$ as a base case.
\item If $v(x) = - 2 \|u\|_\infty \beta(x)$ for $\beta$ is a smooth
  bump function with compact support and $\beta | K = 1$, then the
  Sturm-Liouville comparison theorem implies that $y$ has no sign
  changes.  We can take $y$ strictly positive.  However, in this case,
  the Sturm-Liouville theorem imples that ther are no critical points
  of $y$ either, so $y$ may not tend to zero as $|x|\to \infty$.
\item Hence there should exist an $s$ with $0<s<2\|u\|_\infty$ such
  that if $v(x) = -s \beta(x)$, then $y$ has no sign changes, one
  critical point, and tends to zero as $|x|\to\infty$.  This choice of
  $v$ is what is required.  (The precise details of this argument fall
  under standard Sturm-Liouville theory, which are omitted here.)
\end{itemize}
\end{proof}
\end{prop}

In what follows, we shall not be concerned with the exact form of
$\phi$, but rather we shall assume that the above properties of the
equilibria hold.  Many other choices of $\phi$ will allow a similar
construction.

\begin{lem}
\label{lucid_lem}
The set
\begin{equation}
\label{funnel_eqn}
W=\{v \in C^2(\mathbb{R}) | f_-(x) < v(x) < f_+(x)\text{ for all }x\}
\end{equation}
is a forward invariant set for \eqref{pde1}.  That is, if $u$ is a solution to
\eqref{pde1} and $u(t_0) \in W$, then $u(t) \in W$ for all $t>t_0$.
\begin{proof}
We show that the flow of \eqref{pde1} is inward whenever a timeslice is tangent to either $f_-$ or $f_+$.  To this end, define the
set $B$
\begin{eqnarray*}
B&=&\{v \in C^2(\mathbb{R}) | f_-(x) \le v(x) \le f_+(x)\text{  for all }x,\text{ and there exists an }x_0\\&&\text{ such that }
  v(x_0)=f_+(x)\text{ or } v(x_0)=f_-(x)\}.
\end{eqnarray*}
Without loss of generality, consider a $v\in B$ with a single
point of tangency, $v(x_0)=f_-(x_0)$.  At such a point $x_0$, the
smoothness of $v$ and $f_-$ implies that $\Delta v(x_0) \ge \Delta
f_-(x_0)$ using the maximum principle.  Then, if $u$ is a solution to
\eqref{pde1} with $u(0,x)=v(x)$, we have that
\begin{eqnarray*}
\frac{\partial u(0,x_0)}{\partial t} &=& \Delta v(x_0) - v^2(x_0) +
\phi(x_0)\\ 
&\ge&\Delta f_-(x_0) - f_-^2(x_0) + \phi(x_0)=0,\\
\end{eqnarray*}
hence the flow is inward.  One can repeat the above argument for each
point of tangency, and for tangency with $f_+$ as well.
\end{proof}
\end{lem}

\begin{lem}
\label{nested_funnel_lem}
Solutions to the Cauchy problem 
\begin{equation}
\label{pde_cauchy}
\begin{cases}
\frac{\partial u(t,x)}{\partial t} = \frac{\partial^2 u(t,x)}{\partial
  x^2} - u^2(t,x) + \phi(x),\\
u(0,x) = U(x) \in W_c
\end{cases}
\end{equation}
where
\begin{equation*}
W_c=\{v \in C^2(\mathbb{R}) | g_c(x) < v(x) < f_+(x)\text{ for all }x\}
\end{equation*}
for $c \in [0,1)$ have the property that they lie in $L^1 \cap
L^\infty(\mathbb{R})$ for all $t>0$.  We shall assume that $U$ has
bounded first and second derivatives.

Additionally, when $c \in (0,1)$, solutions to \eqref{pde_cauchy}
cannot have $f_-$ as a limit as $t\to\infty$.
\begin{proof}
  The fact that solutions lie in $L^1 \cap L^\infty(\mathbb{R})$ is
  immediate from Lemma \ref{lucid_lem} and the asymptotic behavior of
  $f_+,f_-$ (Section \ref{series_sec}).  Observe that for each $c\in
  [0,1)$, $W_c$ is forward invariant, and that $W_a \subset W_b$ if
  $a>b$.  Since $f_-$ is not in $W_c$ for $c$ strictly larger than 0,
  the proof is completed.
\end{proof}
\end{lem}

The following is an outline for the rest of the chapter.  All solutions to \eqref{pde_cauchy} have bounded
first and second spatial derivatives.  This implies that all of their
first partial derivatives are bounded (the time derivative is
controlled by \eqref{pde1}).  Using the fact that \eqref{pde1} is
autonomous in time, time translations of solutions are also solutions.
We therefore construct a sequence of solutions $\{u_k\}$ to Cauchy
problems started at $t=0,T_{1}, T_{2}, ...$ which tend to $f_+$ as
$t\to +\infty$, but their initial conditions tend to $f_-$ as $k\to
\infty$.  By Ascoli's theorem, this sequence converges uniformly on
compact subsets to a continuous eternal solution.

\section{Integral equation formulation}
In order to estimate the derivatives of a solution to
\eqref{pde_cauchy}, it is more convenient to work with an integral
equation formulation of \eqref{pde_cauchy}.  This is obtained in the
usual way.

\begin{eqnarray*}
\frac{\partial u}{\partial t} &=& \Delta u - u^2 + \phi\\
\left(\frac{\partial}{\partial t} - \Delta \right) u &=& - u^2 + \phi\\
u &=& \left(\frac{\partial}{\partial t} - \Delta \right)^{-1} ( \phi - u^2)\\
\end{eqnarray*}
\begin{equation}
\label{int_eqn}
u(t,x)=\int_{-\infty}^\infty H(t,x-y) U(y) dy + \int_0^t
\int_{-\infty}^\infty H(t-s,x-y)\left(\phi(y) - u^2(s,y)\right) dy\, ds,
\end{equation}
where $H(t,x)=\frac{1}{\sqrt{4 \pi t}}e^{- \frac{x^2}{4t}}$ is the
usual heat kernel.  

\begin{calc}
\label{first_deriv_short}
We begin by estimating the first derivative of $u$ for a short time.
Let $T>0$ be given, and consider $0\le t\le T$.  The key fact is that
$\int H(t,x) dx = 1$ for all $t$.  Using \eqref{int_eqn}

\begin{eqnarray*}
\left \| \frac{\partial u}{ \partial x} \right \|_\infty &\le& \left \|
\frac{\partial U}{\partial x} \right \|_\infty + \left|\int_0^t \int_{-\infty}^\infty
\frac{\partial}{\partial x}H(t-s,x-y)\left(\phi(y) - u^2(s,y)\right)
dy\, ds\right|\\
&\le& \left \|\frac{\partial U}{\partial x} \right \|_\infty +
\int_0^t \int_{-\infty}^\infty
\left|\frac{\partial}{\partial y}(H(t-s,x-y))\left(\phi(y) -
u^2(s,y)\right) \right| dy\, ds\\
&\le& \left \|\frac{\partial U}{\partial x} \right \|_\infty +
\int_0^t \int_{-\infty}^\infty
\left|H(t-s,x-y)\left(\frac{\partial \phi}{\partial y} - 2u
\frac{\partial u}{\partial y}\right) \right| dy\, ds\\
&\le& \left \|\frac{\partial U}{\partial x} \right \|_\infty +
T \left\|\frac{\partial \phi}{\partial x}\right\|_\infty +
2 \| u\|_\infty \int_0^t \left\|
\frac{\partial u}{\partial x} \right\|_\infty ds.\\
\end{eqnarray*}
This integral equation fence is easily solved to give
\begin{eqnarray*}
\left \| \frac{\partial u}{ \partial x} \right \|_\infty &\le& \left (
\left\| \frac{\partial U}{\partial x}\right \|_\infty + T\left\|
\frac{\phi}{\partial x}\right \|_\infty \right) e^{2 t
  \max\{\|f_+\|_\infty,\|f_-\|_\infty\}}\\
&\le&K_1 e^{K_2 T}.
\end{eqnarray*}
\end{calc}

\begin{calc}
\label{second_deriv_short}
With the same choice of $T$ as above, we find a bound for the second
derivative in the same way:

\begin{eqnarray*}
\left \| \frac{\partial^2 u}{\partial x^2} \right \|_\infty &\le&
\left \|\frac{\partial^2 U}{\partial x^2} \right \|_\infty + T \left \|
\frac{\partial^2 \phi}{\partial x^2}\right \|_\infty + \int_0^t \left
\| \frac{\partial}{\partial y} \left( 2 u \frac{\partial u}{\partial
  y}\right)\right\|_\infty ds\\
 &\le&
\left \|\frac{\partial^2 U}{\partial x^2} \right \|_\infty + T \left \|
\frac{\partial^2 \phi}{\partial x^2}\right \|_\infty + \int_0^t 2
\left \| \frac{\partial u}{\partial x} \right \|^2_\infty + 2 \|
u \|_\infty \left \| \frac{\partial^2 u}{\partial x^2}\right \|_\infty
ds\\
&\le& K_3 e^{K_2 T}
\end{eqnarray*}
for some $K_3$ which depends on $U$, $\phi$, and $T$.
\end{calc}

\begin{calc}
\label{second_deriv_long}
Now, we extend Calculation \ref{second_deriv_short} to handle $t>T$,

\begin{eqnarray*}
\left \| \frac{\partial^2 u}{\partial x^2} \right \|_\infty &\le&
\left \| \frac{\partial^2 U}{\partial x^2} \right \|_\infty +
\left|\frac{\partial^2}{\partial x^2}\int_0^T \int_{-\infty}^\infty
H(t-s,x-y) (\phi(y) - u^2(s,y)) dy\, ds\right | +\\&& 
\left | \frac{\partial^2}{\partial x^2}\int_T^t \int_{-\infty}^\infty
H(t-s,x-y) (\phi(y) - u^2(s,y)) dy\, ds \right |
\\
&\le&
K_3 e^{K_2 T} + 
\int_T^t \left\|\frac{\partial^2}{\partial x^2} H(t-s,x)\right\|_\infty
(\|\phi\|_1 + \|u^2\|_1) ds \\
&\le&
K_3 e^{K_2 T} + K_4 \int_T^t \frac{1}{s\sqrt{s}} ds +K'_4 \int_T^t \frac{1}{s^2\sqrt{s}} ds\\
&\le&
K_3 e^{K_2 T} + K_5 \left(\frac{1}{\sqrt{T}} -
\frac{1}{\sqrt{t}}\right) +K'_5 \left(\frac{1}{T\sqrt{T}} - \frac{1}{t\sqrt{t}}\right)\\
&\le&
K_3 e^{K_2 T} + K_6,\\
\end{eqnarray*}
hence there is a uniform upper bound on $\left \| \frac{\partial^2
  u}{\partial x^2} \right \|_\infty$ which depends only on the initial
  conditions, $\phi$, and $T$.
\end{calc}

\begin{lem}
\label{bdd_deriv}
Let $f\in C^2(\mathbb{R})$ be a bounded function with a bounded second
derivative.  Then the first derivative of $f$ is also bounded, and the
bound depends only on $\|f\|_\infty$ and $\|f''\|_\infty$.
\begin{proof}
The proof is elementary.  The key fact is that
at its maxima and minima, $f$ has a horizontal tangent.  From a
horizontal tangent, the quickest $f'$ can grow is at a rate of
$\|f''\|_\infty$.  However, since $f$ is bounded, there is a maximum
amount that this growth of $f'$ can accrue.  Indeed, a sharp estimate
is
\begin{eqnarray*}
\|f'\|_\infty \le \sqrt{2 \|f\|_\infty \|f''\|_\infty}.
\end{eqnarray*}
\end{proof}
\end{lem}

Using the fact that $u$ is bounded, Lemma \ref{bdd_deriv} implies that
the first spatial derivative of $u$ is bounded.  By \eqref{pde1}, it
is clear that the first time derivative of $u$ is also bounded.

\begin{lem}
\label{bdd_action}
As an immediate consequence of Lemmas \ref{nested_funnel_lem} and
\ref{bdd_deriv}, the action integral
\begin{equation*}
A(u(t))=\int_{-\infty}^\infty \frac{1}{2}\left| \frac{\partial
  u}{\partial x} \right|^2 + \frac{1}{3}u^3(t,x) - u(t,x)\phi(x) dx
\end{equation*}
is bounded.  Therefore, the solutions to the Cauchy problem
\eqref{pde_cauchy} all tend to limits as $t\to\infty$ (Corollary
\ref{limits_to_equilibria}).  By Lemma \ref{nested_funnel_lem}, we
conclude that they all tend to the common limit of $f_+$ when $c>0$.
\begin{proof}
  The latter two terms are bounded due to the fact that $u$ lies in
  $L_1 \cap L^\infty(\mathbb{R})$ for all $t$.  The bound on the first
  term comes from combining the fact that $u$ and its first two
  spatial derivatives are bounded with the asymptotic decay of
  $f_\pm$, and is otherwise straightforward (use L'H\^opital's rule).
\end{proof}
\end{lem}

\section{Construction of an eternal solution}

Let 
\begin{equation*}
U_k(x) = (1-2^{-k-1})g_{1/k}(x) + 2^{-k-1}f_+(x),\text{ for }k \ge 0
\end{equation*}
noting that $U_k \to f_-$ as $k\to \infty$.  Since $U_k$ is a convex
combination of $f_+$ and $g_{1/k}$, it follows that $U_k \in W_{1/k}$ for all
$k$.  Also, since $f_+$ and $f_-$ have bounded first and second
derivatives, the $\{U_k\}$ have a common bound for their first and second
derivatives.

Now consider solutions to the following set of Cauchy problems
\begin{equation}
\label{cauchy_set}
\begin{cases}
\frac{\partial u_k(t,x)}{\partial t} = \frac{\partial^2 u_k(t,x)}{\partial
  x^2} - u_k^2(t,x) + \phi(x),\\
u_k(T_k,x) = U_k(x).
\end{cases}
\end{equation}  

We choose $T_k$ so that for all $k>0$, $u_k(0,0)=u_0(0,0)$.  We can do
this using the continuity of the solution and Lemma \ref{bdd_action}.  As
$k\to \infty$, solutions are started nearer and nearer to the
equilibrium $f_-$, so we are forced to choose $T_k \to -\infty$
as $k\to \infty$.

It's clear that each solution $u_k$ is defined for only $t>T_k$.  
However, for each compact set $S \subset \mathbb{R}^2$, there are
infinitely many elements of $\{u_k\}$ which are defined on it.  The
results of the previous section imply that $\{u_k\}$ is a bounded,
equicontinous family.  As a result, Ascoli's theorem implies that
$\{u_k\}$ converges uniformly on compact subsets to a continous $u$,
which is an eternal solution to \eqref{pde1}.

Our constructed eternal solution will have the value $u(0,0)=u_0(0,0)$,
which is strictly between $f_+$ and $f_-$.  As a result, the eternal
solution we have constructed is not an equilibrium solution.  By Lemma
\ref{bdd_action}, it is a finite energy solution, so it must be a
heteroclinic orbit connecting $f_-$ to $f_+$.

\chapter{Instability of equilibria}
\label{instability_ch}
\section{Introduction}

(This chapter is available on the arXiv as \cite{RobinsonInstability}.)

If we try to apply standard Morse theory to our semilinear parabolic
equation, we encounter a serious difficulty.  In particular, the
stability of the linearization about an equilibrium of our
system is not sufficient to ensure that the equilibrium is stable,
even though there may not be a zero eigenvalue.  The instability of
equilibria for \eqref{pde} is not a new fact, having been studied
carefully in the 1980s.  This chapter is included for completeness,
providing an explicit construction of a sequence of solutions starting
near the equilibrium which all blow up.  Indeed, there is a
complementary result of {\it stability} in certain weighted norms,
described in \cite{Zhang_2000} and \cite{Souplet_2002}.
  
Note that the right side of \eqref{pde} is an operator which has a
spectrum which includes zero, so stability is possible (as in the
unforced heat equation), though not guaranteed.  This is in stark
contrast to the situation in finite-dimensional settings, where
asymptotic stability of the linearized system implies stability of the
equilibrium (in particular, zero is not an eigenvalue of the
linearized operator).  (See \cite{BoyceDiPrima}, for instance.)
Essentially, this difficulty suggests that each critical point is
degenerate.  (This line of reasoning is completed in Chapter
\ref{unstable_ch}.)  This implies that Morse theory (even when
strengthened to its natural infinite-dimensional form
\cite{Palais_1963}) cannot be used to study the dynamics of our
system.

Again, we study the general situation by working with the simpler
Cauchy problem \eqref{limited_pde}.  (The computations we exhibit in
this chapter will carry over {\it mutatis mutandis} to \eqref{pde}.)
However, it is useful to center on an equilibrium solution.  That is,
we apply the change of variables $u(t,x) \mapsto u(t,x) - f(x)$ to
obtain
\begin{equation}
\label{pde_instability}
\begin{cases}
\frac{\partial}{\partial t} u(t,x) = \frac{\partial^2}{\partial x^2} u(t,x) - 2 f(x) u(t,x) -
u^2(t,x)\\
u(0,x)=h(x) \in C^\infty(\mathbb{R})\\
t>0,x\in \mathbb{R},\\
\end{cases}
\end{equation}
where $f \in C_0^\infty(\mathbb{R})$ is a positive function with two
bounded derivatives.  (By $C_0^\infty$, we mean the space of smooth
functions which decay to zero.)  We interpret $f$ as being an
equilibrium of the original problem \eqref{limited_pde}.

The assumption that $f$ be {\it positive} requires some motivation.
With this assumption, the zero function is an asymptotically stable equilibrium
for the linearized problem,
\begin{equation}
\label{linearized_pde}
\frac{\partial}{\partial t} u(t,x) = \frac{\partial^2}{\partial x^2} u(t,x) - 2 f(x) u(t,x).
\end{equation}
by a standard comparison principle argument.  This corresponds neatly
to the case of zero positive eigenvalues of
$\frac{\partial^2}{\partial x^2}-2f$ which was found numerically in
Figure \ref{bif_diag}.  Intuition would suggest that this implies $f$
is a stable equilibrium of \eqref{limited_pde}.  However, using a
technique pioneered by Fujita in \cite{Fujita}, we will show that this
equilibrium is not stable in the nonlinear problem, even if the
initial condition has small $p$-norm for every $1 \le p \le \infty$.

Fujita showed that if $f \equiv 0$, then the zero function is an
unstable equilibrium of \eqref{pde_instability}.  The cause of the
instability in \eqref{pde_instability} is the decay of $f$, for if
$f=\text{const}>0$, then the comparison principle shows that the zero
function is stable.  We extend Fujita's result, so that roughly
speaking, since $f \to 0$ away from the origin, the system is less
stable to perturbations away from the origin.  Another indication that
there may be instability lurking (though not conclusive proof) is that
the decay of $f$ means that the spectrum of the linearized operator on
the right side of \eqref{linearized_pde} includes
zero.

\section{Motivation}
The problem \eqref{pde_instability} describes a reaction-diffusion equation
\cite{FiedlerScheel}, or a diffusive logistic population model with a
spatially-varying carrying capacity.  The choice of $f$ positive means
that the equilibrium $u \equiv 0$ describes a population saturated at
its carrying capacity.  Without the diffusion term, this situation is
well known to be stable.  The decay condition on $f$ means that the
carrying capacity diminishes away from the origin.

The spatial inhomogeneity of $f$ makes the analysis of \eqref{pde_instability}
much more complicated than that of typical reaction-diffusion
equations.  The existence of additional equilibria for \eqref{pde_instability}
is a fairly difficult problem, which depends delicately on $f$.  (See
\cite{Brezis_1984} for a proof of existence of equilibria in a
related setting.)

\section{Instability of the equilibrium}
\label{instability_sec}

Given an $\epsilon>0$, we will construct an initial condition $h \in
C^\infty(\mathbb{R})$ for the problem \eqref{pde_instability}, with $\|h\|_p <
\epsilon$ for each $1\le p \le \infty$, such that $\|u(t)\|_\infty \to
\infty$ as $t \to T < \infty$.  In particular, this implies that
$u\equiv 0$ is not a stable equilibrium of \eqref{pde_instability}, at least
insofar as classical solutions are concerned.  We employ a technique
of Fujita, which provides sufficient conditions for equations like
\eqref{pde_instability} to blow up \cite{Fujita}. (Additionally, \cite{Evans}
contains a more elementary discussion of the technique with a similar
construction.)  Our choice for $h$ can be thought of as a sequence of
progressively shifted gaussians, and we will demonstrate that though
each has smaller $p$-norm than the previous, the solution started at
$h$ still blows up.

\subsection{The technique of Fujita}

The technique of Fujita examines the blow-up behavior of nonlinear
parabolic equations by treating them as ordinary differential equations
on a Hilbert space.  Suppose $u(t)$ solves 
\begin{equation}
\label{fujita_eqn}
\frac{\partial u(t)}{\partial t} = L u(t) + N(u(t),t),
\end{equation}
where $L$ is a linear operator not involving $t$, and $N$ may be
nonlinear and may depend on $t$.  Suppose that $v(t)$ solves 
\begin{equation}
\label{fujita_adj_eqn}
\frac{\partial v(t)}{\partial t} = - L^* v(t),
\end{equation}
where $L^*$ is the adjoint of $L$.  Let $J(t) = \left<v(t),u(t)\right>$.  We
observe that if $|J(t)| \to \infty$ then either $\|v(t)\|$ or
$\|u(t)\|$ also does.  So if $v(t)$ does not blow up, then we can show
that $\|u(t)\|$ blows up, and perhaps more is true.  If we
differentiate $J(t)$, we obtain the identity
\begin{eqnarray*}
\frac{d}{dt}J(t) &=& \frac{d}{dt} \left< v(t), u(t) \right>\\
&=& \left<\frac{dv}{dt},u(t)\right>+\left<v(t),\frac{du}{dt}\right>\\
&=&\left<-L^*v(t),u(t)\right>+\left<v(t),Lu(t)+N(u(t),t)\right>\\
&=&\left<v(t),N(u(t),t)\right>, 
\end{eqnarray*}
where there is typically a technical justification required for the
second equality.  It is often possible to find a bound for
$\left<v(t),N(u(t),t)\right>$ in terms of $J(t)$.  So then the method
provides a fence (in the sense of \cite{HubbardWest}) for $J(t)$,
which we can solve to give a bound on $|J(t)|$.  As a result, the
blow-up behavior of $u(t)$ is controlled by the solution of an {\it
ordinary} differential equation (for $J(t)$) and a {\it linear}
parabolic equation (for $v(t)$), both of which are much easier to
examine than the original nonlinear parabolic equation.

\subsection{Instability in $L^p$ for $1 \le p \le \infty$}

We begin our application of the method of Fujita by working with
$L=\frac{\partial^2}{\partial x^2}-2f$ and $N(u)=- u^2$ in
\eqref{fujita_eqn}.  Since \eqref{fujita_adj_eqn} is then not
well-posed for all $t$, we must be a little more careful than the
method initially suggests.  For this reason, we consider a family of
solutions $v_\epsilon$ to \eqref{fujita_adj_eqn} that have slightly
extended domains of definition.  It will also be important, for technical
reasons, to enforce the assumption that the first and second
derivatives of $f$ are bounded.  

\begin{df}
\label{v_eps_2_df}
Suppose $w=w(t,x)$ solves
\begin{equation}
\label{w_eqn}
\begin{cases}
\frac{\partial w}{\partial t} = \frac{\partial^2 w}{\partial x^2} - 2
f(x) w(t,x)\\
w(0,x) = w_0(x) \ge 0.
\end{cases}
\end{equation}
Define $v_\epsilon(s,x) = w(t-s+\epsilon,x)$ for fixed $t>0$ and
$s<t+\epsilon$.  Notice that by the comparison principle,
$v_\epsilon(s,x) \ge 0$.
\end{df}

\begin{lem}
\label{w_endcond_lem}
Suppose that $w$ solves \eqref{w_eqn}.  Then $w, \frac{\partial
  w}{\partial x} \in C_0(\mathbb{R})$.
\begin{proof}
The standard existence and regularity theorems for linear parabolic
equations (see \cite{ZeidlerIIA}, for example) give that
$w,\frac{\partial w}{\partial x},\frac{\partial^2 w}{\partial x^2} \in
L^2(\mathbb{R})$ and that $w \in C^2(\mathbb{R})$.  The comparison
principle, applied to $\frac{\partial}{\partial t}\frac{\partial
w}{\partial x}$ and $\frac{\partial}{\partial t}\frac{\partial^2
w}{\partial x^2}$ gives that the first and second derivatives of $w$
are bounded for each fixed $t$.  (This uses our assumption that $f$
has two bounded derivatives.)

The lemma follows from a more general result: if $g\in C^1 \cap L^p
(\mathbb{R})$ for $1 \le p < \infty$ and $g' \in
L^\infty(\mathbb{R})$, then $g \in C_0(\mathbb{R})$.  To show this, we
suppose the contrary, that $\lim_{x \to \infty} g(x) \ne 0$ (and
possibly doesn't exist).  By definition, this implies that there is
an $\epsilon >0$ such that for all $x>0$, there is a $y$ satisfying
$y>x$ and $|g(y)|>\epsilon$.  Let $S=\{y|\,|g(y)|>\epsilon\}$, which
is a union of open intervals, is of finite measure, and has $\sup S =
\infty$.  Let $T = \{y | \, |g(y)|>\epsilon/2\}$.  Note that $T$
contains $S$, but since $g'$ is bounded, for each $x \in S$, there is
a neighborhood of $x$ contained in $T$ of measure at least
$\epsilon/\|g'\|_\infty$.  Hence, since $\sup T = \sup S = \infty$,
$T$ cannot be of finite measure, which contradicts the fact that $g
\in L^p(\mathbb{R})$ with $1\le p < \infty$.
\end{proof}
\end{lem}

\begin{lem}
\label{fence_2_lem}
Suppose $u:[0,T)\times \mathbb{R} \to \mathbb{R}$ is a classical
  solution to \eqref{pde_instability} with $u \le 0$ and $u(t)\in
  L^\infty(\mathbb{R})$ for each $t \in [0,T)$.  Then
\begin{equation}
\label{fence_2_eqn}
-\int w(t,x) h(x) dx \le \left( \int_0^t \frac{1}{\|w(s)\|_1}
 ds \right)^{-1},
\end{equation}
where $w$ is defined as in Definition \ref{v_eps_2_df}.

\begin{proof}
Define 
\begin{equation}
J_\epsilon(s)=\int{v_\epsilon(s,x)u(s,x) dx}.
\end{equation}

First of all, we observe that since $u\in L^\infty(\mathbb{R})$,
  $v_\epsilon(s,\cdot)u(s,\cdot)$ is in $L^1(\mathbb{R})$ for each
$s<t$.

Now suppose we have a sequence $\{m_n\}$ of compactly supported
smooth functions with the following properties: \cite{LeeSmooth}
\begin{itemize}
\item $m_n \in C^\infty(\mathbb{R})$,
\item $m_n(x) \ge 0$ for all $x$, \item $\text{supp}(m_n)$ is contained in the interval $(-n-1,n+1)$, and
\item $m_n(x)=1$ for $|x| \le n$.
\end{itemize}
Then it follows that 
\begin{equation*}
J_\epsilon(s)=\lim_{n\rightarrow \infty} \int{v_\epsilon(s,x)u(s,x)m_n(x) dx}.
\end{equation*}

Now 
\begin{eqnarray*}
\frac{d}{ds} J_\epsilon(s) &=& \frac{d}{ds} \lim_{n\rightarrow \infty}
\int{v_\epsilon(s,x)u(s,x)m_n(x) dx} \\
&=& \lim_{h\rightarrow 0} \lim_{n\rightarrow \infty} \frac{1}{h}\int
\left(v_\epsilon(s+h,x)u(s+h,x)-v_\epsilon(s,x)u(s,x)\right)m_n(x) dx.\\
\end{eqnarray*}
We'd like to exchange limits using uniform convergence.  To do this we
show that 
\begin{equation}
\label{big_lim}
\lim_{n\rightarrow \infty} \lim_{h\rightarrow 0} \frac{1}{h}\int
\left(v_\epsilon(s+h,x)u(s+h,x)-v_\epsilon(s,x)u(s,x)\right)m_n(x) dx
\end{equation}
exists and the inner limit is uniform.  We show both together by a
little computation, using uniform convergence and LDCT:

\begin{eqnarray*}
&&\lim_{n\rightarrow \infty} \lim_{h\rightarrow 0} \frac{1}{h}\int
\left(v_\epsilon(s+h,x)u(s+h,x)-v_\epsilon(s,x)u(s,x)\right)m_n(x) dx\\
&=&\lim_{n\rightarrow \infty} \int
\left(\frac{\partial}{\partial s}v_\epsilon(s,x)u(s,x)+v_\epsilon(s,x)\frac{\partial}{\partial s}u(s,x)\right)m_n(x) dx \\
&=&\lim_{n\rightarrow \infty} \int
\left(-\frac{\partial^2}{\partial x^2}
v_\epsilon(s,x)+2f(x)v_\epsilon(s,x)\right) u(s,x) m_n(x) +\\
&&v_\epsilon(s,x)\left(\frac{\partial^2}{\partial x^2} u(s,x) -
u^2(s,x)- 2 f(x) u(x) \right)m_n(x) dx \\
&=&\lim_{n\rightarrow \infty} \int
-v_\epsilon(s,x)u^2(s,x)m_n(x) dx.
\end{eqnarray*}
Minkowski's inequality has that
\begin{equation*}
\left|\int v_\epsilon u m_n dx \right|\le \int v_\epsilon |u| m_n dx \le \left(\int v_\epsilon m_n dx\right)^{1/2}
\left(\int v_\epsilon u^2 m_n dx \right)^{1/2},
\end{equation*}
since $v_\epsilon, m_n\ge 0$.  This gives that

\begin{eqnarray*}
&&\int
-v_\epsilon(s,x)u^2(s,x)m_n(x) dx\\
&\le& - \frac{(\int v_\epsilon u m_n dx)^2 }{\int v_\epsilon m_n dx}\\
&\le& - \frac{\left(\int v_\epsilon u dx \right)^2}{\int v_\epsilon m_1 dx}, \\
\end{eqnarray*}
hence the inner limit of \eqref{big_lim} is uniform.  On the other
hand, 
\begin{equation*}
|v_\epsilon(s,x)u^2(s,x)m_n(x)| \le v_\epsilon(s,x) \|u(s)\|_\infty^2
 \in L^1(\mathbb{R})
\end{equation*}
so the double limit of \eqref{big_lim} exists by dominated
convergence.  Thus we have the fence
\begin{equation}
\label{jeps_2_eqn}
\frac{d J_\epsilon(s)}{ds} \le -\frac{(J_\epsilon(s))^2}{\|v_\epsilon(s)\|_1}.
\end{equation}

We solve the fence \eqref{jeps_2_eqn} to obtain (note $J_\epsilon \le 0$)
\begin{eqnarray*}
\frac{1}{\|v_\epsilon(s)\|_1} &\le& - \frac{d J_\epsilon(s)}{ds}
\frac{1}{(J_\epsilon(s))^2}\\
\int_0^t \frac{1}{\|v_\epsilon(s)\|_1} ds &\le& \frac{1}{J_\epsilon(t)}
- \frac{1}{J_\epsilon(0)}\\
\int_0^t \frac{1}{\|v_\epsilon(s)\|_1} ds &\le& -\frac{1}{J_\epsilon(0)}.
\end{eqnarray*}
Taking the limit as $\epsilon \to 0$ of both sides of the inequality
yields
\begin{equation*}
-\int w(t,x) h(x) dx \le \left( \int_0^t \frac{1}{\|w(t-s)\|_1} ds
 \right)^{-1}
 =  \left( \int_0^t \frac{1}{\|w(s)\|_1} ds \right)^{-1},
\end{equation*}
as desired.
\end{proof}
\end{lem}

\begin{rem}
\label{l_inf_stability}
Since we are interested in proving the instability of the zero
function in \eqref{pde_instability}, consider $u(0,x) = h(x) = -
\epsilon$ for $\epsilon > 0$.  Then \eqref{fence_2_eqn} takes on the
simple form
\begin{equation}
\label{fence_2_improved_easy}
\epsilon \int_0^t \frac{\|w(t)\|_1}{\|w(s)\|_1} ds \le 1.
\end{equation}
So in particular, $\|u(t)\|_\infty$ blows up if there exists a $T>0$
such that $\epsilon \int_0^T \frac{\|w(T)\|_1}{\|w(s)\|_1} ds > 1.$

The stability of the zero function in \eqref{pde_instability} depends
on the stability of the zero function in \eqref{w_eqn} -- the
linearized problem.  If the zero function in the linearized problem is
very strongly attractive, say $\|w(t)\|_1 \sim e^{-t}$, then
\begin{equation*}
\int_0^t \frac{e^{-t}}{e^{-s}} ds = (1-e^{-t}) < 1,
\end{equation*}
and so a small choice of $\epsilon<1$ does not cause blow-up via a
violation of \eqref{fence_2_improved_easy}.  On the other hand,
blow-up occurs if it is less attractive, say $\|w(t)\|_1 \sim
t^{-\alpha}$ for $\alpha \ge 0$.  Because then
\begin{equation*}
\int_0^t \frac{s^\alpha}{t^\alpha} ds = \frac{t}{\alpha + 1},
\end{equation*}
whence blow-up occurs before $t=\frac{\alpha + 1}{\epsilon}$.

In the particular case of $f(x)=0$ for all $x$, we note that $w$ is
simply a solution to the heat equation, which has
$\|w(t)\|_1=\|w_0\|_1$ for all $t$ (by direct computation
using the fundamental solution, say), so blow up occurs.  Thus
we can recover a special case of the original blow-up result of Fujita
in \cite{Fujita}.
\end{rem}

\begin{thm}
\label{fence_2_improved_violation}
Suppose a sufficiently small $\epsilon>0$ is given.  Then for a
certain choice of initial condition $h(x)$ with
$\|h\|_p < \epsilon$ for all $1\le p \le \infty$, there exists a
$T>0$ for which $\lim_{t\to T^-}\|u(t)\|_\infty = \infty$.
\begin{proof}
First, it suffices to choose $\|u(0)\|_1<\epsilon$ and $\|u(0)\|_\infty<\epsilon$, since
\begin{equation*}
\|u\|_p = \left(\int |u|^p dx\right)^{1/p} \le \|u\|_\infty^{(p-1)/p}
\|u\|_1^{1/p} < \epsilon.
\end{equation*}  
We assume, contrary to what is to be proven, that $\|u(t)\|_\infty$
does not blow up for any finite $t$.  In other words, assume that $u:
[0,\infty)\times \mathbb{R} \to \mathbb{R}$ is a classical solution to
\eqref{pde_instability}, with $\|u(t)\|_\infty < \infty$ for all $t$.  We make
several definitions:
\begin{itemize}
\item Choose $0<\beta < \min\left\{\epsilon,\frac{\epsilon^4}{16 \pi^2}\right\}$.
\item Choose $\gamma>0$ small enough so that 
\begin{equation}
\label{vio_1}
\frac{\beta}{27 \gamma^2} = K,
\end{equation}
for some some arbitrary $K > 1$.
\item Since $0 \le f \in C_0^\infty(\mathbb{R})$, we can choose an
  $x_1$ such that
\begin{equation}
\label{vio_2}
f(x) \le \gamma \text{ when } x < x_1.
\end{equation}
\item Next, we choose $x_0<x_1$ so that
\begin{equation}
\label{vio_3}
\sqrt{t}\|f\|_\infty \left ( 1 - \text{erf} \left (
\frac{x_1-x_0}{2\sqrt{t}} \right ) \right ) < \gamma 
\end{equation}
for all $0 < t < \frac{1}{4 \gamma^2}$.  Notice that any choice less
than $x_0$ will also work.
\item Choose the initial condition for \eqref{pde_instability} to be
\begin{equation}
\label{initial_choice_eqn}
u(0,x)=h(x)=- \beta e^{\beta^{3/2} (x-x_0)^2}.
\end{equation}
This choice of initial condition has $\|u(0)\|_\infty =
\beta<\epsilon$, $\|u(0)\|_1 = 2 \pi^{1/2} \beta^{1/4}<\epsilon$, and $\left \|
\frac{\partial^2 u(0)}{\partial x^2}\right \|_\infty = \mu = 2
\beta^{5/2}.$  (The value of $\mu$ will be important shortly.)
\item Finally, let $w_0(y)=\delta(y-x_0)$ (the Dirac
  $\delta$-distribution), and suppose that $w$ solves \eqref{w_eqn}.
  In other words, choose $w$ to be the fundamental solution to
  \eqref{w_eqn} concentrated at $x_0$.  Note that the maximum
  principle ensures both that $w(t,x) \ge 0$ for all $t>0$ and $x \in
  \mathbb{R}$ and that $\|w(t)\|_1 \le \|w(0)\|_1 = 1$ for all $t>0$.
  This allows us to rewrite \eqref{fence_2_eqn} as
\begin{equation}
\label{fence_2_improved}
- t \int w(t,x) h(x) dx \le 1.
\end{equation}
\end{itemize}

Now we estimate the integral in \eqref{fence_2_improved}.  Notice that
\begin{eqnarray*}
\frac{d}{dt}\int w(t,x) \left(-h(x)\right) dx &=& 
\int \left(\frac{\partial^2 w}{\partial x^2} - 2 f(x) w(t,x)\right)
\left(-h(x)\right) dx \\
&=&\int \left(-\frac{\partial^2 u}{\partial x^2} + 2 f(x) h(x) \right) w(t,x)dx, \\
\end{eqnarray*}
where Lemma \ref{w_endcond_lem} eliminates the
boundary terms.  Now suppose $z$ solves the heat equation with the
same initial condition as $w$, namely
\begin{equation}
\label{z_eqn}
\begin{cases}
\frac{\partial z}{\partial t} = \frac{\partial^2 z}{\partial x^2}\\
z(0,x) = w_0(x) = \delta(x-x_0).
\end{cases}
\end{equation}
The comparison principle estabilishes that $z(t,x) \ge w(t,x)$ for all
$t>0$ and $x \in \mathbb{R}$, since $f,w \ge 0$.  As a result, we have
that
\begin{eqnarray*}
\frac{d}{dt}\int w(t,x) \left(-h(x)\right) dx &\ge&
\int \left(-\left |\frac{\partial^2 u}{\partial x^2}\right | + 2 f(x) h(x) \right)
z(t,x)dx \\
&\ge&
- \mu - 2 \beta \int f(x) z(t,x) dx,\\
\end{eqnarray*}
where $\mu=\left\|\frac{\partial^2 u}{\partial x^2}(0)\right\|_\infty$ and
$\beta=\|u(0)\|_\infty$, which is an integrable equation.  As a result,
\begin{equation}
\label{w_1_bnd}
\int w(t,x) \left(-h(x)\right) dx \ge \beta -\mu t - 2 \beta \int_0^t
\int \int f(x) \frac{1}{\sqrt{4 \pi s}} 
e^{-\frac{(x-y)^2}{4s}}w_0(y) dy\, dx\, ds.
\end{equation}
On the other hand using our choice for $w_0$,
\begin{eqnarray*}
\int_0^t \int \int & f(x)& \frac{1}{\sqrt{4 \pi s}}
e^{-\frac{(x-y)^2}{4s}}w_0(y) dy\, dx\, ds = \int_0^t \int f(x) \frac{1}{\sqrt{4 \pi s}}
e^{-\frac{(x-x_0)^2}{4s}} dx\, ds\\
&\le& \int_0^t \frac{1}{\sqrt{4\pi s}} \left(\gamma
\int_{-\infty}^{x_1} e^{-\frac{(x-x_0)^2}{4s}} dx + \|f\|_\infty
\int_{x_1}^\infty e^{-\frac{(x-x_0)^2}{4s}} dx \right) ds\\
&\le& \frac{\gamma \sqrt{t}}{4} + \frac{1}{2}\|f\|_\infty \int_0^t 1 -
\text{erf}\left( \frac{x_1-x_0}{2\sqrt{s}}\right) ds\\
&\le& \frac{\gamma \sqrt{t}}{4} + \frac{1}{2}\|f\|_\infty \int_0^t 1 -
\text{erf}\left( \frac{x_1-x_0}{2\sqrt{t}}\right) ds\\
&\le& \frac{\gamma \sqrt{t}}{4} + \frac{1}{2}t \|f\|_\infty \left( 1 - \text{erf}
\left( \frac{x_1-x_0}{2\sqrt{t}} \right ) \right )\\
&\le& \frac{3 \gamma \sqrt{t}}{4} \le \gamma \sqrt{t},\\
\end{eqnarray*}
we have used \eqref{vio_2}, \eqref{vio_3}, and assumed that $0 < t <
\frac{1}{4 \gamma^2}$.  Then \eqref{fence_2_improved} becomes
\begin{equation*}
1 \ge t \int w(t,x) \left(-h(x)\right) dx \ge \beta t - \mu t^2 - 2
\beta \gamma t \sqrt{t} = -2 \beta^{5/2} t^2 - \frac{2 \beta^{3/2}
  t^{3/2}}{\sqrt{27 K}} + \beta t,
\end{equation*}
using our choices of $\mu$, $\gamma$, and initial condition.  Maple
reports that the maximum of $A(t)=-2 \beta^{5/2} t^2 - \frac{2 \beta^{3/2}
  t^{3/2}}{\sqrt{27 K}} + \beta t$ is unique, occurs at $0 <
t_0 < \frac{1}{4 \gamma^2}$, and has the asymptotic expansion
\begin{equation*}
A(t_0) \sim K - 18 K \sqrt{\beta} + 432 K^3 \beta + O(\beta^{3/2}).
\end{equation*}
Thus for all small enough $\epsilon > \beta$, we obtain a
contradiction to \eqref{fence_2_improved} since $K>1$.  Thus, for some
$T<t_0<\infty$, $\lim_{t\to T^-} \|u(t)\|_\infty = \infty$.
\end{proof}
\end{thm}

\section{Discussion}
Theorem \ref{fence_2_improved_violation} gives a fairly strong
instability result.  No matter how small an initial condition to
\eqref{pde_instability} is chosen, even with all $p$-norms chosen small, solutions
can blow up so quickly that they fail to exist for all $t$.  This
precludes any kind of stability for classical solutions.  Like the
analogous result in Fujita's paper, the kind of initial conditions
which can be responsible for blow up are of the nicest kind imaginable
-- gaussians in either case!  

It must be understood that the argument in Theorem
\ref{fence_2_improved_violation} depends crucially on the decay of
$f$.  Without it, the lower bound on $\int w(t,x) (-h(x)) dx$
decreases too quickly.  Indeed, if $f=\text{const}>0$ and $h(x)>-f$,
then the comparison principle demonstrates that the zero function is
asymptotically stable.  On the other hand, any rate of decay for $f$
satisfies the hypotheses of Theorem \ref{fence_2_improved_violation},
and so will cause \eqref{pde_instability} to exhibit instability.

Finally, although we have examined the case where the nonlinearity in
\eqref{pde_instability} is due to $u^2$, there is no obstruction to extending the
analysis to any nonlinearity like $|u|^k$, with degree $k$ greater than
2.  A higher-degree nonlinearity would result in a somewhat different
form for \eqref{fence_2_eqn}, but this presents no further
difficulties to the argument.  Indeed, by analogy with Fujita's work,
higher-degree nonlinearities would result in significantly faster
blow-up.

\chapter{Cell complex structure for the space of heteroclines}
\label{unstable_ch}
\section{Introduction}

In this chapter, we determine that all unstable manifolds of \eqref{limited_pde}
\begin{equation}
\label{limited_pde_unstable}
\frac{\partial u(t,x)}{\partial t}=\frac{\partial^2 u(t,x)}{\partial x^2} - u^2(t,x)
+ \phi(x)
\end{equation}
are finite dimensional.  This is not a particularly new result, indeed
Theorem 5.2.1 in \cite{Henry} can easily be made to apply with the
Banach spaces we shall choose.  Theorem 5.2.1 in \cite{Henry} shows
the existence of a smooth finite dimensional unstable manifold locally
at an equilibrium.  One can then use the iterated time-1 map of the
flow for \eqref{limited_pde_unstable} to extend this local manifold to
a maximal unstable manifold.  There are also finite Hausdorff
dimensional attractors for the forward Cauchy problem on bounded
domains \cite{JRobinson_2001}.  However, we shall exhibit a more
global approach to the finite dimensionality of the unstable
manifolds.  This approach allows us to examine the finite
dimensionality of the space of heteroclinic orbits connecting a pair
of equilibria, which is a new result in the spirit of
\cite{Floer_gradient}.  The techniques used here depend rather
delicately on both the degree of the nonlinearity (quadratic) and the
spatial dimension (1).  Both of these are important in the standard
methodology as well, as the portion of the spectrum of the
linearization in the right half-plane needs to be bounded away from
zero.  In the case of \eqref{limited_pde_unstable}, the spectrum in
the right-half plane is discrete and consists of a finite number of
points.

\section{The linearization and its kernel}

We begin by considering an equilibrium solution $f$ to
\eqref{limited_pde_unstable}.  As discussed in Chapter \ref{nonauto_ch}, this
solution has asymptotic behavior which places it in $C^2 \cap L^1 \cap
L^\infty(\mathbb{R})$.  We are particularly interested in solutions
which lie in the $\alpha$-limit set of $f$, those solutions which
are defined for all $t<0$ and tend to $f$.  As in previous chapters,
center on this equilibrium by applying the change of variables
$u(t,x) \mapsto u(t,x) - f(x)$ to obtain
\begin{equation}
\label{pde_unstable}
\begin{cases}
\frac{\partial}{\partial t} u(t,x) = \frac{\partial^2}{\partial x^2}
u(t,x) - 2 f(x) u(t,x) - u^2(t,x)\\ u(0,x)=h(x) \in
C^2(\mathbb{R})\\
\lim_{t\to -\infty} u(t,x) = f(x)\\
 t<0,x\in \mathbb{R}.\\
\end{cases}
\end{equation}
Thus we have a final value problem for our nonlinear equation.  All
solutions to \eqref{pde_unstable} will tend to zero as $t\to -\infty$
uniformly by Theorem \ref{unif_limits_to_equilibria}.  Of course,
\eqref{pde_unstable} is ill-posed.  We show that there is only a
finite dimensional manifold of choices of $h$ for which a solution
exists.

\subsection{Backward time decay}

The decay of solutions to zero is a crucial part of the analysis, as
it provides the ability to perform Laplace transforms.  In the forward
time direction, one obtains upper bounds for solutions by way of
maximum principles, and lower bounds for the upper bounds by way of
Harnack estimates.  In the backward time direction, these tools
reverse roles.  Harnack estimates provide upper bounds, while the maximum
principle provides lower bounds for the upper bound.  In the proof of
Theorem \ref{unif_limits_to_equilibria}, the latter was used to some
advantage.  In this section, we briefly apply a standard Harnack
estimate to obtain an exponentially decaying upper bound.

Harnack estimates for a very general class of parabolic equations are
discussed in \cite{Kurihara_1967} and \cite{Aronson_1967}.  In those articles, the authors
examine positive solutions to
\begin{equation*}
\text{div }{\bf A}(x,t,u,\nabla u) - \frac{\partial u}{\partial t} = B(x,t,u,\nabla u),
\end{equation*}
where $x\in\mathbb{R}^n$, and ${\bf
  A}:\mathbb{R}^{2n+2}\to\mathbb{R}^n$ and $B:\mathbb{R}^{2n+2}\to
  \mathbb{R}$ satisfy
\begin{eqnarray*}
|{\bf A}(x,t,u,p)| &\le& a|p| + c|u| + e\\
|B(x,t,u,p)| &\le& b|p| + d|u| + f\\
p \cdot {\bf A}(x,t,u,p) &\ge& \frac{1}{a} |p|^2 - d|u|^2 - g,\\
\end{eqnarray*}
for some $a>0$ and $b,...g$ are measurable functions.  For a solution
$u$ defined on a rectangle $R$, the authors define a pair of
congruent, disjoint closed rectangles $R^+,R^- \subset R$ with $R^-$
being a backward time translation of $R^+$.  The main result is the
Harnack inequality
\begin{equation}
\label{harnack_eqn}
\max_{R^-} u \le \gamma \left( \min_{R^+}u + L \right ),
\end{equation}
where $\gamma>0$ depends only on geometry and $a$ (but not $b,...g$)
and $L$ is a linear combination of $e,f,g$ whose coefficients depend
on geometry.  

In the case of \eqref{pde_unstable}, or indeed of the analogous
equation with higher degree terms, we have that \eqref{harnack_eqn}
will apply with $L=0$.  Notice that the conditions on $A,B$ are
satisfied because any solution to \eqref{pde_unstable} is
automatically a finite energy solution, and therefore is bounded and
has bounded first derivatives.  The only difficulty is that
\eqref{harnack_eqn} applies for {\it positive} solutions, while
\eqref{pde_unstable} may have solutions with negative portions.  However, one can
pose the problem for the (weak) solution of
\begin{eqnarray*}
\frac{\partial |u|}{\partial t} &=& \text{sgn } (u) \left(\Delta u - u^2 - 2 f u\right)\\
&=&\Delta |u| - u|u| - 2 f |u|\\
&\ge& \Delta |u| - |u|^2 - 2 |f| |u|\\
\end{eqnarray*}
for which we only get positive solutions.  By iterating
\eqref{harnack_eqn} we have that solutions to \eqref{pde_unstable} decay
exponentially as $t\to -\infty$.  

\subsection{Topological considerations}

\begin{df}
Let $Y_a(X)$ be the subspace of $C^1(X,C^{0,\alpha}(\mathbb{R}))$ which
consists of functions which decay exponentially to zero like $e^{at}$,
where $0<\alpha \le 1$.  We define the weighted norm
\begin{equation*}
\|u\|_{Y_a} = \left\|e^{-at}\|u(t)\|_{C^{0,\alpha}(\mathbb{R})}\right\|_{C^1} 
\end{equation*}
 and the space
\begin{equation*}
Y_a(X)=\left\{u=u(t,x)\in C^1(X,C^{0,\alpha}(\mathbb{R}))|
\|u\|_{Y_a} < \infty \right\}.
\end{equation*}
In a similar way, we can define the weighted Banach space $Z_a(X)$ as a
subspace of $C^0(X,C^{0,\alpha}(\mathbb{R}))$.  It is quite important
that $Y_a$ and $Z_a$ are Banach algebras under pointwise
multiplication.
\end{df}

In light of the previous section, solutions to \eqref{pde_unstable}
are zeros of the densely defined nonlinear operator $N:Y_a((-\infty,0]) \to Z_a((-\infty,0])$ given by
\begin{equation}
\label{N_def}
N(u)=\frac{\partial u}{\partial t} - \frac{\partial^2 u}{\partial x^2}
+ u^2 + 2fu.
\end{equation}
About the zero function, the linearization of $N$ is the densely
defined linear map $L:Y_a((-\infty,0]) \to Z_a((-\infty,0])$ given by
\begin{equation}
\label{L_def}
L = \frac{\partial }{\partial t} - \frac{\partial^2 }{\partial x^2}
+ 2 f = \frac{\partial }{\partial t} - H,
\end{equation}
where we define $H=\frac{\partial^2 }{\partial x^2} - 2 f$.  Also note
that $L$ is the Frech\'et derivative of $N$, which follows from the
fact that $Y_a$ and $Z_a$ are Banach algebras.

\begin{rem}
We are using $C^{0,\alpha}(\mathbb{R})$ instead of $C^0(\mathbb{R})$
to ensure that $N$ and $L$ be densely defined.  We could use space of
continous functions which decay to zero, or the space of uniformly
continous functions equally well.
\end{rem}

\begin{conv}
We shall conventionally take $a>0$ to be smaller than the smallest
eigenvalue of $H$.
\end{conv}

We show two things: that the kernel of $L$ is finite
dimensional, and that $L$ is surjective.  These two facts enable us to
use the implicit function theorem to conclude that the space of
solutions comprising the $\alpha$-limit set of an equilibrium is a
finite dimensional submanifold of $Y_a((-\infty,0])$.

\subsection{Dimension of the kernel}

\begin{lem}
\label{eq_findim}
If $f$ is an equilibrium solution, then the operator
$L:Y_a((-\infty,0])\to
Z_a((-\infty,0])$ in \eqref{L_def} has a finite
dimensional kernel.
\begin{proof}
Notice that the operator $L$ is separable, so we try the usual
separation $h(t,x)=T(t)X(x)$.  Substituting into \eqref{L_def} gives
\begin{eqnarray*}
0&=&Lh=\left(\frac{\partial}{\partial t} - \frac{\partial^2}{\partial
  x^2} + 2 f \right)h\\
&=&T'X+T\left(-\frac{\partial^2}{\partial x^2} + 2f\right)X\\
\frac{T'}{T}&=&\frac{\left(\frac{\partial^2}{\partial x^2} - 2f\right)X}{X}=\lambda\\
\end{eqnarray*}
for some $\lambda\in\mathbb{C}$.  The separated equation for $T$
yields $T=C_x e^{\lambda t}$.  Since we are looking for the kernel of
$L$ in $Y_a \subset L^\infty(\mathbb{R}^2)$, we must conclude
that $\lambda$ must have nonnegative real part.  On the other hand,
the spectrum of $H=\left(\frac{\partial^2}{\partial x^2} - 2f\right)$
is strictly real, so $\lambda \ge 0$.  Indeed, there are finitely many
positive possibilities for $\lambda$ each with finite-dimensional
eigenspace.  This is a standard fact about the Schr\"{o}dinger
operator $H$ since $f$ is an equilibrium (Proposition
\ref{finitely_many_eigenvalues}).  Thus $L$ has a finite dimensional
kernel.
\end{proof}
\end{lem}

\subsection{Surjectivity of the linearization}

In order to show the surjectivity of $L$, we will construct a map
$\Gamma:Z_a((-\infty,0]) \to
  Y_a((-\infty,0])$ for which $L \circ \Gamma =
    \text{id}_{Z_a}$.  That is, we construct a right-inverse to $L$,
    noting of course that $L$ is typically not injective.
We shall derive a formula for $\Gamma$ using the Laplace
transform $v \mapsto \overline{v}$
\begin{equation*}
\overline{v}(s,x) = \int_{-\infty}^0 e^{st} v(t,x) dt,
\end{equation*}
where $\Re(s)>-a$ and $v\in Z_a((-\infty,0])$.

Since Lemma \ref{eq_findim} essentially solves \eqref{pde_unstable},
we will be solving the inhomogeneous problem with zero final condition
\begin{equation}
\label{eq_findim_zero}
\begin{cases}
\frac{\partial v(t,x)}{\partial t}-\frac{\partial^2 v(t,x)}{\partial x^2} + 2
f(x) v(t,x) = -w(t,x) \in Z_a((-\infty,0])\\
v(0,x)=0\\
\end{cases}
\end{equation}
for $t<0$.  The Laplace transform of this problem is 
\begin{eqnarray*}
s\overline{v}(s,x)+\frac{\partial^2 \overline{v}(s,x)}{\partial x^2} - 2 f(x)
\overline{v}(s,x)&=&\overline{w}(s,x)\\
(H+s)\overline{v}(s,x)&=&\overline{w}(s,x).
\end{eqnarray*}
Choose a vertical contour $C$ with $0>\Re(s) > -a$, so that the
Laplace transforms are well-defined, and that the contour remains
entirely in the resolvent set of $-H$.  Then we can invert to obtain
\begin{equation*}
\overline{v}(s,x)=(H+s)^{-1}\overline{w}(s,x).
\end{equation*}
Using the inversion formula for the Laplace transform yields
\begin{eqnarray*}
v(t,x)&=& \frac{1}{2\pi i} \int_C e^{-st} (H+s)^{-1} \overline{w}(s,x)
ds\\
&=&\frac{1}{2\pi i} \int_C e^{-st} (H+s)^{-1} \int_t^0
e^{s\tau}w(\tau,x)d\tau\,ds\\
&=&\int_t^0 \left(\frac{1}{2\pi i} \int_C e^{s(\tau-t)} (H+s)^{-1} 
ds \right)w(\tau,x)d\tau.\\
\end{eqnarray*}

\begin{figure}
\begin{center}
\includegraphics[height=2.5in]{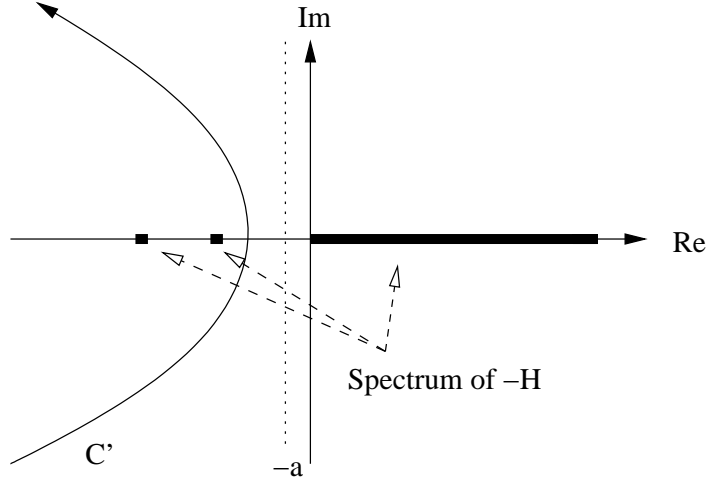}
\end{center}
\caption{Definition of the contour $C'$}
\label{Cp_contour_fig}
\end{figure}

We can obtain operator convergence of the operator-valued integral
in parentheses if we deflect the contour $C$.  Choose instead the
portion $C'$ of the hyperbola (See Figure \ref{Cp_contour_fig}) 
\begin{equation}
\left(\Re(s)\right)^2-\left(\Im(s)\right)^2 = \frac{1}{4}(\lambda-a)^2
\end{equation}
(where $\lambda$ is the smallest magnitude eigenvalue of $-H$) which
lies in the left half-plane as our new contour.  Then, since $-H:C^{0,\alpha}
\to C^{0,\alpha}$ is sectorial about $(\lambda-a)/2$ (Proposition
\ref{sectorial_H_prop}), Theorem 1.3.4 in \cite{Henry} implies that
the integral
\begin{equation*}
\left(\frac{1}{2\pi i} \int_{C'} e^{s(\tau-t)} (H+s)^{-1} 
ds \right)
\end{equation*}
defines an operator-valued semigroup $e^{-H(\tau-t)}$, so the formula
  for $\Gamma$ is given by
\begin{equation}
\Gamma(w)(t,x)= \int_t^0 e^{-H(\tau-t)} w(\tau,x)d\tau.
\end{equation}
It remains to show that the image of $\Gamma$ is in fact $Y_a$,
as it is easy to see that its image is in $L^\infty$.  That the image
is as advertised is not immediately obvious because the contour
deflection $C \to C'$ changes the domain of the Laplace transform.  In
particular, the derivation given above is no longer valid with the
new contour.

Therefore, we must estimate $\|v\|_{Z_a}$ (recall that $\lambda$ is the smallest magnitude eigenvalue of $-H$)

\begin{eqnarray*}
\|e^{-at} v(t,x)\|_{C^0} &=& \left \| \frac{1}{2\pi i}\int_{C'} (s+H)^{-1}
 \int_t^0 e^{-(s+a)(t-\tau)} e^{a\tau}w(\tau,x) d\tau\,ds \right \|_{C^0}\\
&\le&\frac{1}{2\pi}\int_{C'} \frac{K_1}{|s-\lambda|}e^{-\Re(s+a)t}
 \int_t^0 e^{\Re(s+a)\tau} \|w\|_{Z_a} d\tau\,ds\\
&\le&\frac{K_1\|w\|_{Z_a}}{2\pi}\int_{C'} \frac{1}{|s-\lambda|}e^{-\Re(s+a)t}
 \frac{1}{\Re(s+a)}\left(1-e^{\Re(s+a)t}\right)ds\\
&\le&\frac{K_1\|w\|_{Z_a}}{\pi}\int_{C'}
 \frac{ds}{|s-\lambda||\Re(s+a)|}\\
&\le&K_2 \|w\|_{Z_a},\\
\end{eqnarray*}
where $0<K_1,K_2<\infty$ are independent of $t$ and $w$.  We have made
use of the estimate in Proposition \ref{sectorial_H_prop} of the norm
of $(H+s)^{-1}:C^{0,\alpha} \to C^{0,\alpha}$ when $s$ is in the resolvent set of
$-H$.  In particular, note that the choice of $C'$ being to the left
of $-a$ is crucial to the convergence of the integrals.  Thus the
image of $\Gamma$ lies in $Z_a$.  The backward-time decay of
$\frac{\partial v}{\partial t}$ is immediate from the Harnack
inequality, so in fact the image of $\Gamma$ lies in $Y_a$.

\begin{thm}
\label{unstable_thm}
The linear map $L:Y_a((-\infty,0])\to
  Z_a((-\infty,0])$ is surjective and has a
  finite dimensional kernel.  Therefore the set $N^{-1}(0)$ is a
  finite dimensional manifold, which is the unstable manifold of the
  equilibrium $f$.  The dimension of $N^{-1}(0)$ is precisely the
  dimension of the positive eigenspace of $H$.
\begin{proof}
The only thing which remains to be shown is that the domain $Y_a$
splits into a pair of closed complementary subspaces: the kernel of
$L$ and its complement.  That its complement is closed follows
immediately from a standard application of the Hahn-Banach theorem.
(Extend $\text{id}_{\text{ker }L}$ to all of $Y_a$.)
\end{proof}
\end{thm}

Combining the fact that an equilibrium solution can have an
empty unstable manifold (we numerically computed the dimension of
the eigenspaces of $L$ in Chapter \ref{nonauto_ch}) and
is yet unstable, we have proven the following result.

\begin{thm}
All equilbrium solutions to \eqref{limited_pde_unstable} are degenerate
critical points in the sense of Morse.
\end{thm}

\section{Linearization about heteroclinic orbits}

We can extend the technique of the previous section to the
linearization about a heteroclinic orbit.  The resulting
generalization of Theorem \ref{unstable_thm} is that the connecting
manifolds of \eqref{limited_pde_unstable} are all finite dimensional.

Suppose that $u$ is a heteroclinic orbit of \eqref{limited_pde_unstable}.  
Let $f_-,f_+$ be the equilibrium solutions of \eqref{pde} to which
$u$ converges as $t\to -\infty$ and $t\to +\infty$ respectively.  

Suppose that $\lambda_0:\mathbb{R} \to (0,\infty)$ is the smallest
positive eigenvalue of $H(t)$.  It is easy to see that $\lambda_0$ is
piecewise $C^1$, for instance, see Proposition I.7.2 in
\cite{Kielhoefer_2004}.  Propostion \ref{spectral_control} ensures
that $\lambda_0$ is a bounded function.  We will define a pair of
bounded, piecewise $C^1$ functions $\lambda_1$ and $\lambda_2$ which
will aid us in defining a two more pairs of function spaces.  Let
$\lambda_1:\mathbb{R} \to (0,\infty)$ be a bounded, piecewise $C^1$
function with bounded derivative which has the following properties:
\begin{itemize}
\item $\lambda_1(t)$ is never an eigenvalue of $H(t)$, 
\item $\lim_{t\to\infty}\frac{\lambda_1(t)}{\lambda_0(t)}<1$,
\item $\lim_{t\to -\infty}\frac{\lambda_1(t)}{\lambda_0(t)}<1$, and
\item since $u\to
f_\pm$ uniformly, for a sufficiently large $R>0$, $\lambda_1$ can be
chosen so that there are no jumps on its restriction to $\mathbb{R} -
[-R,R]$.
\end{itemize}

Defining $\lambda_2$ is a somewhat more delicate problem.  We would
like to exclude the solutions which lie in the unstable manifold of
$f_+$, since they cannot lie in the space of heteroclines from $f_-
\to f_+$.  We do this by separating the eigenvalues corresponding to
the intersection of the unstable manifolds of $f_-$ and $f_+$ from those
which lie in the stable manifold of $f_+$.  However, there is an
obstruction to this technique.  In particular, the eigenvalues of
$H(t)=\frac{\partial^2}{\partial x^2} - 2 u(t)$ vary with time, and
can bifurcate.  To avoid this issue, we need some kind of regularity
for the eigenvalues to prevent them from bifurcating.  We follow Floer
\cite{Floer_relative} in the following way:

\begin{conj}
\label{baire_conj}
There is a generic subset (a Baire subset) of choices for the
coefficients $a_i$ in \eqref{pde} so that if $u$ is a heteroclinic
orbit, all of the eigenvalues of $H(t)$ are simple.  
\end{conj}

Numerical evidence, as exhibited in Chapters \ref{nonauto_ch} and
\ref{example_ch} suggests that the above Conjecture is true.  When we
assume that all of the eigenvalues of $H(t)$ are simple, and therefore
do not undergo any bifurcations other than passing through zero, we
shall say $u$ is a heterocline contained in $U_{reg}$.

\begin{figure}
\begin{center}
\includegraphics[height=2.5in]{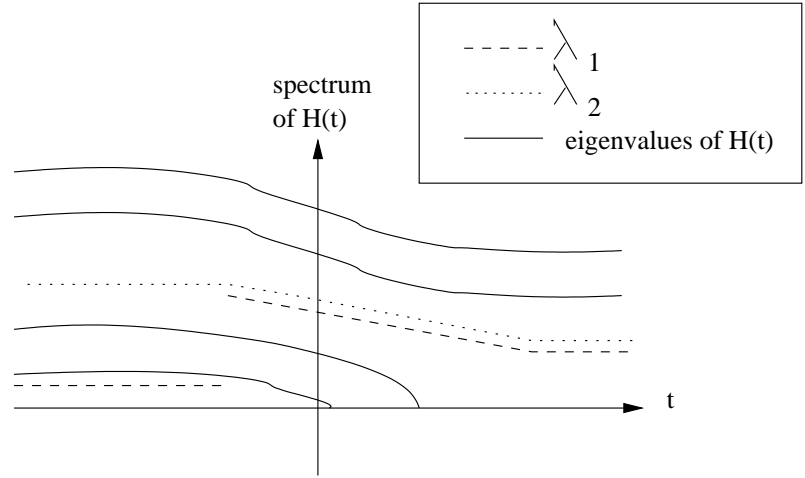}
\end{center}
\caption{Definition of $\lambda_1$ and $\lambda_2$}
\label{lambda_defs_fig}
\end{figure}

Let $\lambda_2$ be in $C^1(\mathbb{R})$ such that
\begin{itemize}
\item $\lambda_2=\lambda_1$ on $[R,\infty)$, and
\item $\lambda_2(t)$ is not an eigenvalue of
$H(t)$ for any $t$.
\end{itemize}
We can do this when $u \in U_{reg}$.  See Figure
\ref{lambda_defs_fig}.

\begin{df}
Define the Banach algebra $Y_{\lambda_i}(X)$ (for $i=1,2$) to
be the set of $u$ in $C^1(X,C^{0,\alpha}(\mathbb{R}))$ such that the
norm
\begin{equation*}
\left\|e^{-\int_0^t \lambda_i(\tau) d\tau} \|u(t)\|_{C^{0,\alpha}}
\right\|_{C^1}<\infty,
\end{equation*}
where $X$ is an interval containing zero.  Likewise, we can define the
spaces $Z_{\lambda_i}(X) \subset
C^0(X,C^{0,\alpha}(\mathbb{R}))$ in a similar way.  That these are
Banach spaces follows from the boundedness of the $\lambda_i$.  It is
also elementary to see that these are Banach algebras.
\end{df}

We then consider $N_i,L_i$ as $Y_{\lambda_i}(\mathbb{R})\to
Z_{\lambda_i}(\mathbb{R})$, where $L_i$ is the linearization of
$N_i$ about $u$ for $i=1,2$.  (Again, since $Y_{\lambda_i}$ and
$Z_{\lambda_i}$ are Banach algebras, $L_i$ is the Frech\'et
derivative of $N_i$.)  For a $i\in \{1,2\}$, consider the
restriction $L^-_i$ of $L_i$ to a map
$Y_{\lambda_i}((-\infty,0])\to
Z_{\lambda_i}((-\infty,0])$.  We rewrite
\begin{equation}
\label{lminus}
L^-_i=\left(\frac{\partial}{\partial t} - \frac{\partial^2}{\partial x^2}+2f_- \right)+(2f_- - 2u).
\end{equation}
Likewise, we can define $L^+_i:Y_{\lambda_i}([0,\infty))\to Z_{\lambda_i}([0,\infty))$.

We define the positive eigenspaces $V^+$ for the equilibria as well
\begin{equation}
V^+(f_\pm)=\text{span }\left \{v\in C^{0,\alpha}(\mathbb{R})| \text{ there
  is a }\lambda>0 \text{ with } \left(\frac{\partial^2
  }{\partial x^2} -2f_\pm \right)v = \lambda v \right \}.
\end{equation}
Note in particular that $\text{dim }V^+(f_\pm)<\infty$.

\begin{lem}
\label{findim_reqs}
If $u\in U_{reg}$ is a heterocline that converges to $f_\pm$ as
$t\to\pm\infty$, then the operator $L_i$ has a finite
dimensional kernel for $i\in\{1,2\}$, and in particular
\begin{equation*}
\lim_{t\to -\infty}\text{dim }V^+(u(t)) - \lim_{t\to +\infty}\text{dim }V^+(u(t))  \le \text{dim ker }L_i
\le\text{dim ker }L^-_i <\infty.
\end{equation*}
(The condition $u\in U_{reg}$ is only necessary for the $i=2$ case.)
\begin{proof}
Notice that the first term of \eqref{lminus} has finite dimensional
kernel by Lemma \ref{eq_findim} and closed image by Theorem \ref{unstable_thm}.
The second term of \eqref{lminus} is a compact operator since $u\to f_-$ uniformly.
Thus $L^-_i$ has a finite dimensional kernel.  Let $\text{span}\{v_m\}_{m=1}^{M}=\ker L^-_i$
and consider the set of Cauchy problems 
\begin{equation}
\label{cpset}
\begin{cases}
\frac{\partial h}{\partial t}=\frac{\partial^2 h}{\partial x^2}-2u h \text{ for }t>0\\
h(0,x)=v_m(0,x).
\end{cases}
\end{equation}
Standard parabolic theory gives uniqueness of solutions to \eqref{cpset}, and that
a solution $h$ lies in the kernel of $L^+_i$, the restriction of $L_i$ to $[0,\infty)\times\mathbb{R}$.
Therefore $\text{dim ker }L_i \le \text{dim ker }L^-_i<\infty$.  

For the other inequality, modify $u$ outside of
$[-R,R]\times\mathbb{R}$ to get a $\bar{u}$ so that the linearization $\overline{L_i}$
of $N$ about $\bar{u}$ satisfies
\begin{itemize}
\item $\text{ker }\overline{L_i}$ is isomorphic to $\text{ker }L_i$ as vector
  spaces,
\item $\bar{u}|_{(-\infty,-R)\times\mathbb{R}} = f_-$, and
\item $\bar{u}|_{(R,\infty)\times\mathbb{R}} = f_+$.
\end{itemize}
We can do this for a sufficiently large $R$, since $u$ tends uniformly
to equilibria.  Then the flow of
\begin{equation*}
\frac{\partial h}{\partial t}=\frac{\partial^2 h}{\partial x^2} +
2\bar{u} h
\end{equation*}
defines an injective linear map from the timeslice at $-R$ to the
timeslice at $R$.  (That is, it gives an injective map from
$C^{0,\alpha}(\mathbb{R})$ to itself -- injectivity being an expression
of the uniqueness of solutions.)  Each element $v$ of the kernel
of $\overline{L_i}$ evidently must have $v(-R)\in V^+(f_-)$ and
$v(R) \notin V^+(f_+)$.  Therefore, the injectivity ensures that the
intersection of the image under the flow of $V^+(f_-)$ with the
complement of $V^-(f_+)$ has at least dimension $\text{dim }V^+(f_-) -
\text{dim }V^+(f_+)$. 
\end{proof}
\end{lem}

\begin{rem}
Multiplication by $u$, $C^1(\mathbb{R}^2,C^{0,\alpha}(\mathbb{R}))\to C^0(\mathbb{R}^2)$ is not a compact operator, in particular note that $\text{dim ker }L^+_i=\infty$.
\end{rem}

\begin{thm}
\label{finite_dim_connecting_thm}
Let $u$ be a heterocline of \eqref{limited_pde_unstable}
which connects equilibria $f_\pm$.  There exists a union $\bigcup M_u$ of finite dimensional
submanifolds $M_u$ of $C^1(\mathbb{R},C^{0,\alpha}(\mathbb{R}))$ which
\begin{itemize}
\item contains $u$ and
\item consists of heteroclines connecting $f_-$ to $f_+$.
\end{itemize}
If $u\in U_{reg}$, then $M_u$ has dimension $\lim_{t\to -\infty}\text{dim }V^+(u(t)) - \lim_{t\to
  \infty}\text{dim }V^+(u(t))$, and this is maximal among such
  submanifolds $M_u$.
\begin{proof}
Observe that $L_1$ is surjective, since it is easy to show that the
formula
\begin{equation*}
\Gamma_1(w)(t)=\int_t^0 e^{-\int_0^{T-t} H(\tau) d\tau} w(T,x) dT
\end{equation*}
is a well defined right inverse of $L_1$.  This involves showing that 
\begin{equation*}
e^{-\int_0^{t} H(\tau) d\tau}=\frac{1}{2\pi i} \int_{C(t)} e^{st}
(H(t)+s)^{-1} ds
\end{equation*}
converges, where we note that the contour changes with time.  As it
happens, the computation in \cite{Henry} goes through with the only
change that at $t=0$, we deflect the contour to the right, rather than
the left (as in Figure \ref{Cp_contour_fig}).  Since Lemma
\ref{findim_reqs} shows that $L_1$ has finite dimensional kernel, then
it follows that $M_u=N_1^{-1}(0)$ is a union of finite dimensional
manifolds, with a finite maximal dimension.  It is obvious that $M_u$
consists entirely of heteroclinic orbits and contains $u$.  

It remains to show that the dimension of $M_u$ is as advertised and
maximal.  Observe that $L_2$ is a compact perturbation of an operator
$L'_2:Y_{\lambda_2}(\mathbb{R})\to
Z_{\lambda_2}(\mathbb{R})$ which is time-translation invariant.
This follows from the precise choice of $\lambda_2$ being continous
and not intersecting the eigenvalues of $H$.  $L_2$ and $L'_2$ are
both surjective by exactly the same reasoning as for $L_1$.  $L'_2$ is
injective by using separation of variables as in Lemma \ref{eq_findim}
(noting that all nontrivial solutions blow up in the
$Y_{\lambda_2}$ norm).  Therefore the Fredholm index of $L'_2$,
hence $L_2$ is zero.  However, this implies that $L_2$ is injective.

Since $L_2$ is bijective, any solution to $L_2 u = 0$ which decays
faster than $e^{\int \lambda_2(t) dt}$ as $t\to -\infty$ ends up
growing faster than $e^{\int \lambda_2(t) dt}$ as $t\to +\infty$, and
in particular does not tend to zero.  As a result, such a solution
cannot be in $\text{ker }L_1$.  This implies that $\text{dim ker }L_1
\le \lim_{t\to -\infty}\text{dim }V^+(u(t)) - \lim_{t\to
  \infty}\text{dim }V^+(u(t))$, which with the estimate in Lemma
\ref{findim_reqs} completes the proof.
\end{proof}
\end{thm}

\begin{rem}
Even if $u \notin U_{reg}$ (when there exist nonsimple eigenvalues of
$H(t)$), the function $\lambda_1$ can still be constructed.  As a
result, we {\it always} get that the connecting manifold $M_u$ is
finite-dimensional.
\end{rem}

\begin{cor}
\label{cell_complex_cor}
The space of heteroclinic orbits has the structure of a cell complex
with finite dimensional cells.  This cell complex structure is
evidently finite dimensional if there exist only finitely many
equilibria for \eqref{limited_pde_unstable}.
\end{cor}

\section{Conclusions}

We have shown that the tangent space at an equilibrium splits into a
finite dimensional unstable subspace, and infinite dimensional center
and stable subspaces.  However, it is quite clear by Chapter
\ref{instability_ch} that the center subspace is nonempty and large.
Indeed, considering the work of \cite{Souplet_2002}, the center and
stable subspaces are not closed complements of each other.
Additionally, we have given conditions for the space of heteroclinic
orbits to have a finite dimensional cell complex structure.

\chapter{An extended example}
\label{example_ch}
\section{Introduction}
Consider the following equation

\begin{equation}
\label{example_pde}
\frac{\partial u}{\partial t}=\frac{\partial^2 u}{\partial x^2} - u^2
+ (x^2-c)e^{-x^2/2},
\end{equation}
where the choice of $\phi$ in \eqref{limited_pde} has been fixed.  The
bifurcation diagram for the equilibria of \eqref{example_pde} can be
found in Figure \ref{bif_diag}.  Based on the Theorem
\ref{unstable_thm}, the number of positive eigenvalues shown in Figure
\ref{bif_diag} corresponds exactly to the dimension of the unstable
manifold of each equilibrium.

\section{Frontier of the stable manifold}

\begin{figure}
\begin{center}
\includegraphics[height=1.33in]{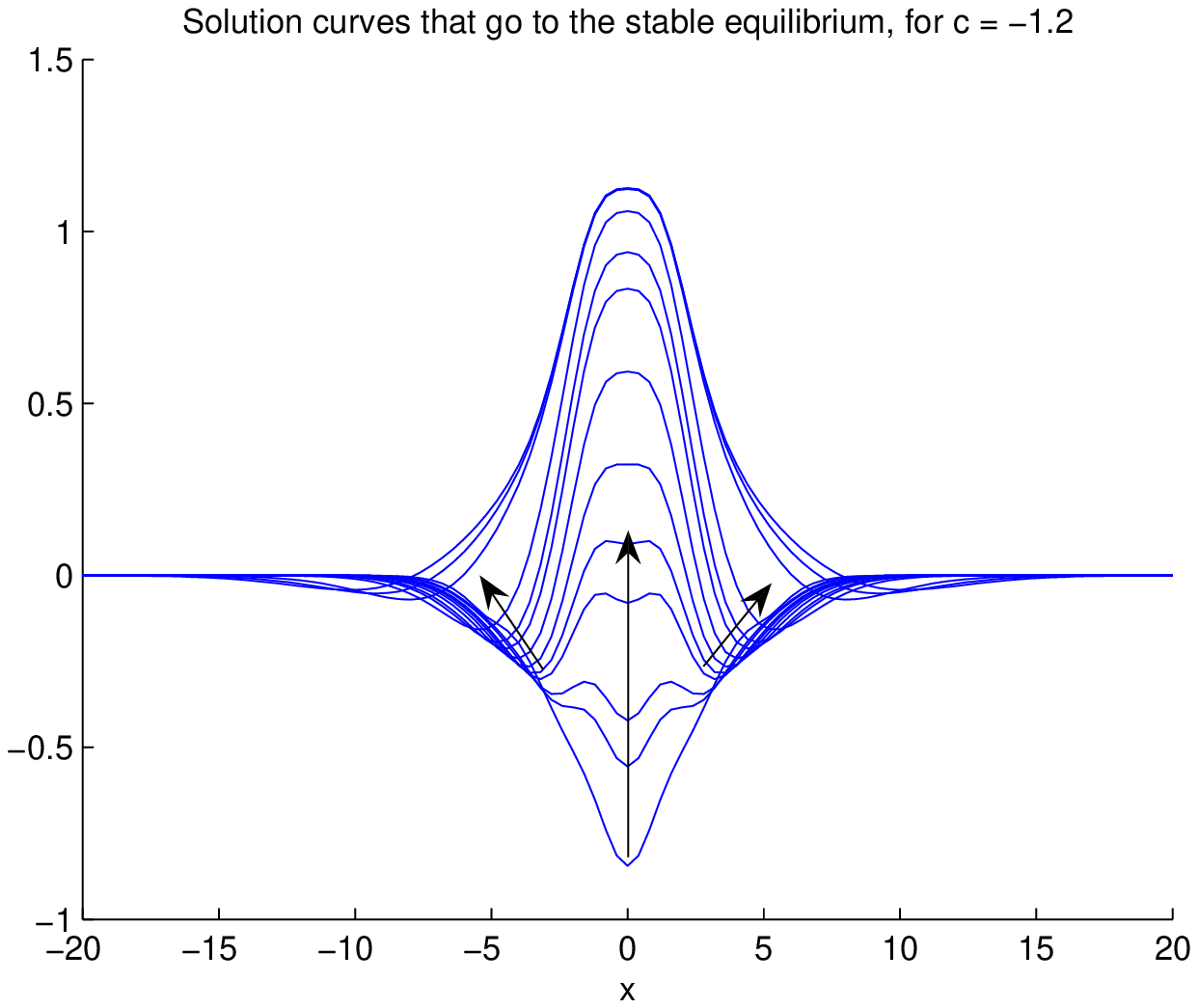}
\includegraphics[height=1.33in]{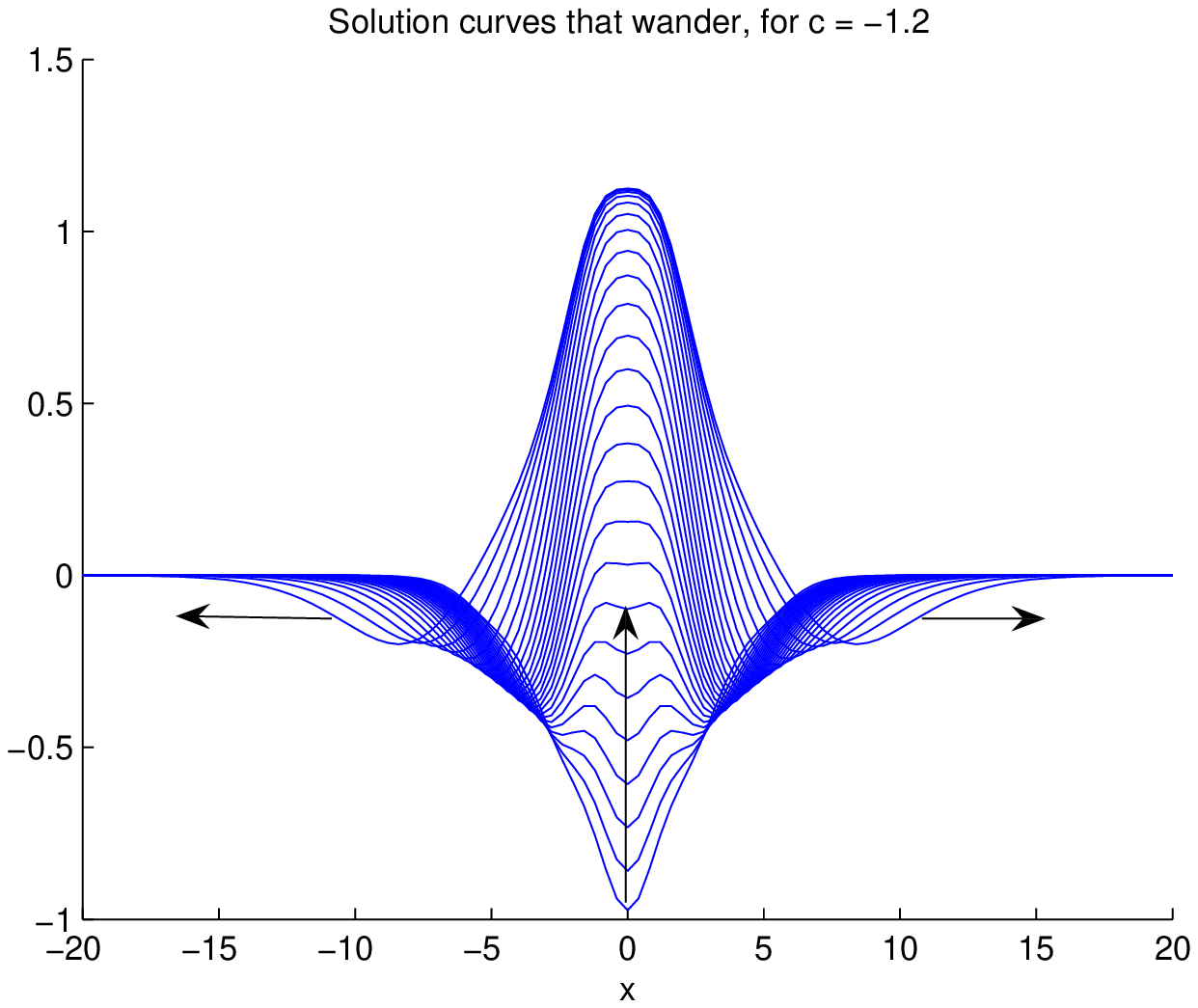}
\includegraphics[height=1.33in]{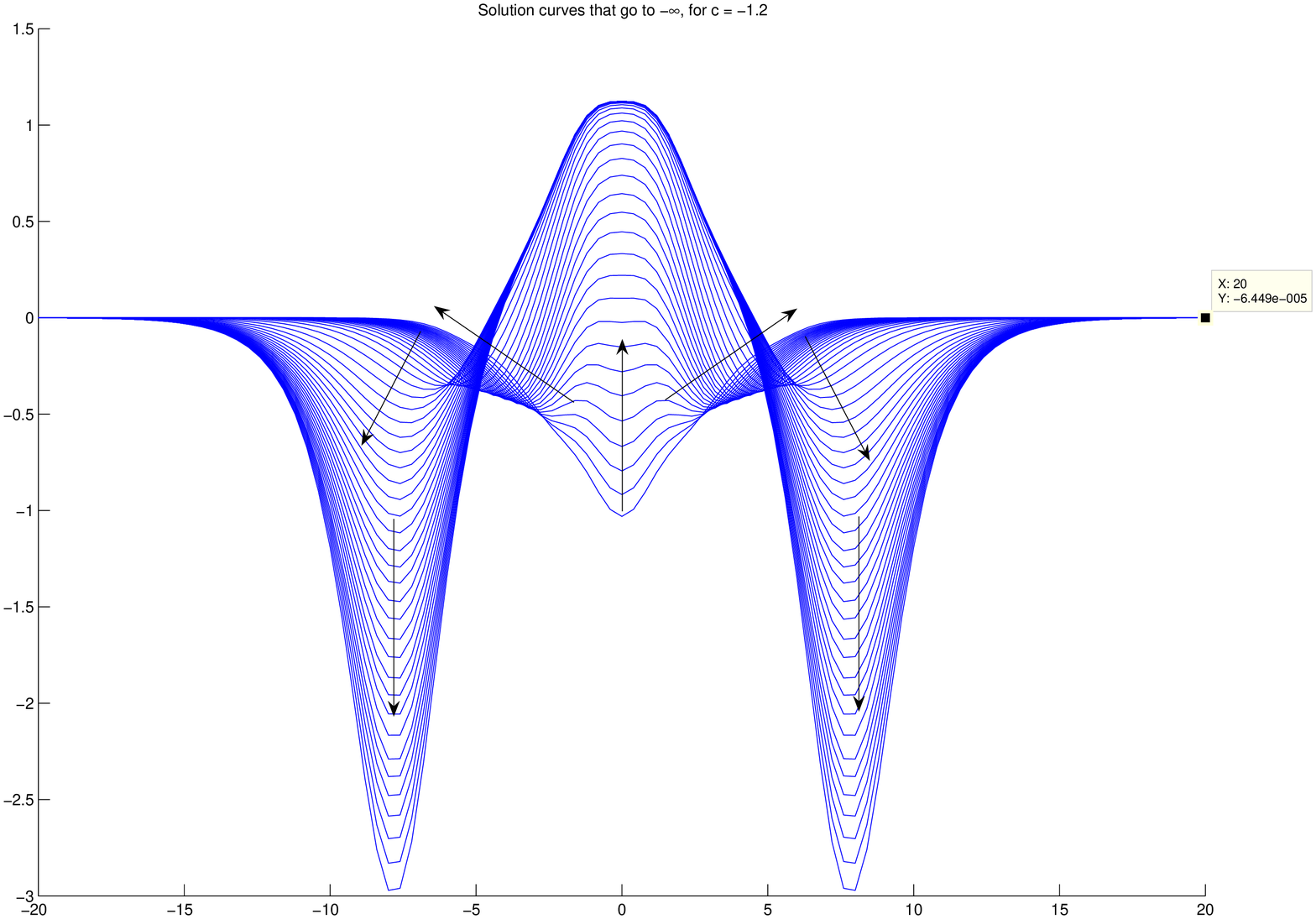}
\end{center}
\caption{Behavior of solutions near the frontier of the stable
  manifold of $f_0$ (horizontal axis is $x$)}
\label{frontier_fig}
\end{figure}

According to Figure \ref{bif_diag}, when $c=-1.2$, there is only one
equilibrium, $f_0$.  It has empty unstable manifold, though of course
it is asymptotically unstable (as is shown in Chapter
\ref{instability_ch}).  On the other hand, $f_0$ has an infinite
dimensional stable manifold, which is not all of
$C^{0,\alpha}(\mathbb{R})$, as a consequence of the asymptotic
instability.  As a result, its stable manifold has a frontier in
$C^{0,\alpha}(\mathbb{R})$ (which may not be a boundary in the sense
of a manifold with boundary).  We are interested in the qualitative
behavior of solutions near and along this frontier.  We know by
Theorem \ref{unif_limits_to_equilibria} that if they tend to $f_0$
uniformly on compact subsets, then they do so uniformly.  It is
enlightening to use a numerical procedure to this end.  We start
solutions at the following family of initial conditions
\begin{equation}
\label{start_eq}
u_A(x) = f_0(x) + A e^{-x^2/10}.
\end{equation}
Using the Fujita technique (exactly as shown in Chapter
\ref{instability_ch}), we can show that for sufficiently negative $A$,
the solution started at $u_A$ will not be eternal.  As a result, the
family of initial conditions $u_A$ intersects the frontier of the
stable manifold of $f_0$.  An approximation to the value of $A$ which
corresponds to the frontier can be easily found using a binary search.
Some typical such solutions are shown in Figure \ref{frontier_fig},
and the approximate value of $A$ corresponding to the frontier is
$A\approx -2.15$

The qualitative behavior shown in Figure \ref{frontier_fig} indicates
that there is some kind of traveling disturbance in the frontier
solutions, which seems like a traveling wave.  However, such a
solution also appears to tend uniformly on compact subsets to $f_0$,
so in fact it converges uniformly.  (The uniform convergence is not
obvious from the figure, due to the numerical solution being truncated
at a finite time.)  The leading edge of this disturbance collapses to
$-\infty$ in finite time for solutions just outside the stable
manifold of $f_0$.

\section{Flow near equilibria with two-dimensional unstable manifolds}

\begin{figure}
\begin{center}
\includegraphics[height=2.5in]{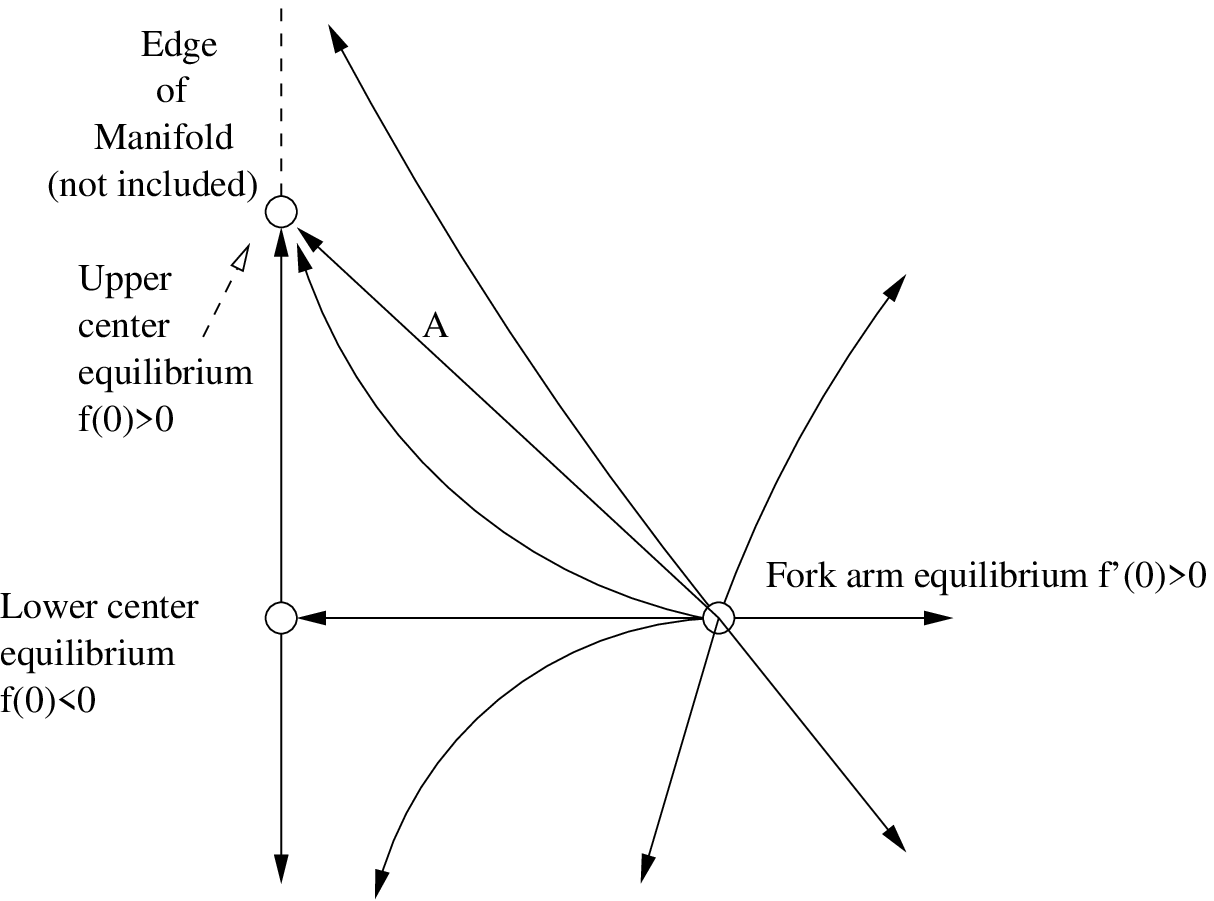}
\includegraphics[height=2.5in]{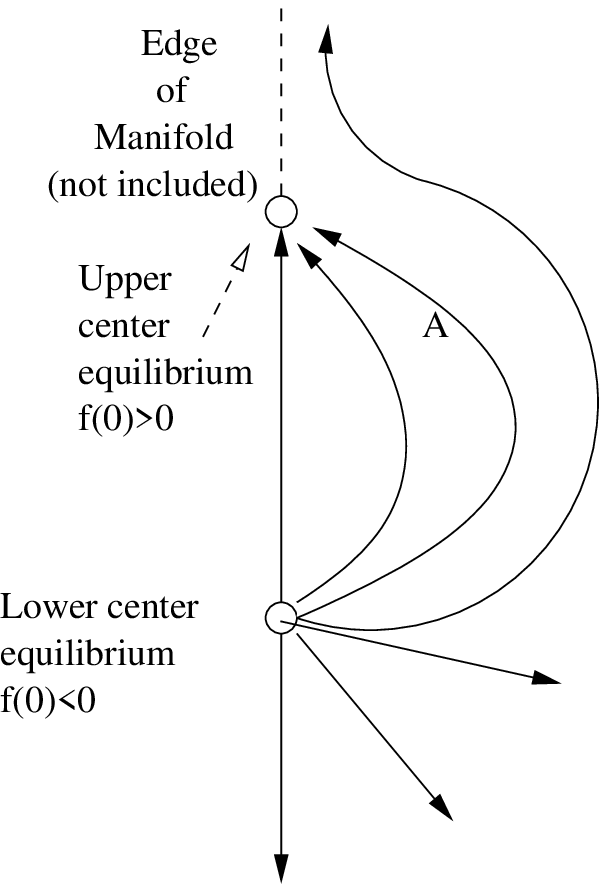}
\end{center}
\caption{Flow in the unstable manifold of a ``fork arm.''  $c=0.0600$
  (left); $c=0.0501$ (right)}
\label{flow_fig}
\end{figure}

\begin{figure}
\begin{center}
\includegraphics[height=2in]{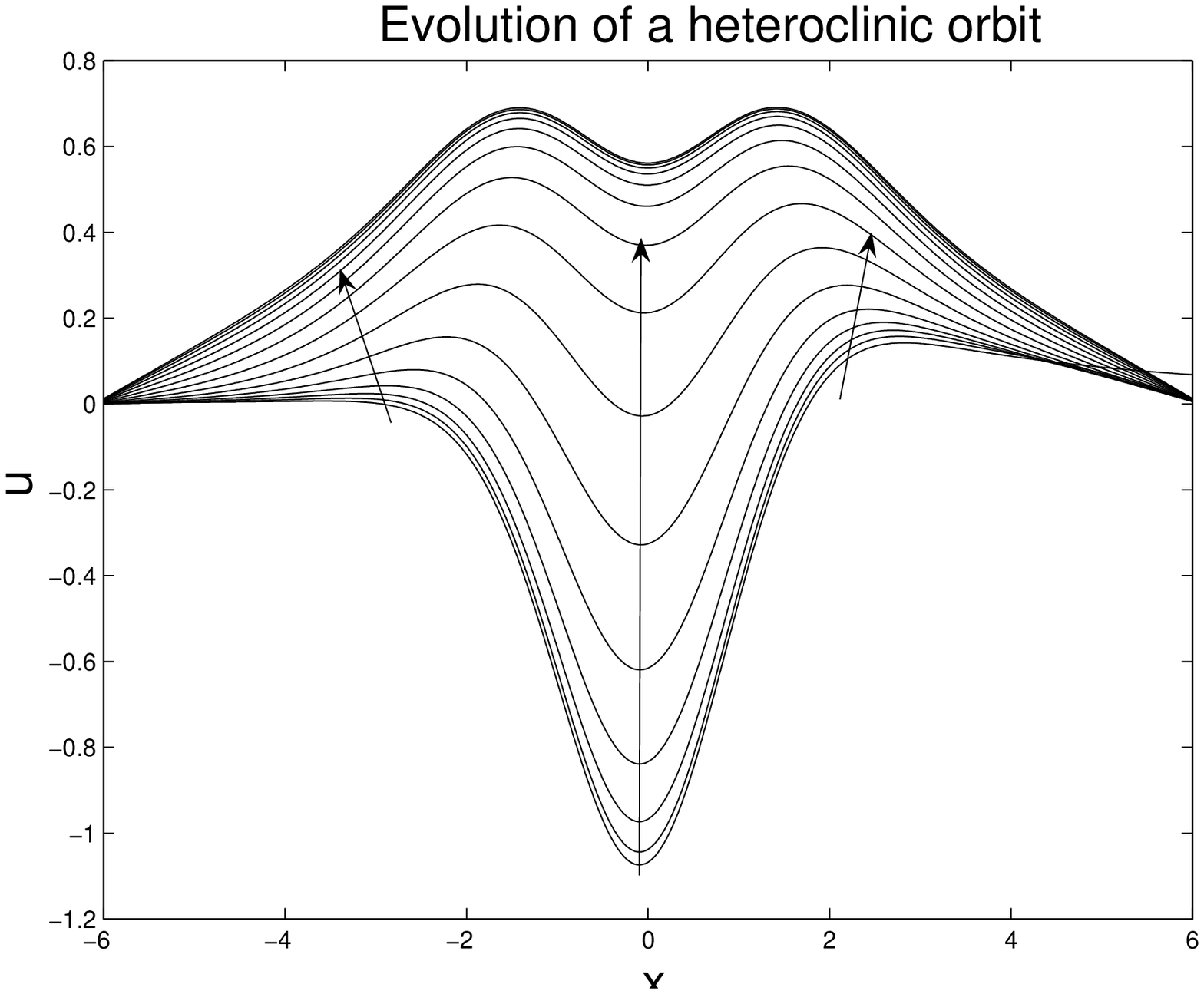}
\includegraphics[height=2in]{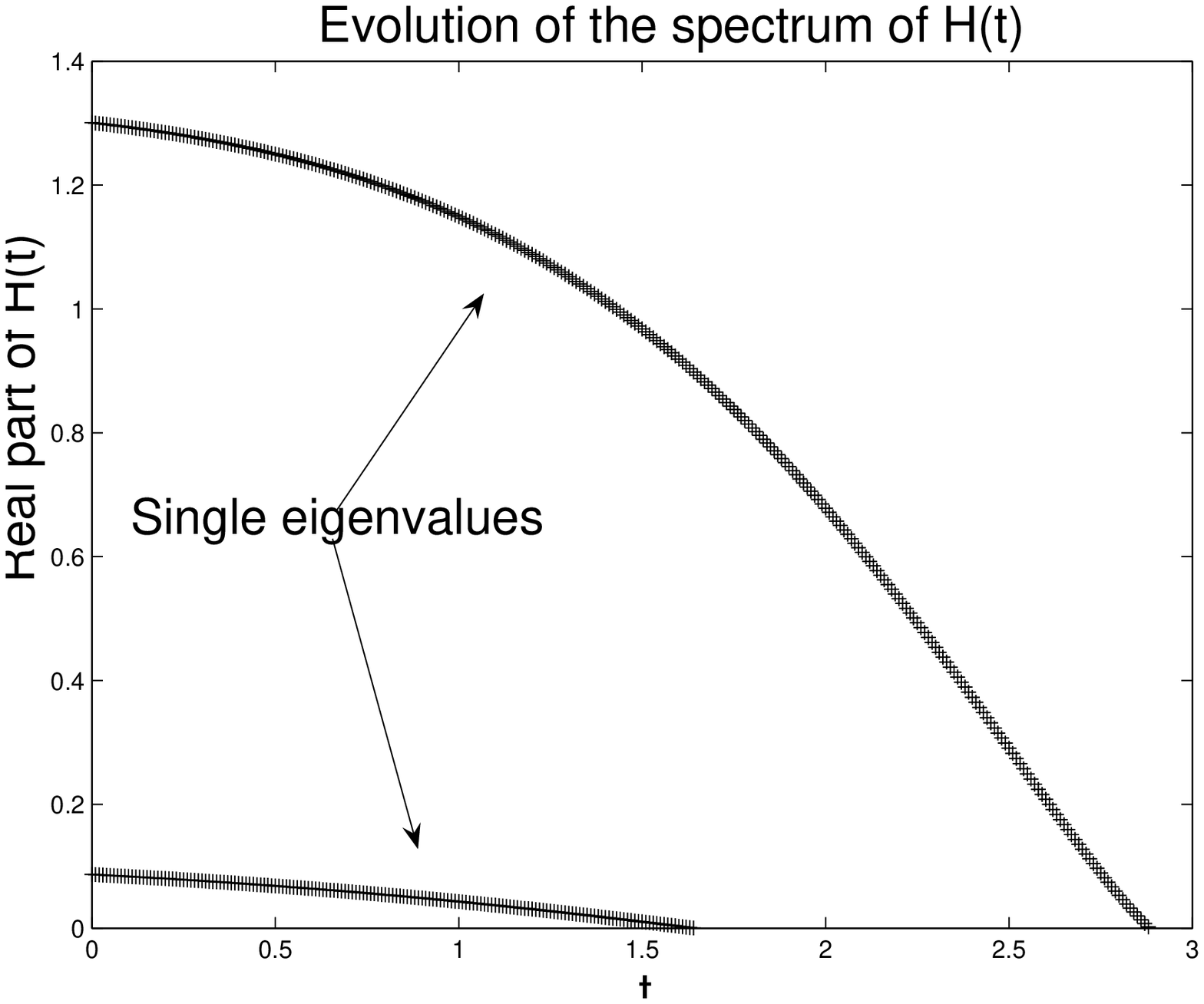}
\end{center}
\caption{A typical heteroclinic orbit to the left of boundary A, with
  the spectrum of $H(t)$ as a function of $t$.}
\label{solb_fig}
\end{figure}

Also of interest is the structure of the flow in the unstable
manifold of the ``fork arms'' which occur at $c=0.0740$, as they
approach the pitchfork bifurcation at $c=0.0501$.  Figure \ref{flow_fig} shows a
schematic of the flow based on numerical evidence.  Of particular
interest is the behavior near the boundary marked A.  Solutions to the
right of the boundary are not eternal solutions -- they fail to exist
for all $t$.  Solutions to the left of A are heteroclinic orbits
connecting the equilibrium with an unstable manifold of dimension 2 to the
equilibrium with an unstable manifold of dimension zero.  A typical
such solution is shown in Figure \ref{solb_fig}.

\begin{figure}
\begin{center}
\includegraphics[width=4in]{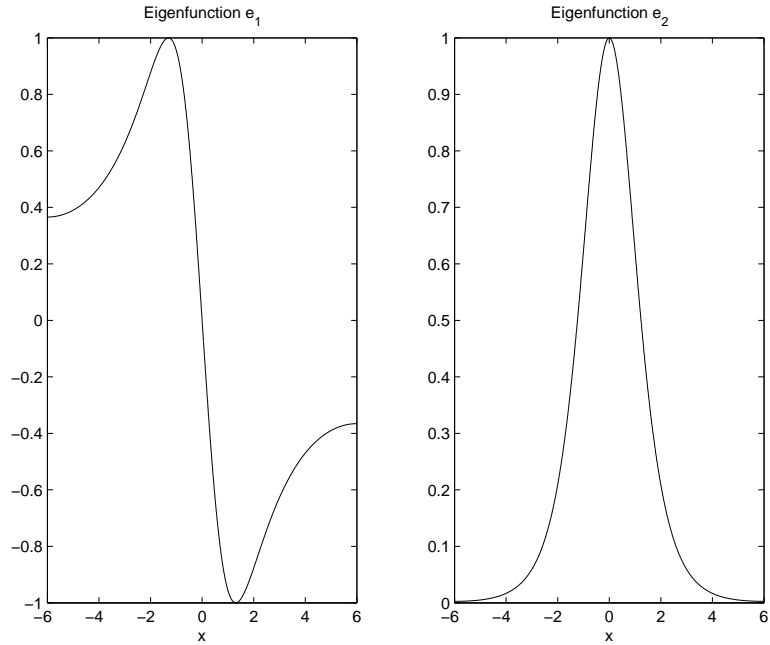}
\end{center}
\caption{Eigenfunctions describing unstable directions at $f_1$}
\label{eigenfuncs_fig}
\end{figure}

To examine solutions near the boundary A, we center our attention on
the case $c=0$, which has two equilibria, one of which (call it $f_1$)
has a 2-dimensional unstable manifold.  (This corresponds to the right
pane of Figure \ref{flow_fig}.)  If we linearize about $f_1$, the
operator $H=\frac{\partial^2}{\partial
x^2}-2f_1:C^{0,\alpha}(\mathbb{R})\to C^{0,\alpha}(\mathbb{R})$ has a
pair of simple eigenvalues, as is easily seen in the right pane of
Figure \ref{solb_fig} at $t=0$.  One of these eigenvalues is smaller,
to which is associated the eigenfunction $e_1$ in Figure
\ref{eigenfuncs_fig}.  The eigenfunction $e_2$ is associated to the
larger eigenvalue.  In Figure \ref{flow_fig}, $e_1$ corresponds to the
horizontal direction, and $e_2$ corresponds to the vertical direction.
From the proof of Lemma \ref{eq_findim}, it is clear that
$\{e_1,e_2\}$ spans the tangent space of the unstable manifold at
$f_1$.  Therefore, we specify initial conditions $u_{A,\theta}(x)$ for
a numerical solver using
\begin{equation}
\label{ic_eq}
u_{A,\theta}(x)=f_1(x) + A \left(e_1(x)\cos \theta + e_2(x) \sin\theta\right).
\end{equation}

\begin{figure}
\begin{center}
\includegraphics[width=5.75in]{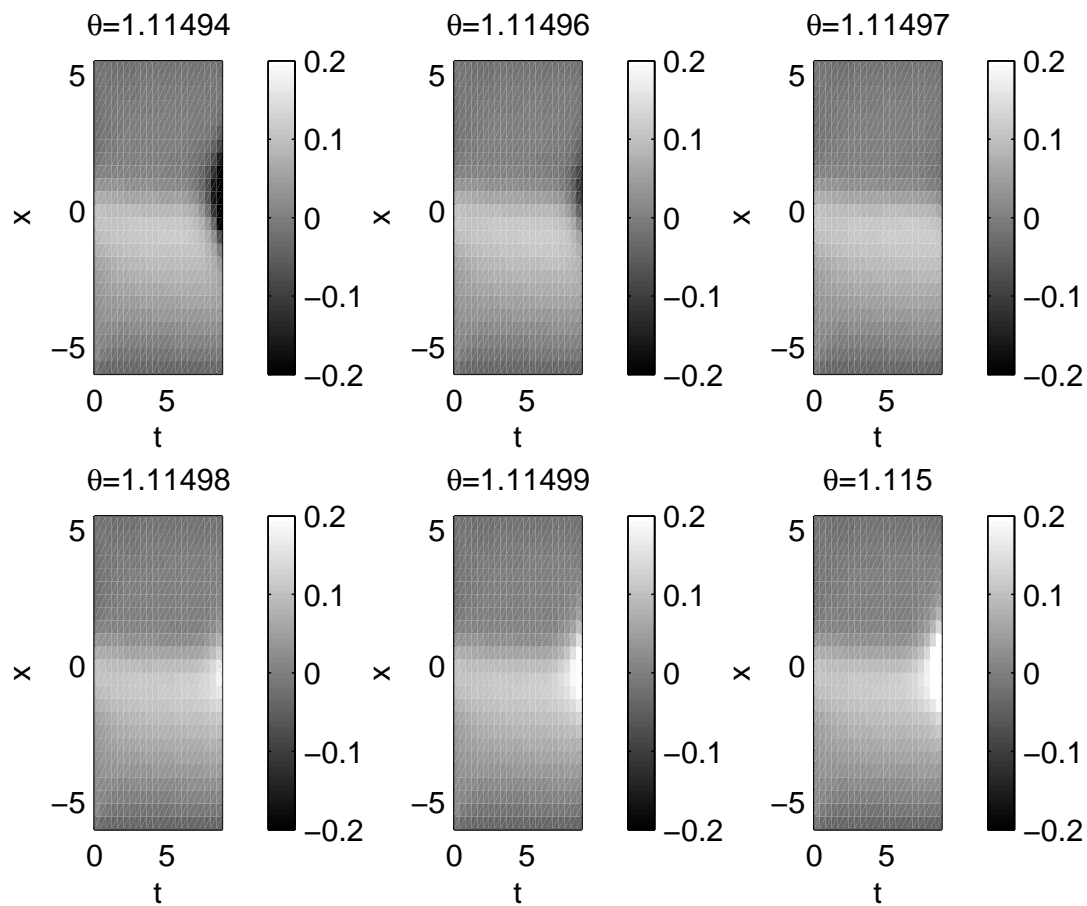}
\end{center}
\caption{Difference between equilibrium $f_1$ and the numerical
  solution started at $u_{A,\theta}$, where black indicates a value of
  -0.2, and white indicates 0.2.  The horizontal axis represents $t$,
  and the vertical axis represents $x$.  $A=0.1$ in all figures.
  Starting from the upper left, $\theta=1.11494, 1.11496, 1.11497,
  1.11498, 1.11499, 1.115.$ }
\label{perturb_fig}
\end{figure}

(Taking $A$ small allows us to approximate solutions which tend to
$f_1$ in backwards time.)  Since the perturbations along $e_1,e_2$ are
quite small, and indeed the eigenvalue associated to $e_1$ is much
smaller than that associated to $e_2$, examining the numerical results
of evolving $u_{A,\theta}$ is quite difficult.  The behavior along the
boundary occurs at a much smaller scale than $f_1$, yet is crucial in
determining the long-time behavior of the solution.  To remedy this,
the boundary behavior is better emphasized by plotting
$u_{A,\theta}(t,x) - f_1(x)$ instead.  Figure \ref{perturb_fig} shows
the results of evolving initial conditions \eqref{ic_eq} for $A=0.1$
and various values of $\theta$.

Solutions in Figure \ref{perturb_fig} show a similar kind of behavior
as in the case of the frontier of $f_0$.  There is a traveling front,
which moves very slowly in the negative $x$-direction.  However, the
behavior is quite a bit more delicate.  The determining factor in
locating the frontier of $f_0$ is the perturbation in a direction
roughly like $e_2$, which has a large eigenvalue.  On the other hand,
for $f_1$, Figure \ref{flow_fig} indicates that such a direction is
not parallel to the boundary of the connecting manifold.  (The
boundary direction is some linear combination of $e_1$ and $e_2$, with
a numerical value for the angle $\theta$ being roughly 1.114975
radians.) The eigenvalue associated to $e_1$ is roughly ten times
smaller, and therefore perturbations in that direction are much more
sensitive.  Additionally, the action of the flow is therefore
primarily in the direction of $e_1$, which tends to mask effects in
other directions.  For this reason, it was visually necessary to
postprocess the numerical solutions by subtracting $f_1$ from them.
Otherwise the presence of the traveling front was unclear.

\chapter{Conjectures and future work}
\label{conj_ch}
\section{Conjectures about the present problem}

\subsection{Analytical conjectures}

It seems that under reasonable conditions on the coefficients of
\eqref{pde}, all eternal solutions ought to be heteroclinic orbits.
An easy calculation with the formula \eqref{action_eqn} for the action
$A(u(t))$ shows that if $|A(u(t))|$ blows up, then one of the
following is true:
\begin{enumerate}
\item $\|u(t)\|_1 \to \infty$,
\item $\|u(t)\|_\infty \to \infty$, or
\item $\|Du(u)\|_2 \to \infty$.
\end{enumerate}
Essentially, eternal solutions which are not heteroclinic orbits are
big in some sense.  Of course, traveling fronts satisfy the first
condition.  On the other hand, the Harnack inequality seems to imply
that eternal solutions do not blow up in the $\infty$-norm as $t\to
-\infty$.  More intriguingly, \cite{Zhang_2000} and
\cite{Souplet_2002} show that under certain conditions on the
coefficients of \eqref{pde}, global solutions to the forward Cauchy
problem have a universal bound on their $\infty$-norm.  However, these
results are obtained under the hypothesis that the solution $u$ is
strictly negative, a condition that is essential to their analysis.
Relaxing this condition leads to currently open problems.

\begin{conj}
Suppose all of the coefficients $a_i$ in \eqref{pde} decay
sufficiently fast as $|x|\to\infty$.  Then all eternal solutions are
bounded in the $\infty$-norm by a universal bound, which depends only
on the $a_i$.
\end{conj}

More ambitious is the following (which involves proving a universal
bound for the 1-norm as well):

\begin{conj}
Suppose all of the coefficients $a_i$ in \eqref{pde} decay
sufficiently fast as $|x|\to\infty$.  Then all eternal solutions are
heteroclinic orbits.
\end{conj}

Related to both of these conjectures is the conjecture that under
suitable decay conditions on the $a_i$, there exist only finitely many
equilibria (Conjecture \ref{finitely_many_equilib_conj}).

\subsection{Conjectures related to the topology of the space
  of heteroclinic orbits}

Much of what remains to be understood about the space of heteroclinic
orbits of \eqref{pde} and \eqref{limited_pde} involves a more precise
understanding of the gluing maps between the cells in its cell complex
structure.  The eventual goal is to construct a homology theory,
called a Floer homology, for the space of heteroclinic orbits.  This
would allow the space of heteroclinic orbits to be decomposed as a
complex of connecting manifolds (without boundary) and boundary maps
which associate higher dimensional manifolds to lower dimensional
ones.  From the outset, degeneracy in the sense of Morse provides the
biggest obstacle to this kind of theory.  In particular, {\it
  nondegeneracy} allows one to show that generically, connecting
manifolds can only have boundaries of one dimension lower.  However,
in the example of the previous chapter, namely that of
\eqref{limited_pde} with $\phi=(x^2-c)e^{-x^2/2}$, such a statement is
still true.  Perhaps it is possible that one can find conditions for
connecting manifolds to have codimension-1 boundaries, even in the
face of degeneracy in the equilibria.  Or put another way,

\begin{conj}
\label{degeneracy_is_ok}
When the flow of \eqref{pde} is restricted to the space of
heteroclinic orbits $\mathcal{H}$, all of the equilibria become nondegenerate
critical points in the sense of Morse.
\end{conj}

Another obstacle is that there needs to be some kind of compactness
result for the space of heteroclinic orbits modulo time translation
(or perhaps modulo action of the flow).  In Floer's case, he was able
to employ Gromov's compactness results for pseudoholomorphic curves.
This leads to his ``no bubbling theorem''.  However, no such result is
known in our case.  The closest available results are those of
\cite{Zhang_2000} and \cite{Souplet_2002}, which only hold for {\it
  positive} solutions to \eqref{pde} that have been centered on an
equilibrium.  It is easy to show that if there exists a positive
equilibrium for \eqref{limited_pde}, then their results suffice to show
compactness, but the situation of general sign is currently an open problem.

To summarize, we have the following conjectures:

\begin{conj}
  The space of heteroclinic orbits of \eqref{pde} modulo time
  translation is compact in $Y_\lambda$.
\end{conj}

\begin{conj}
There is a generic subset (a Baire subset) of choices for the
coefficients $a_i$ in \eqref{pde} so that if $u$ is a heteroclinic
orbit, all of the eigenvalues of $H(t)$ are simple.  
\end{conj}

We then define $C_k(R)$ to be the free $R$-module generated by the
$k$-dimensional connecting manifolds of \eqref{pde}.  Probably it is
best to think of $R=\mathbb{Z}/2$, at least to fix ideas.

\begin{conj}
For a generic subset of coefficients $a_i$ in \eqref{pde}, there is a
collection of maps $\partial_k: C_{k} \to C_{k-1}$ such that
\begin{itemize}
\item $\partial_k$ is an $R$-module homomorphism, for each $k \ge 0$,
\item $\partial_0 = 0$,
\item $\partial_{k-1} \circ \partial_k = 0$ for each $k \ge 0$, and
\item $v = \partial_k(u)$ if and only if roughly speaking $v$ is a sum
  of boundary elements of $u$, obtained by a deformation retraction of
  a neighborhood of $v$ in $u$ onto $v$.  In Floer's case, as is
  likely in ours, this is a rather involved construction called the
  ``gluing theorem'' described in \cite{Floer_spheres}.
\end{itemize}
\end{conj}

This would turn $(C_*(R),\partial_*)$ into what is most reasonably
called a ``Floer complex.''  One can then define the Floer homology
modules, $F_k= \text{ker }\partial_k / \text{im }\partial_{k+1}$ and
formulate the following (reasonable) conjecture:

\begin{conj}
Let $E$ be the space of heteroclinic orbits of \eqref{pde}.  For a
generic subset of $a_i$ in \eqref{pde}, $F_k \cong H_k(E;R)$.  That is,
the Floer complex computes the homology of the space of heteroclinic orbits.
\end{conj}

\section{Future work on related problems}

\subsection{Higher spatial dimensions, with decay conditions enforced}

Of course, the most obvious dependence on 1-dimensional space is the
equilibrium analysis of Chapter \ref{nonauto_ch}.  The analogous
nonlinear elliptic problem \eqref{eq_pde} is not well understood.
Indeed, very little is known about \eqref{eq_pde} at all, especially
if the solutions are allowed to be of general sign.

One thing is likely: the spatial decay of heteroclinic orbits is much
slower -- not in $L^1(\mathbb{R}^n)$.  Worse, Sturm-Liouville theory is no
longer available to control the eigenvalues of $H=\Delta-2u$.
Therefore, there might be infinitely many positive eigenvalues of the
operator $H$, which accumulate at zero.  As a result, equilibria might
have infinite dimensional unstable manifolds.  The analysis of the
structure of the connecting manifolds will therefore not work, though
there should be a filtration structure based on Lyapunov exponent for
the unstable manifolds, which should allow for an infinite dimensional
cell complex with finite dimensional cells.  Additionally, the slower
spatial decay will disrupt the finite energy classification scheme in
Chapter \ref{classify_ch}.

\subsection{Relaxation of decay conditions on the coefficients}

If we no longer require that the coefficients $a_i$ decay to zero as
$|x|\to \infty$, then \eqref{pde} can support traveling wave
solutions.  Indeed, there can be extremely complicated and delicate
traveling wave structures if the spatial dimension is also greater
than 1.  This will remove the uniform convergence to equilibria, of
course.  Also, likely is that what will be found is that the space of
heteroclinic orbits is an infinite-dimensional cell complex, perhaps
were the ``cells'' are Banach manifolds.  The resulting dynamics can
therefore be expected to become extremely complicated.

\appendix
\chapter{Spectrum of Schr\"odinger operators}
\label{schrod_ch}
\section{Introduction}

This appendix recounts a few standard facts about the structure
of the spectrum of the Laplacian and Schr\"odinger operators $\Delta,H: C^2(\mathbb{R})
\to C^0(\mathbb{R})$ respectively.  Nothing in this appendix is
original, but it is useful to have the facts and the requisite
calculations available for reference.

\subsection{Spectrum of the Laplacian operator}
\begin{prop}
The spectrum of $\Delta: C^2(\mathbb{R})\to
C^0(\mathbb{R})$ is the closure of the negative real axis.
\begin{proof}
The spectrum of $\Delta$ contains all $\lambda \in \mathbb{C}$ for which $(\Delta -
\lambda)$ is not injective.  In other words, it contains the solutions
to the equation
\begin{equation*}
u''-\lambda u = 0.
\end{equation*}
An elementary calculation yields that $u=c_1 e^{\sqrt{\lambda} x}+c_2
e^{-\sqrt{\lambda} x}$.  If $\lambda$ is real and nonpositive, there are
nontrivial bounded solutions (which are oscillatory or constant).
Otherwise the nontrivial $C^2$ solutions are unbounded.  Hence
the spectrum must contain the closed negative real axis.

Next, we show $(\Delta - \lambda)^{-1}$ exists and is bounded away from
  the closed negative real axis.
To show that $(\Delta - \lambda)^{-}$ exists, we find an inversion formula, which is valid when $\lambda$ does not
lie on the closed negative real axis.  To this end one can solve
\begin{equation*}
u''-\lambda u = f
\end{equation*}
using a slightly modified version of
Calculation \ref{invert_calc}, 
\begin{equation*}
u(x)=\frac{1}{2\sqrt{\lambda}}\int f(y) \left\{\begin{matrix}
(-1) e^{-\sqrt{\lambda} (y-x)} & \text{ if }
    \Re(\sqrt{\lambda})(y-x)>0\\
e^{-\sqrt{\lambda} (y-x)} & \text{ otherwise}\\
\end{matrix}\right\} dy. 
\end{equation*}

This inversion formula defines a bounded inverse for $(\Delta
-\lambda)$ when $\lambda$ is not on the closed negative real axis.  To
see this, simply observe that
\begin{equation*}
\frac{1}{2|\sqrt{\lambda}|} \int e^{-|\Re(\sqrt{\lambda})||y-x|} dy
\end{equation*}
is independent of $x$.  This fact implies that
$\|u\|_{C^2} \le K_\lambda \|f\|_{C^0}$ for some finite $K$ by
differentiation under the integral.  Therefore the complement of the
closed negative real axis is in the complement of the spectrum of
$\Delta$.
\end{proof}
\end{prop}

\section{Spectrum of Schr\"odinger operators}
The previous section can be generalized to the case of Schr\"odinger
operators $H=(\Delta -V): C^2(\mathbb{R})\to C^0(\mathbb{R})$ to
obtain a few results of interest.  Assume that $V$ is a smooth
function which satisfies $\lim_{|x|\to\infty} V(x) = A$.

\begin{prop}
\label{spectral_control}
The spectrum of $H$ contains the portion of the real axis less than or
equal to $-A$.  All of the eigenvalues of $H$ are contained in the
portion of the real axis less than or equal to the supremum of $V$.
\begin{proof}
Of course, the eigenvalues $\lambda$ are those where there are
nontrivial solutions to the equation
\begin{equation}
\label{schrod_inj_eq}
u''(x)-(V(x)+\lambda)u(x)=0.
\end{equation}
Recast \eqref{schrod_inj_eq} as a first-order system, namely
\begin{equation}
\label{schrod_fo_eq}
\frac{d}{dx}\begin{pmatrix}u\\u'\end{pmatrix}=\begin{pmatrix}0& 1\\
  V(x)+\lambda & 0\end{pmatrix}\begin{pmatrix}u\\u'\end{pmatrix}=T(x)\begin{pmatrix}u\\u'\end{pmatrix}.
\end{equation}
Observe that the eigenvalues of $T(x)$ are purely imaginary if
$\lambda$ is real and less than $-V(x)$.  Thus, the flow
restricted to a plane of constant $x$ consists of periodic orbits if
and only if $\lambda \le - V(x)$ (ignoring the origin, of course).  

Now suppose that $\lambda$ is real and $\lambda < -A$.  Then there
exists an $R>0$ such that for all $|x|>R$, $|f(x)-A| < \frac{1}{2}
|\lambda + A|$.  So on $\mathbb{R} - [-R,R]$, solutions to
\eqref{schrod_inj_eq} will all tend to limiting cycles.  On $[-R,R]$,
solutions grow exponentially fast, at a rate of no more than
$\sqrt{\|v\|_\infty+|\lambda|}$.  Thus there exist nontrivial bounded
solutions to \eqref{schrod_inj_eq}, which are obviously in $C^2$.
Hence the spectrum of $H$ contains $(-\infty,-A]$.

On the other hand, if $\lambda$ is real and $A \le \text{sup }V <
\lambda$, the origin is always a saddle point for the first order
system \eqref{schrod_fo_eq}.  Thus all solutions to
\eqref{schrod_inj_eq} are unbounded.  Likewise, if $\lambda \in
\mathbb{C} - \mathbb{R}$, then the system \eqref{schrod_fo_eq} is a spiral, so 
$(H-\lambda)$ is injective for $\lambda \notin (-\infty,\text{sup
}V]$.
\end{proof}
\end{prop}

It is best to treat the portion of the spectrum lying between $-A$ and $\text{sup }V$ 
using Sturm-Liouville theory.  Indeed, \eqref{schrod_inj_eq} on
$[0,\infty)$ with boundary conditions
\begin{equation*}
u(0)=a,\;\lim_{x\to\infty} u(x)=0
\end{equation*}
is a classic Sturm-Liouville problem.  It is known that the
eigenvalues of this problem are discrete and accumulate only at zero.

\begin{prop}
\label{finitely_many_eigenvalues}
Suppose $V \in C^0(\mathbb{R})$ is positive outside a compact
interval.  Then the operator $H$ has finitely many
eigenvalues greater than $-A$.
\begin{proof}
It suffices to show that any solution $u$ to
\begin{equation*}
u''-(V+A)u = 0
\end{equation*}
has finitely many zeros.  Notice that (assuming by hypothesis that
$\lambda>-A$) 
\begin{equation*}
-V-\lambda < - V-A
\end{equation*}
so the Sturm-Liouville comparision theorem states that any
eigenfunction $v$ of $H$ with eigenvalue $\lambda$ has strictly fewer zeros
than solutions to the $\lambda=-A$ case.  Also, the zeros of solutions
which are bounded for half intervals $(-\infty,x_0)$ are monotonic in
$\lambda$.

By hypothesis, $V$ is positive on $\mathbb{R} - [-R,R]$ for some $R>0$.
By comparision with the case of $V \equiv 0$ on $\mathbb{R} - [-R,R]$,
zeros of $u$ only occur on $[-R,R]$.  By comparison with
$\|V\|_\infty$, the number of zeros is proportional to
$R\sqrt{\|V\|_\infty}$, and therefore finite.
\end{proof}
\end{prop}

It is Proposition \ref{finitely_many_eigenvalues} that ensures that
the finite dimensionality results of Chapter \ref{unstable_ch} hold.
In the case of higher spatial dimensions, Sturm-Liouville theory does
not apply (at least if there is no assumed symmetries in the
equilibria).  It is therefore possible that there is no finite
dimensionality for the space of heteroclines for higher spatial
dimensions.

It will be technically important in Chapter \ref{unstable_ch} that
$(H+s)=(\Delta - (V + s))$ is sectorial about any real $s$ not in the
spectrum of $H$.  It is then useful to consider the densely
defined operator $(H+s):C^{0,\alpha}(\mathbb{R}) \to
C^{0,\alpha}(\mathbb{R})$ intead of $C^2 \to C^0$, where $0<\alpha \le 1$.

\begin{prop}
\label{sectorial_H_prop}
If $s \in \mathbb{R}$ is not in the spectrum of $H=(\Delta - V)$, then $(H+s):C^{0,\alpha}(\mathbb{R}) \to C^{0,\alpha}(\mathbb{R})$ is sectorial for $0<\alpha \le 1$.
\begin{proof}
We have already shown that the structure of the spectrum is favorable
for the sectoriality of $(H+s)=(\Delta - (V+s))$, and the operator is
densely defined if $\alpha > 0$.  (It is {\it not} densely defined if
$\alpha=0$.)  What remains is that the operator norm of the inverse
must decay like $K/|s|$, and that its image must consist of continous
functions.  For the former:

\begin{eqnarray*}
u'' - (V+s) u &=& f\\
s\left ( \frac{1}{s} \Delta - I \right) u &=& V u + f \\
u &=& \frac{1}{s} \left ( \frac{1}{s} \Delta - I \right)^{-1} \left (
V u + f \right)\\
u(x)&=&-\frac{1}{2\sqrt{s}} \int e^{-|x-y|\sqrt{s}} (V(y)u(y)+f(y)) dy,\\
\end{eqnarray*}
where $s$ is not in the spectrum.  Using Calculation \ref{bound_calc}, 
\begin{eqnarray*}
\|u\|_\infty &\le& \frac{1}{|s|} \left \| \left ( \frac{1}{s} \Delta - I
\right )^{-1} (Vu+f) \right \|_\infty \\
&\le& \frac{1}{|s|} \|Vu\|_\infty + \frac{1}{|s|} \|f\|_\infty\\
&\le& \frac{\frac{1}{|s|}}{1-\frac{1}{|s|} \|V\|_\infty} \|f\|_\infty.
\end{eqnarray*}

Since $V$ and $f$ are assumed to decay to zero and the kernel
$e^{-|x-y|\sqrt{s}}$ is in $L^1$, it is immediate that $u$ must also
decay to zero.  The following calculation shows that image of $(\Delta
- (V+s))^{-1}$ consists of Lipschitz functions.  Assume $x_1 < x_0$,
so that

\begin{eqnarray*}
|u(x_0)-u(x_1)| &=& \frac{1}{2\sqrt{s}} \left | \int_{-\infty}^{x_0}
 e^{(y-x_0)\sqrt{s}} (Vu+f)dy - \int_{x_0}^\infty e^{(x_0-y)\sqrt{s}}
 (Vu+f)dy \right .\\&& \left . - \int_{-\infty}^{x_1}
 e^{(y-x_1)\sqrt{s}} (Vu+f)dy + \int_{x_1}^\infty e^{(x_1-y)\sqrt{s}}
 (Vu+f)dy \right|\\
&=& \frac{1}{2\sqrt{s}} \left| \int_{-\infty}^{x_1}
 \left(e^{(y-x_0)\sqrt{s}} - e^{(y-x_1)\sqrt{s}}\right)(Vu+f) dy
 \right. \\&&+
 \int_{x_1}^{x_0} e^{(y-x_0)\sqrt{s}}(Vu+f) dy \\&&+ \int_{x_0}^{\infty}
 \left(e^{(x_1-y)\sqrt{s}} - e^{(x_0-y)\sqrt{s}}\right)(Vu+f) dy \\&&+ \left.
 \int_{x_1}^{x_0} e^{(x_1-y)\sqrt{s}}(Vu+f) dy \right|\\
&\le& \frac{1}{2\sqrt{s}} \|Vu+f\|_\infty \left( \int_{-\infty}^{x_1}
 \left | e^{(y-x_0)\sqrt{s}} - e^{(y-x_1)\sqrt{s}}\right| dy \right.\\&&+
 \left. \int_{x_0}^{\infty}
 \left|e^{(x_1-y)\sqrt{s}} - e^{(x_0-y)\sqrt{s}}\right| dy + \frac{2}{\sqrt{s}}(1-e^{(x_1-x_0)\sqrt{s}})
\right)\\
&\le& \frac{1}{2\sqrt{s}} \|Vu+f\|_\infty \left( \frac{1}{\sqrt{s}}e^{x_1\sqrt{s}}
 \left | e^{-x_0\sqrt{s}} - e^{-x_1\sqrt{s}}\right| \right.\\&&
 \left. +
 \frac{1}{\sqrt{s}}^{-x_0\sqrt{s}}
 \left|e^{x_1\sqrt{s}} - e^{x_0\sqrt{s}}\right| + \frac{2}{\sqrt{s}}(1-e^{(x_1-x_0)\sqrt{s}})
\right)\\
&\le& \frac{2}{|s|} \|Vu+f\|_\infty  \left | e^{(x_1-x_0)\sqrt{s}} - 1\right|\\
&&\to 0 \text{ as } |x_0 - x_1| \to 0.\\
\end{eqnarray*}
But since we've chosen $x_1 < x_0$, the above calculation proves that
$u$ is Lipschitz.  Looking at the Lipschitz constant and the bound on
$u$, it is immediate that the $C^{0,1}$-operator norm of $(H+s)^{-1}$
decays like $1/|s|$.  Thus $(H+s)$ is sectorial.
\end{proof}
\end{prop}

It is important to remark that the above proof shows that in fact
$(\Delta - (V+s))^{-1}$ is a bounded operator $C^0(\mathbb{R}) \to
C^2(\mathbb{R})$.  However, the norm of the operator does not decay as
$\Im(s) \to \infty$ as required for a sectorial operator.  In
particular, 
\begin{equation*}
\left | \frac{du}{dx} \right| \le \frac{K_1 \|Vu+f\|_\infty}{\sqrt{s}},
\end{equation*}
and
\begin{equation*}
\left | \frac{d^2u}{dx^2} \right| \le K_2 \|Vu+f\|_\infty,
\end{equation*}
for $K_1,K_2$ independent of $s$ and $f$.  Examples can be constructed
to show that these bounds are tight.

\bibliography{b_exam_bib}

\end{document}